\documentclass[a4paper]{amsart}

\usepackage[latin1]{inputenc}
\usepackage[T1]{fontenc}
\usepackage{palatino}

\usepackage[english]{babel}
\usepackage{amsthm, amssymb, amsmath, amsfonts}
\usepackage{hhline}
\usepackage{float}
\usepackage{colortbl}
\usepackage{enumitem}

\usepackage{tikz}		
\usetikzlibrary{arrows, fit, calc, shapes.misc, shapes.geometric, decorations.markings, backgrounds, matrix}

\usepackage{pgfplots}
\usepackage{pgf}
\usepackage{hyperref}
\usepackage[capitalise]{cleveref}
\usepackage[fixamsmath,disallowspaces]{mathtools}
\usepackage{cite,placeins}
\usepackage[arrow, matrix, curve, color]{xy}
\usepackage{mathdots}

\newtheorem{thm}{Theorem}[section]
\newtheorem{lem}[thm]{Lemma} 
\newtheorem{cor}[thm]{Corollary}
\newtheorem{prop}[thm]{Proposition}
\newtheorem{conj}[thm]{Conjecture}

\newtheorem{thm*}{Theorem}
\newtheorem{lem*}[thm*]{Lemma} 
\newtheorem{cor*}[thm*]{Corollary}
\newtheorem{prop*}[thm*]{Proposition}
\newtheorem{conj*}[thm*]{Conjecture}

\theoremstyle{definition}
\newtheorem{defn}{Definition}[section]
\theoremstyle{remark}
\newtheorem{rem}[defn]{Remark}

\newenvironment{prf}{\begin{proof}[Proof]}{\end{proof}}

\newcommand{\Bl}[1]{{\mathbb{#1}}}

\newcommand{\DZ}{\Bl{Z}}
\newcommand{\DN}{\Bl{N}}

\newcommand{\DR}{\Bl{R}}
\newcommand{\DC}{\Bl{C}}

\newcommand{\DV}{{\Bl{V}}}

\newcommand{\DL}{\Bl{L}}

\newcommand{\op}[1]{{\operatorname{#1}}}

\newcommand{\uH}{\underline{H}}

\newcommand{\oBx}{\overline{B_x}}
\newcommand{\oB}{\overline B}

\newcommand{\sbim}{\mathcal S\mathcal B}
\newcommand{\smod}{\mathcal S}
\newcommand{\bsbim}{{\mathcal B\mathcal S}}
\newcommand{\hecke}{\textbf H }
\newcommand{\gp}[1]{\mathcal K_0 \left(#1 \right)}
\newcommand{\gps}[1]{\mathcal K_0^s \left(#1 \right)}

\mathchardef\mhyphen="2D
\newcommand{\mhy}{\mhyphen}
\newcommand{\Pm}{\textbf{P}_m}
\newcommand{\LP}{\DZ[v^{\pm 1}]}
\newcommand{\scol}{{\textcolor{red}{s}}}
\newcommand{\tcol}{{\textcolor{blue}{t}}}

\newcommand{\un}[2]{{}_{#1}\underline{ #2}}
\newcommand{\unr}[2]{\underline{ #2}_{#1}}

\tikzset{%
  highlight/.style={rectangle,rounded corners,fill=red!15,draw=red,
    fill opacity=0.5,thick},
  bendBelow/.style={bend left=70, looseness=2},
  bendAbove/.style={bend right=70, looseness=2},
  object/.style={circle, fill, inner sep=1.5pt, outer sep=0mm},
  labeling/.style={outer sep=0mm, inner sep=0mm},
  1morph/.style={->, shorten >= 0.5pt, >=stealth'},
  2morph/.style={-implies,double,double equal sign distance,
                 shorten >=2pt, shorten <=3pt},
  spot/.style={color=black, thin, dashed},
  sline/.style={color=red, line width=1.5pt},
  tline/.style={color=blue, line width=1.5pt},
  uline/.style={color=green, line width=1.5pt},
  randline/.style={color=pink, line width=1.5pt},
  sdot/.style={color=red, thin, fill},
  tdot/.style={color=blue, thin, fill},
  udot/.style={color=green, thin, fill},
  randdot/.style={color=pink, thin, fill},
  bbline/.style={color=black, line width = 1.5pt},
  bbdot/.style={color=black, thin, fill},
  bline/.style={color=black, line width=0.8pt}
}

\author{Marc Sauerwein}
\title{Koszul Duality and Soergel Bimodules for Dihedral Groups}
\address{Department of Mathematics, University of Bonn, 53115 Bonn, Germany}
\email{sauerwein@math.uni-bonn.de}
\begin{document}

\pgfdeclarelayer{background}
\pgfdeclarelayer{foreground}
\pgfsetlayers{background,main,foreground}

\setcounter{tocdepth}{1}

\setcounter{page}{1}\pagenumbering{arabic}
\begin{abstract}
Every Coxeter group $(W,S)$ gives rise to an associated Hecke algebra $\hecke_{(W,S)}$ which can be categorified by the additive monoidal category of Soergel bimodules $\sbim$. Under this isomorphism the Kazhdan-Lusztig basis $\{\underline{H}_x\}_{x\in W}$ corresponds to certain indecomposable Soergel bimodules $\{B_x\}_{x\in W}$ (up to shift). In this thesis we study the structure of the endomorphism algebra (of maps of all degrees) $\mathcal A:= \op{End}^\bullet_{\sbim}\left(\bigoplus_{x\in W} B_x \right)\otimes_R\DR$. Via category $\mathcal{O}$ it has been proven for all Weyl groups (see \cref{thm:sd}) that $\mathcal A$ is a self-dual Koszul algebra. We extend this result to all dihedral groups by purely algebraic methods using representation theory of quivers and Soergel calculus.
\end{abstract}

\maketitle
\tableofcontents

\section{Introduction}
\noindent
To any Coxeter system $(W,S)$ one has an associated Hecke algebra $\hecke_{(W,S)}$. The Hecke algebra may be categorified by Soergel bimodules $\sbim$, an additive monoidal category of bimodules over a polynomial ring. The indecomposable bimodules $\{B_x\}$ in $\sbim$ are (up to grading shift) parametrised by the group $W$. The main object in this paper is the endomorphism algebra (consisting of maps of all degrees) of $\textbf B:= \bigoplus_{x\in W} B_x$ where the right action of polynomials of positive degree is trivialised:
\[\mathcal A := \op{End}^\bullet_{\sbim}\left( \textbf B \right)\otimes_R \DR. \] We prove the following result via purely algebraic methods:

\begin{thm*}\label{thm:goal}
For a dihedral group $(W,S)$ the $\DR$-algebra $\mathcal A$ is a self-dual Koszul algebra.
\end{thm*}

\subsection{Motivation.}
Let $\mathfrak g$ be a complex semisimple Lie algebra. It turns out that the category $\mathfrak g\mhy \op{Mod}$ of all $\mathfrak g$-modules is far too large to be understood algebraically. Seminal for the further study of the representation theory of $\mathfrak g$ was the introduction of category $\mathcal O$ by Bernstein, Gelfand and Gelfand (see \cite{BGG}).
\noindent Fix a Borel $\mathfrak b$, a Cartan $\mathfrak h$ in $\mathfrak g$ and define $\mathcal O := \mathcal O(\mathfrak g, \mathfrak b,\mathfrak h)$ to be the full subcategory of $\mathfrak g\mhy\op{Mod}$ whose elements $M$ are finitely generated over $\mathfrak g$, $\mathfrak h$-semisimple and locally $\mathfrak b$-finite. In particular, all finite dimensional modules and Verma (= standard) modules $\Delta(\lambda)$ ($\lambda \in \mathfrak{h}^*$) lie in $\mathcal O$. 
This restriction made it easier to handle the category and led to beautiful new results such as BGG reciprocity \cite[Prop. 1]{BGG} and the Kazhdan-Lusztig conjectures \cite[Conj. 1.5.]{KL:C}. Within the principal block $\mathcal O_0\subset \mathcal O$ let $L:=\bigoplus_{w\in W}L(w\cdot 0)$ be the direct sum of the simple modules and $P:=\bigoplus_{w\in W}P(w\cdot 0)$ the direct sum of their corresponding (indecomposable) projective covers, i.e. $P$ is a projective generator. We have the following result due to Soergel in \cite{S:CO}:

\begin{thm*}[Koszul self-duality for the principal block $\mathcal O_0$]\label{thm:sd}
There exists an isomorphism of finite dimensional $\DC$-algebras 
\[A:=\op{End}_{\mathcal O_0}(P)\cong \op{Ext}^\bullet_{\mathcal O_0}(L,L),\]
where the right hand side is a ring via the cup product. Furthermore, $\op{Ext}^\bullet_{\mathcal O_0}(L,L)$ is a Koszul algebra.
\end{thm*}

Although $\op{End}_{\mathcal O_0}(P)$ is not obviously graded, it inherits a grading from the naturally graded $\op{Ext}$-algebra. The first glimpse of Koszul duality was discovered earlier when mathematicians were investigating composition series of Verma modules  in category $\mathcal O_0$ and found formulas of the form
\[[\Delta(x\cdot 0):L(y\cdot 0)] = \sum_i \op{Ext}^i (\Delta(w_0x\cdot 0),L(w_0y\cdot 0)).\]
These formulas can be explained by Koszul self-duality on the level of derived categories (see \cite[Theorem 1.2.6.]{BGS}). The existing proofs of Koszul self-duality are difficult and rely heavily on geometric techniques.

\noindent Using $\mathcal O_0 \cong \op{Mod}\mhy A$ one obtains a $\DZ$-graded version of $\mathcal O_0$ as $\mathcal O_0^\DZ := \op{gMod}\mhy \op{A}$. We have the following isomorphisms of $\DZ[v^{\pm 1}]$-modules (see \cite{S:CO} and \cite[Theorem 7.1.]{St:COG}):
\begin{align*}
\gp{\mathcal O_0^\DZ} & \overset{\sim}{\longrightarrow} \gps{\smod}\\
P(x\cdot 0)\langle i \rangle &\longmapsto \DV P(x\cdot 0)\langle i\rangle = \oBx\langle i \rangle,
\end{align*} 
where $\oB := B \otimes_R \DR$ is the Soergel module corresponding to the Soergel bimodule $B$. Let $\smod$ denote the category of graded Soergel modules and $\DV$ Soergel's combinatorial functor. 
Translating \cref{thm:sd} into the setting of Soergel (bi)modules via the above identification frees the result from geometry and yields
\[\mathcal A \overset{\text{Thm. }\ref{thm:sd}}{\cong} \op{Ext}_{\mathcal O_0}^\bullet (L,L) \cong E(\mathcal A),\]
where $E(\mathcal A)$ denotes the Koszul dual of a Koszul algebra.

\subsection{Structure of the paper.} This paper contains two parts.
\begin{enumerate}[label=\textbf{Part} \arabic*:,labelwidth=1in,listparindent=0.5in]
\item In the first four section we provide the necessary background on the Hecke algebra, Soergel bimoldues, Sorgel calculus and Koszul algebras.
\item In the last two sections we realise the endomorphism ring of Soergel bimodules as a path algebra of a quiver and show its Koszul self-duality. 
\end{enumerate}

 \subsection{Acknowledgements}
This paper started as my masther's thesis under the supervision of Geordie Williamson. I would like to thank him for his support and this very interesting project. I would also like to thank Catharina Stroppel for very helpful discussions.

\section{Preliminaries}
\subsection{Basic Definitons}
Let $(W,S)$ be a Coxeter system. Recall that $W$ is equipped with the Bruhat order $\leq$ and the length function $\ell: W\to \DN_0$ which counts the number of simple reflections in a reduced expression. For an arbitrary sequence $\underline{w}=(s_1,s_2,\ldots, s_n) $ in $S$ we denote the product $s_1s_2\cdots s_n$ by $w$, viewed as an element in $W$. For such a sequence $\underline w$ its length is defined by $\ell(\underline w)=n$. Observe that $\ell(\underline w) \geq \ell(w)$ with equality if and only if $\underline w$ is a reduced expression for $w$. By abuse of notation we write $\underline w = s_1s_2\cdots s_n$. It is crucial to distinguish between $w$ and $\underline w$, since the latter denotes a distinct sequence of simple reflections whereas $w$  is their product in $W$. 

Following Soergel's normalisation as in \cite{S:KLP} define the associated Hecke algebra $\hecke = \hecke_{(W,S)}$ as the unital, associative $\LP$-algebra generated by $\{H_s\}_{s\in S}$ with relations
\begin{align}
H_s^2 & = 1 + (v^{-1}-v)H_s, \\
\underbrace{H_sH_tH_s \ldots}_{m_{st} \text{ factors}} &	= \underbrace{H_tH_sH_t \ldots}_{m_{st} \text{ factors}},\label{eqn:braid}
\end{align}
for all $s\neq t \in S$. Here $v$ is just an indeterminant.

Given a reduced expression $\underline{w}=s_1s_2\cdots s_n$ we set $H_{w} := H_{s_1}\cdots H_{s_n}$ which is well-defined by the Lemma of Matsumoto (see \cite{Ma:WG}). The elements $\{H_w\}_{w\in W}$ form the \emph{standard basis} of $\hecke$ as a $\LP$-module. An easy calculation shows

\begin{lem}\label{lem:mult}
Let $w\in W$ and $s\in S$. For a basis element $H_w$ and a generator $H_s$ there is the multiplication rule:
\begin{align}
H_w H_s \enspace & = \enspace  
\begin{cases}
H_{ws} & \text{ if  } ws > w \\
H_{ws} + (v^{-1}-v)H_{w} & \text{ if  } ws< w
\end{cases}.
\end{align}
\end{lem}

Each $H_s$ for $s\in S$ is invertible with inverse $H_s + (v-v^{-1})$ and thus all standard basis elements are units. There is a unique $\DZ$-linear involution $\overline{\cdot}: \hecke \to \hecke$ such that $v \mapsto v^{-1}$ and $H_s \mapsto H_s^{-1}$. This involution is called \emph{duality} and it is easily checked that $H_w \mapsto H_{w^{-1}}^{-1}$ for $w\in W$.

\begin{thm}[\cite{KL:C}]\label{thm:klbasis}
There exists a unique $\LP$-basis $\{\uH_w\}_{w\in W}$ of $\hecke$ consisting of self-dual elements such that
\[\uH_w = H_w + \sum_{x <w}h_{x,w}H_x\]
where $h_{x,w} \in v\DZ[v]$. 
\end{thm}
This basis is called the \emph{Kazhdan-Lusztig basis} and the $h_{x,w}$ are the \emph{Kazhdan-Lusztig polynomials}. Note that the $h_{x,w}$ are not the originally defined Kazhdan-Lusztig polynomials $p_{x,w}$ (see \cite{KL:C}) but there is the following relation (see \cite[Rem. 2.6.]{S:KLP}):
\[h_{x,w}(v) = v^{\ell(w)-\ell(x)}p_{x,w}(v^{-2}).\]

Mimicking the notion of a trace form from linear algebra we mean by a trace on $\hecke$ a $\LP$-linear map $\varepsilon: \hecke \to \LP$ satisfying $\varepsilon(h_1h_2)=\varepsilon(h_2h_1)$ for all $h_1,h_2 \in \hecke$. The \emph{standard trace} on $\hecke$ is defined via $\varepsilon(H_w):= \delta_{w,e}$. 

\begin{lem}\label{lem:stdtr}
For $w$ and $w^\prime \in W$ we have $\varepsilon(H_wH_{w^\prime})=\delta_{w,(w^\prime)^{-1}}$.
\end{lem}

\begin{prf}
Induction over the length of $w^\prime$ combined with the multiplication rule from \cref{lem:mult}.
\end{prf}

For a sequence $\underline w= s_1s_2 \cdots s_n$ we call $\textbf{e}=\textbf{e}_1\textbf{e}_2\cdots \textbf{e}_n$ with $\textbf{e}_i\in \{0,1\}$ a \emph{$01$-(sub)sequence} of $\underline w$ which picks out the subsequence $\underline{w}^{\textbf{e}}:=s_1^{\textbf{e}_1}s_2^{\textbf{e}_2}\cdots s_n^{\textbf{e}_n}$. Given such a $01$-sequence its \emph{Bruhat stroll} is the sequence $e,x_1,\ldots,x_n=\underline{w}^{\textbf e}$ where
\begin{align}\label{eqn:bruhatstroll}
x_i:= s_1^{\textbf{e}_1}s_2^{\textbf{e}_2}\cdots s_i^{\textbf{e}_i}.
\end{align}
This stroll allows us to decorate each step of the $01$-sequence $\textbf{e}$ with either U(p) or D(own) encoding the path in the Bruhat graph. We assign U to the index $i$ if $x_{i-1}s_i > x_{i-1}$ and D if $x_{i-1}s_i < x_i$.

There is a $\LP$-linear anti-involution $\iota$ satisfying $\iota(H_s) = H_s$ for $s\in S$. It is easily checked that $\iota(H_w) = H_{w^{-1}}$ for $w\in W$. Note that $\iota$ and $\overline{\cdot}$ commute, we denote their composition as $\omega$. It follows that $\omega$ is a $\DZ$-linear anti-involution on $\hecke$ satisfying $\omega(v)=v^{-1}$ and $\omega(H_w)=H_w^{-1}$ for $w\in W$. Using the aforementioned standard trace we can define the \emph{standard pairing} $\hecke \times \hecke \to \LP$ by $(h,h^\prime):= \varepsilon(\omega(h)h^\prime)$.

\subsection{Notations (following \cite{E:DC})}
Let $(W,S)$ be a \emph{dihedral group of type $I_2(m)$} $(m\geq 3)$, that is a Coxeter system $(W,\{\scol,\tcol\})$ where $(\scol\tcol)^m=e$. The elements in $S=\{ \scol,\tcol \}$ are called \emph{simple reflections} or \emph{colours}. As before, we denote a sequence of simple reflections by $\underline w$ and shorten expressions of length $\geq 1$ by
\[\un{\scol}{k} := \underbrace{\scol\tcol\scol\ldots}_{k \text{ factors}},\quad\quad \unr{\scol}{k}:= \underbrace{\ldots \scol\tcol\scol}_{k \text{ factors}}, \]
similarly for $\tcol$. Omitting the underline means the corresponding product in $W$. We write $e= {}_\scol 0 ={}_\tcol 0$ for the identity and $w_0= {}_\scol m = {}_\tcol m$ for the longest element in $W$.
	
The restriction to dihedral groups makes it possible to obtain a closed formula for all Kazhdan-Lusztig basis elements simultaneously.

\begin{prop}\label{prop:smooth}
Let $(W,\{s,t\})$ be a dihedral group and $w\in W$, then
\[\uH_w = \sum_{x \leq w}v^{\ell(w)-\ell(x)}H_x.\]
\end{prop}

\begin{prf}
Induction over the length of $w$ (see \cite[Bsp. 2.3.]{He:KL}).
\end{prf}

For arbitrary Coxeter groups there is no such formula and the computation is inductively following for example the proof of \cref{thm:klbasis} in \cite[Theorem 2.1.]{S:KLP}. However, for a finite Coxeter group with longest element $w_0$ its corresponding Kazhdan-Lusztig basis element $\uH_{w_0}$ can always be computed with the formula in \cref{prop:smooth} (see \cite{KL:C}). Using \cref{prop:smooth} we can easily deduce the following lemma.

\begin{lem}\label{lem:omega}
For $w\in W$ we have $\omega(\uH_w)= \uH_{w^{-1}}$.
\end{lem}

\section{Soergel Bimodules}
Let $(W,S)$ be a Coxeter system of type $I_2(m), m\geq 3$. Recall the associated \emph{geometric representation} $\mathfrak{h}=\mathbb R \alpha_\scol^{\vee} \oplus \mathbb R \alpha_\tcol^{\vee}$ with its Cartan matrix (see \cite{H:CG}): 

\begin{align}\label{eqn:grep}
\begin{pmatrix}
2 & -2 \cos \left(\frac{\pi}{m} \right)\\
-2 \cos \left(\frac{\pi}{m} \right)&2
\end{pmatrix}.
\end{align}

We fix this realisation once and for all. For $\mathfrak h$ consider
$R:= S(\mathfrak h^*)= \bigoplus_{i\geq 0} S^i(\mathfrak h^*),$
the symmetric algebra on $\mathfrak h^*$, a graded $\DR$-algebra such that $\op{deg} (\mathfrak h^*) =2$. By construction, $W$ acts on $\mathfrak h$ and hence it acts on $\mathfrak h^*$ via the contragredient representation. Extending this action by grading preserving automorphism yields an action of $W$ on $R$. The ring of invariants of a single simple reflection $s\in S$ is denoted by $R^s \subseteq R$.

The two main module categories in this thesis are $R\mhy\op{Bim}$ and $R\mhy\op{gBim}$, the category of finitely generated $R$-bimodules and graded $R$-bimodules respectively (the latter with grading preserving morphisms). The category $R\mhy\op{gBim}$ is considered as a graded category with the grading shift down denoted by $\langle n \rangle$. For $M=\bigoplus M_i$ we define $M\langle n\rangle_i := M_{i+n}$. Moreover, for two graded bimodules $M$ and $N$ we write $\op{Hom}^\bullet(M,N):= \bigoplus_{n\in\DZ}(M,N \langle n \rangle)$ for the bimodule homomorphisms between $M$ and $N$ of all degrees.

\begin{rem}
Realisations of Coxeter groups can be defined in more generality as modules over a commutative domain and do not have to be symmetric either. For the more general case see \cite{EW:SC}.
\end{rem}

For $s\in S$ define the graded $R$-bimodule $B_s:= R \otimes_{R^s} R\langle 1\rangle $ which is a fixed graded lift of the $R$-bimodule $R\otimes_{R^s} R$. We often write $\otimes_s := \otimes_{R^s}$ and the tensor product structure $\otimes_R$ is denoted as juxtaposition. For a given sequence $\underline{w}=s_1s_2\cdots s_n$ define the corresponding \emph{Bott-Samelson bimodule} by
\[B_{\underline w}:=B_{s_1}\otimes_R B_{s_2}\otimes_R \cdots \otimes_R B_{s_n} \enspace = B_{s_1}B_{s_2}\cdots B_{s_n}. \]
The full monoidal subcategory $\bsbim$ of $R\mhy\op{Bim}$ generated by $B_s$ for $s\in S$ is called the \emph{category of Bott-Samelson bimodules}. Since we chose a fixed graded lift for $R\otimes_{R^s} R$, we have a graded lift for every Bott-Samelson bimodule. Finally, the \emph{category of Soergel bimodules} $\sbim$ is defined to be the Karoubian envelope of the additive closure of this graded version of $\bsbim$. It is crucial to distinguish that $\bsbim$ is a subcategory of $R\mhy\op{Bim}$, whilst $\sbim$ is a subcategory of $R\mhy\op{gBim}$ and therefore morphisms between Soergel bimodules are grading preserving. Observe that $\sbim$ is additive but not abelian. Soergel classified the indecomposable objects in $\sbim$ (see \cite[Theorem 6.14.]{S:KLP}).

\begin{thm}[Classification of indecomposable Soergel bimodules]\label{thm:clsbim}
Given any reduced expression $\underline{w}$ of $w \in W$, the Bott-Samelson $B_{\underline{w}}$ contains up to isomorphism a unique indecomposable summand $B_w$ which does not occur in $B_{\underline{y}}$ for any expression $\underline{y}$ of $y \in W$ with $\ell(y) < \ell(w)$. In addition, $B_w$ does up to isomorphism not depend on the reduced expression $\underline{w}$.
A complete set of representatives of the isomorphism classes of all indecomposable Soergel bimodules is given by
\[ \{ B_w\langle m\rangle \; \vert \; w \in W \text{ and } m \in \DZ \}. \]
\end{thm}

Note that the split Grothendieck group $\gps{\mathcal C}$ of an additive, monoidal and graded category $\mathcal C$ inherits a $\LP$-algebra structure. Soergel proved that the category of Soergel bimodules $\sbim$ categorifies the Hecke algebra $\hecke$ (see \cite[Theorem 1.10.]{S:BMP}) for certain (difficult) reflection faithful realisations over infinite fields of characteristic $\neq 2$. However, Libedinsky showed in \cite{L:Equiv} that this categorification works for the (easier) geometric representation as well.

\begin{thm}[Soergel's Categorification Theorem]\label{thm:scat}
For the geometric realisation $\mathfrak h$ there is a unique isomorphism of $\LP$-~algebras given by
\begin{align*}
\varepsilon: \; &\hecke \overset{\sim}{\longrightarrow} \gps{\sbim} \\
& \uH_s \mapsto [B_s] \text{.}
\end{align*}
\end{thm}

Using the standard pairing on $\hecke$ it is possible to describe the graded rank of the homomorphism space between two Soergel bimodules (see \cite[Theorem 5.15.]{S:BMP}).

\begin{thm}[Soergel's $\op{Hom}$-Formula]\label{thm:shom}
Given any two Soergel bimodules $B$ and $B'$, the homomorphism space $\op{Hom}_{\sbim}^{\bullet}(B, B')$ is free as a left (resp. right) $R$-module and its graded rank is given by $(\varepsilon^{-1} [B], \varepsilon^{-1} [B'] )$ where $(-, -)$ denotes the standard pairing on the Hecke algebra.
\end{thm}

Soergel constructed an inverse map to $\varepsilon$ which he called the character map. However, the construction is not explicit and he only conjectured what the pre-images of the indecomposable Soergel bimodules are for an arbitrary Coxeter group.

\begin{conj}[Soergel's Conjecture]\label{conj:s}
If $k$ is a field of characteristic $0$ then $\op{ch}(B_w)=\uH_w$.
\end{conj}

Soergel himself proved this for Weyl groups and dihedral groups (see \cite[Theorem 2]{S:HC}). The case of universal Coxeter groups was shown by Fiebig in \cite{F}. This conjecture is a very deep result and implies the Kazhdan-Lusztig conjectures (\cite{KL:C}). Elias and Williamson recently gave the first algebraic proof for an arbitrary Coxeter group with a fixed reflection-faithful representation over $\DR$ (see \cite{EW:HT}). By a result of Libedinsky in \cite{L:Equiv} this includes every finite Coxeter group with its geometric realisation.
In this paper we only consider dihedral groups and their geometric realisations over $\DR$ and thus we can use Soergel's Conjecture for our calculations to determine the dimensions of homomorphism spaces between Soergel bimodules.

Recall that we denote $R= \bigoplus_{i\geq 0} S^i(\mathfrak h^*)$ with maximal ideal $S^+ := \bigoplus S^i(\mathfrak h^*)$. Therefore we can view $\DR \cong R/S^+$ as an $R$-bimodule. To each Soergel bimodule $M$ we can associate the \emph{Soergel module} $\overline{M}:= M\otimes \DR$ where the right action is trivialised. Let $\mathcal S$ denote the \emph{category of Soergel modules}. For a Weyl group Soergel proved that there is an isomorphism
\begin{align}\label{rem:iso}
\op{Hom}^\bullet_{\mathcal S}(\oB,\oB^\prime) \cong \op{Hom}^\bullet_{\sbim}(B,B^\prime)\otimes_R\DR
\end{align}
for $B,B^\prime$ (see \cite[Thm. 2, Part 4.]{S:HC}). Moreover, he conjectured this isomorphism for every finite Coxeter group (see \cite[Thm. 2, Part 5.]{S:HC}) which he proved recently (\cite{S:Un}). For simplicity we write by abuse of notation $\op{Hom}_\smod(B,B^\prime)$ instead of $\op{Hom}_\smod(\oB,\oB^\prime)$ for Soergel bimodules $B,B^\prime$. The important results for this paper are summarised in the following theorem which is an immediate consequence of the above:

\begin{thm}\label{thm:dim}
Given a dihedral group $W$, we have for $x,y\in W$:
\begin{itemize}
\item $\op{grdim}_\DR\; \op{Hom}^\bullet_{\smod} (B_x,B_y)= (\uH_x,\uH_y)=\sum_{a\leq x,y} v^{\ell(x)+\ell(y)-2\ell(a)}$ 
\item $\op{dim}_\DR\; \op{Hom}^\bullet_{\smod} (B_x, B_y) = |W_{\leq x}\cap W_{\leq y}|$
\end{itemize}
where $W_{\leq x}:=\{w\in W\mid w\leq x\}$. In particular $\op{Hom}^\bullet_{\smod} (B_x,B_y)$ is concentrated in non-negative degrees $d$ such that $0\leq d\leq \ell(y)-\ell(x)$ for $x\leq y$ with $\dim_\DR \op{Hom}_{\smod} (B_x,B_y)= \delta_{x,y}$.
\end{thm}

\section{Soergel Calculus (in the dihedral case)}
By construction the category of Soergel bimodules $\sbim$ is the Karoubian envelope of $\bsbim$ and therefore it is enough to describe $\bsbim$ by planar graphs and identify the idempotents. In this section we introduce a diagrammatic approach to the category of Bott-Samelson bimodules $\bsbim$ following \cite{EW:SC}. These techniques are what we refer to as \emph{Soergel calculus}.

\subsection{Generators}
Recall that we fixed the geometric realisation $\mathfrak h$ for our given dihedral group $(W,\{\scol,\tcol\})$ of type $I_2(m)$ for $m\geq 3$. The Bott-Samelson bimodule $B_{\underline{w}}=B_{\scol} \otimes B_{\tcol} \otimes \cdots \otimes B_{\scol}$ for $\underline w = \scol\tcol\cdots\scol$ is completely determined by an ordered sequence of colours (or colured dots on a horizontal line). A morphism between two Bott-Samelson bimodules from  $B_{\underline w}$ to $B_{\underline{w^\prime}}$ is given by a linear combination of isotopy classes of decorated graphs with coloured edges in the planar stripe $\DR \times [0,1]$ such that the edges induce sequences of coloured dots on the bottom boundary $\DR \times \{0\}$ (resp. the top boundary $\DR \times \{1\}$) in each summand corresponding to $\underline w$ (resp. $\underline{w^\prime}$). In particular, these diagrams represent morphisms from the bottom sequence to the top sequence and therefore should be read from bottom to top.

\begin{table}
\begin{tabular}{c l l l}
   \centering
   \begin{tikzpicture}[baseline=-0.5ex]
      \draw[spot] (0,0) circle (0.5cm);
      \draw[sline] (0,-0.5) -- (0,0);
      \draw[sdot] (0,0) circle (2pt);
   \end{tikzpicture} & deg $1$ & $B_{\scol} \longrightarrow R$ & $a \otimes b \mapsto ab $ \\[0.5cm]
   \begin{tikzpicture}[baseline=-0.5ex]
      \draw[spot] (0,0) circle (0.5cm);
      \draw[sline] (0,0.5) -- (0,0);
      \draw[sdot] (0,0) circle (2pt);
   \end{tikzpicture} & deg $1$ & $R \longrightarrow B_{\scol}$ & 
      $1 \mapsto \frac{1}{2}(\alpha_s \otimes 1 + 1 \otimes \alpha_s ) $ \\[0.5cm]
   \begin{tikzpicture}[baseline=-0.5ex]
      \draw[spot] (0,0) circle (0.5cm);
      \clip (0,0) circle (0.5cm);
      \draw[sline] (-1,-1) -- (0,0);
      \draw[sline] (1,-1) -- (0,0);
      \draw[sline] (0,0) -- (0,1);
   \end{tikzpicture} & deg -$1$ & $B_{\scol} B_{\scol} \longrightarrow B_{\scol}$ 
      & $1 \otimes g \otimes 1 \mapsto \partial_s(g) \otimes 1 $ \\[0.5cm]
   \begin{tikzpicture}[baseline=-0.5ex]
      \draw[spot] (0,0) circle (0.5cm);
      \clip (0,0) circle (0.5cm);
      \draw[sline] (1,1) -- (0,0);
      \draw[sline] (-1,1) -- (0,0);
      \draw[sline] (0,0) -- (0,-1);
   \end{tikzpicture} & deg -$1$ & $B_{\scol} \longrightarrow B_{\scol} B_{\scol}$ 
      & $1 \otimes 1 \mapsto 1 \otimes 1 \otimes 1 $ \\[0.5cm]
   \begin{tikzpicture}[baseline=-0.5ex]
      \draw[spot] (0,0) circle (0.5cm);
      \clip (0,0) circle (0.5cm);
      \node (f) at (0,0) {$f$};
   \end{tikzpicture} & deg$f$ & $R \longrightarrow R$ & $1 \mapsto f$ \\[0.5cm]
   \begin{tikzpicture}[baseline=-0.5ex, outer sep=0.1mm, inner sep=0mm]
      \draw[spot] (0,0) circle (0.5cm);
      \clip (0,0) circle (0.5cm);         
      \draw[sline] (0,0) -- (210:0.5cm);
      \draw[tline] (0,0) -- (230:0.5cm);
      \draw[sline] (0,0) to (250:0.5cm);
      
      \draw[tline] (0,0) -- (150:0.5cm);
      \draw[sline] (0,0) -- (130:0.5cm);
      \draw[tline] (0,0) to (110:0.5cm);
      
      \draw[dotted, thick] (0, -0.3) to (0.25, -0.3);
      \draw[dotted, thick] (0, 0.3) to (0.25, 0.3);
      
   \end{tikzpicture} & deg $0$ & $\underbrace{B_{\scol} B_{\tcol} B_{\scol} \cdots}_{m_{st} \text{ factors}}
      \longrightarrow \underbrace{B_{\tcol} B_{\scol} B_{\tcol} \cdots}_{m_st \text{ factors}} $ &     
\end{tabular}
\newline 
\caption{Generating morphisms and their degrees in $\mathcal{D}$}\label{tab:genmorph}
\end{table}

\begin{defn}
For a dihedral group $(W,S)$ define $\mathcal D = \mathcal D_{(W,S)}$ to be the $\DR$-linear monoidal category as follows: The objects are sequences $\underline w$ in $S$ (which are denoted sometimes by $B_{\underline w}$). The empty sequence $\emptyset$ is often denoted by $\op{\textbf{1}}$. The $\op{Hom}$-spaces are $\DZ$-graded $\DR$-vector spaces generated by the diagrams in \cref{tab:genmorph} modulo local relations. The monoidal structure is the concatenation of sequences.
\end{defn}

For $s\in S$ the \emph{Demazure operator} $\partial_s: R \to R^s$ in \cref{tab:genmorph} is defined as $\partial_s(f):= \frac{f-sf}{\alpha_s}$. The first two morphisms in \cref{tab:genmorph} are called \emph{dots}, whereas the second two morphisms are called \emph{trivalent vertices} and the last morphism is called \emph{$2m_{\scol,\tcol}$-valent vertex}. The explicit formula for the $2m_{\scol,\tcol}$-valent vertex is very difficult. Therefore we only explain what the morphism does. In the case of a dihedral group of type $I_2(m)$ the longest element $w_0$ can be expressed as $st\cdots = ts\cdots$ with $m$ factors on each side and by \cref{thm:clsbim} both $B_{\scol}B_{\tcol}\cdots$ and $B_{\tcol}B_{\scol}\cdots$ contain $B_{w_0}$ as summand with multiplicity $1$. The $2m_{\scol,\tcol}$-vertex is the projection and inclusion of this summand and therefore uniquely determined up to a scalar.

Libedinsky showed in \cite{L:SB} that the morphisms from \cref{tab:genmorph} generate all morphisms in $\bsbim$. For the compositions of a trivalent vertex with a dot we define \emph{caps} and \emph{cups} as follows:

\begin{alignat}{5}
\begin{tikzpicture}[baseline=(current  bounding  box.center)]
\draw[sline] (0,0) arc (0:180:0.5cm);
\end{tikzpicture}
\enspace & := \enspace 
\begin{tikzpicture}[baseline=(current  bounding  box.center)]
\draw[sline] (0,0) arc (0:180:0.5cm);
\draw[sline] (-0.5,0.5) to (-0.5,0.8);
\draw[sdot] (-0.5,0.8) circle (2pt);
\end{tikzpicture}
& \quad\quad\quad\quad \text{ and }\quad\quad\quad\quad &
\begin{tikzpicture}[baseline=(current  bounding  box.center)]
\draw[sline] (0,0) arc (180:360:0.5cm);
\end{tikzpicture}
\enspace & := \enspace 
\begin{tikzpicture}[baseline=(current  bounding  box.center)]
\draw[sline] (0,0) arc (180:360:0.5cm);
\draw[sline] (0.5,-0.5) to (0.5,-0.8);
\draw[sdot] (0.5,-0.8) circle (2pt);
\end{tikzpicture}.
\end{alignat}

\subsection{Relations}
Since we work with dihedral groups, only one- and two-colour relations can occur. For Coxeter groups of rank $\geq 3$ there exist the three-colour or \emph{Zamolodzhikov} relations which are more difficult (see \cite{EW:SC}).

\subsubsection{The one-colour relations}
The object $B_{\scol}$ is a Frobenius object in $R\mhy\op{Bim}$ (see \cite{EK}) where the dots correspond to unit and counit while the trivalent vertices correspond to the multiplication and comultiplication. Altogether we obtain that any one-coloured diagram is isotopy invariant and we have the following three non-polynomial relations:
\begin{alignat}{3}
   \label{eqn:genass}
   \begin{tikzpicture}[baseline=-0.5ex, 
                       sline, inner sep=0mm, outer sep=0mm,scale=0.7]
      \draw (-0.25,-0.75) .. controls (-0.25,-0.5) .. (0,-0.25);
      \draw (0.25,-0.75) .. controls (0.25,-0.5) .. (0,-0.25);                
      \draw (-0.25,0.75) .. controls (-0.25,0.5) .. (0,0.25);
      \draw (0.25,0.75) .. controls (0.25,0.5) .. (0,0.25);  
      \draw (0,-0.25)to (0,0.25);               
   \end{tikzpicture} 
   \enspace & = \enspace \begin{tikzpicture}[baseline=-0.5ex, 
                       sline, inner sep=0mm, outer sep=0mm,scale=0.7]
      \draw (-0.5,-0.75) .. controls (-0.5,-0.5) .. (-0.25,0);
      \draw (-0.5,0.75) .. controls (-0.5,0.5) .. (-0.25,0);
      \draw (0.5,-0.75) .. controls (0.5,-0.5) .. (0.25,0);
      \draw (0.5,0.75) .. controls (0.5,0.5) .. (0.25,0);
      \draw (-0.25,0) to (0.25,0);
   \end{tikzpicture}
   &, \quad\quad\quad
   \begin{tikzpicture}[baseline=-0.5ex, 
                       sline, inner sep=0mm, outer sep=0mm,scale=0.7]
      \draw (0, -0.75) to (0,0.75);
      \draw (0, 0) to (0.25, 0);
      \draw[fill, thin] (0.25, 0) circle (2pt);
   \end{tikzpicture} 
   \enspace & = \enspace \begin{tikzpicture}[baseline=-0.5ex, 
                       sline, inner sep=0mm, outer sep=0mm,scale=0.7]
      \draw (0, -0.75) to (0, 0.75);
   \end{tikzpicture} 
   &, \quad\quad\quad 
   \begin{tikzpicture}[baseline=-0.5ex, 
                       inner sep=0mm, outer sep=0mm,scale=0.7]      
      \draw[spot] (0,0) circle (0.75cm);
      \draw[sline] (0,0) circle (0.25cm);
      \draw[sline] (0,-0.75) to (0,-0.25);
   \end{tikzpicture}\enspace & = \enspace 0.
\end{alignat}
We refer to these relations as \emph{general associativity}, \emph{general unit} and \emph{the needle}. Furthermore, there are two relations involving polynomials:
\begin{alignat}{2}
    \label{eqn:poly}
   \begin{tikzpicture}[baseline=-0.5ex, 
                       inner sep=0mm, outer sep=0mm,scale=0.7]      
      \draw[spot] (0,0) circle (0.75cm);
      \draw[sline] (0,-0.3) to (0,0.3);
      \draw[sdot] (0,0.3) circle (2pt);
      \draw[sdot] (0,-0.3) circle (2pt);
   \end{tikzpicture} \enspace & = \enspace
   \begin{tikzpicture}[baseline=-0.5ex, 
                       inner sep=0mm, outer sep=0mm,scale=0.7]  
      \draw[spot] (0,0) circle (0.75cm);
      \node (v) at (0,0) {$\alpha_{\scol}$};
   \end{tikzpicture}
&,\quad\quad\quad\quad
   \begin{tikzpicture}[baseline=-0.5ex, 
                       inner sep=0mm, outer sep=1mm,scale=0.7]  
      \draw[spot] (0,0) circle (0.75cm);
      \clip (0,0) circle (0.75cm);
      \draw[sline] (0,-1) to node[left, color=black] (f) {\small$f$} (0,1);
   \end{tikzpicture} \enspace & = \enspace
   \begin{tikzpicture}[baseline=-0.5ex, 
                       inner sep=0mm, outer sep=1mm,scale=0.7]  
      \draw[spot] (0,0) circle (0.75cm);
      \clip (0,0) circle (0.75cm);
      \draw[sline] (0,-1) to node[right, color=black] (sf) {\small$\scol f$} (0,1);
   \end{tikzpicture} \enspace + \enspace
   \begin{tikzpicture}[baseline=-0.5ex, 
                       inner sep=0mm, outer sep=1mm,scale=0.7]  
      \draw[spot] (0,0) circle (0.75cm);
      \clip (0,0) circle (0.75cm);
      \draw[sline] (0,-1) to (0,-0.4);
      \draw[sdot] (0,-0.4) circle (2pt);
      \draw[sline] (0,1) to (0,0.4);
      \draw[sdot] (0,0.4) circle (2pt);
      \node (parsf) at (0,0) {\small$\partial_{\scol}(f)$};
   \end{tikzpicture}.
\end{alignat}   
As a direct consequence we have:

\begin{align}
\begin{tikzpicture}[baseline=-0.5ex, 
                       inner sep=0mm, outer sep=1mm,scale=0.7]  
      \draw[spot] (0,0) circle (0.75cm);
      \clip (0,0) circle (0.75cm);
      \draw[sline] (0,-1) to (0,-0.4);
      \draw[sdot] (0,-0.4) circle (2pt);
      \draw[sline] (0,1) to (0,0.4);
      \draw[sdot] (0,0.4) circle (2pt);
   \end{tikzpicture}
   \enspace & = \enspace \frac{1}{2}\enspace 
   \left(
   \begin{tikzpicture}[baseline=-0.5ex, 
                       inner sep=0mm, outer sep=1mm,scale=0.7]  
      \draw[spot] (0,0) circle (0.75cm);
      \clip (0,0) circle (0.75cm);
      \draw[sline] (0,-1) to node[left, color=black] (f) {\small$\alpha_\scol$} (0,1);
   \end{tikzpicture} 
   \enspace + \enspace
   \begin{tikzpicture}[baseline=-0.5ex, 
                       inner sep=0mm, outer sep=1mm,scale=0.7]  
      \draw[spot] (0,0) circle (0.75cm);
      \clip (0,0) circle (0.75cm);
      \draw[sline] (0,-1) to node[right, color=black] (f) {\small$\alpha_\scol$} (0,1);
   \end{tikzpicture} 
   \right). \label{eqn:polyslide2}
\end{align}

\subsubsection{The two-colour relations}
There are three two-colour relations. The first is the cyclicity of the $2m_{\scol,\tcol}$-vertex. The other two relations specify how the $2m_{\scol,\tcol}$-vertex interacts with trivalent vertices and dots. For $m=m_{\scol,\tcol}$ we have the \emph{two-colour associativity} depending on the parity of $m$:

\begin{alignat}{2}
m \text{ even:}\quad
\begin{tikzpicture}[baseline=-0.5ex, 
                       inner sep=0mm, outer sep=1mm, scale=0.5]
 \draw[tline] (-0.2,0) to (1.2,0);
 \draw[tline] (0.8,-1.5) to (0.8,1.5);
 \draw[tline] (1.2,0) .. controls (2.2,-1) .. (2.2,-1.5);
 \draw[tline] (1.2,0) ..  controls (2.2,1) .. (2.2,1.5);
 \draw[sline] (0.4,-1.5) .. controls (0.4,-1) .. (0.8,0); 
 \draw[sline] (0.8,0) .. controls (1.4,1) .. (1.4,1.5);
 \draw[sline] (0.4,1.5) .. controls (0.4,1) .. (0.8,0); 
 \draw[sline] (0.8,0) .. controls (1.4,-1) .. (1.4,-1.5);
 \draw[dotted,thick] (0.98,-1.3) to (1.22,-1.3);
 \draw[dotted,thick] (0.98,1.3) to (1.22,1.3);
\end{tikzpicture} 
\enspace & = \enspace \begin{tikzpicture}[baseline=-0.5ex, 
                       inner sep=0mm, outer sep=1mm,scale=0.5]
 \draw[tline] (0.8,-1.5) to (0.8,1.5);
 \draw[tline] (-0.8,0) to (-0.6,0);
 \draw[tline] (-0.6,0) .. controls (0,0.7) .. (0.8,0.7);
 \draw[tline] (0.8,0.7) .. controls (1.6,0.7) and (2.2,1) .. (2.2,1.5);
 \draw[tline] (-0.6,0) .. controls (0,-0.7) .. (0.8,-0.7);
 \draw[tline] (0.8,-0.7) .. controls (1.6,-0.7) and (2.2,-1) .. (2.2,-1.5);
 \draw[sline] (0.4,-1.5) .. controls (0.4,-1) .. (0.8,-0.7);
 \draw[sline] (0.8,-0.7) .. controls (1.4,0) .. (0.8,0.7);
 \draw[sline] (0.8,0.7) .. controls (0.4,1) .. (0.4,1.5);
 \draw[sline] (1.4,-1.5) .. controls (1.4,-1) .. (0.8,-0.7);
 \draw[sline] (0.8,-0.7) .. controls (0.4,0) .. (0.8,0.7);
 \draw[sline] (0.8,0.7) .. controls (1.4,1) .. (1.4,1.5);
 \draw[dotted,thick] (0.98,-1.3) to (1.22,-1.3);
 \draw[dotted,thick] (0.9,0) to (1.14,0);
 \draw[dotted,thick] (0.98,1.3) to (1.22,1.3);
 \end{tikzpicture}
  & \qquad\qquad m \text{ odd:} \quad
  \begin{tikzpicture}[baseline=-0.5ex, 
                       inner sep=0mm, outer sep=1mm,scale=0.5]
 \draw[tline] (-0.2,0) to (0.8,0);
 \draw[tline] (0.4,-1.5) .. controls (0.4,-1) .. (0.8,0);
 \draw[tline] (0.4,1.5) .. controls (0.4,1) .. (0.8,0);
 \draw[tline] (0.8,0) .. controls (1.6,1) .. (1.6,1.5);
 \draw[tline] (0.8,0) .. controls (1.6,-1) .. (1.6,-1.5);
 \draw[sline] (0.8,0) to (1.2,0);
 \draw[sline] (0.8,0) .. controls (1.2,1) .. (1.2,1.5);
 \draw[sline] (0.8,0) .. controls (1.2,-1) .. (1.2,-1.5);
 \draw[sline] (1.2,0) .. controls (2,-1) .. (2,-1.5);
 \draw[sline] (1.2,0) ..  controls (2,1) .. (2,1.5);
 \draw[sline] (0,-1.5) .. controls (0,-1) .. (0.8,0); 
 \draw[sline] (0,1.5) .. controls (0,1) .. (0.8,0); 
 \draw[dotted,thick] (0.68,-1.3) to (0.92,-1.3);
 \draw[dotted,thick] (0.68,1.3) to (0.92,1.3);
\end{tikzpicture} 
\enspace & = \enspace \begin{tikzpicture}[baseline=-0.5ex, 
                       inner sep=0mm, outer sep=1mm,scale=0.5]
 \draw[tline] (0.4,-1.5) .. controls (0.4,-1) .. (0.8,-0.7);
 \draw[tline] (0.8,0.7) .. controls (0.4,1) .. (0.4,1.5);
 \draw[tline] (0.8,-0.7) .. controls (0.4,0) .. (0.8,0.7);
 \draw[tline] (-0.8,0) to (-0.6,0);
 \draw[tline] (-0.6,0) .. controls (0,0.7) .. (0.8,0.7);
 \draw[tline] (-0.6,0) .. controls (0,-0.7) .. (0.8,-0.7);
 \draw[tline] (0.8,-0.7) .. controls (1.6,0) .. (0.8,0.7);
 \draw[tline] (1.6,-1.5) .. controls (1.6,-1) .. (0.8,-0.7);
 \draw[tline] (0.8,0.7) .. controls (1.6,1) .. (1.6,1.5);
 \draw[sline] (0.8,-0.7) .. controls (1.2,0) .. (0.8,0.7);
 \draw[sline] (1.2,-1.5) .. controls (1.2,-1) .. (0.8,-0.7);
 \draw[sline] (0.8,0.7) .. controls (1.2,1) .. (1.2,1.5);
 \draw[sline] (0.8,0.7) .. controls (1.6,0.7) and (2,1) .. (2,1.5);
 \draw[sline] (0.8,-0.7) .. controls (1.6,-0.7) and (2,-1) .. (2,-1.5);
 \draw[sline] (0,-1.5) .. controls (0,-1) .. (0.8,-0.7);
 \draw[sline] (0.8,0.7) .. controls (0,1) .. (0,1.5);
 \draw[sline] (0.8,-0.7) .. controls (0,0) .. (0.8,0.7);
 \draw[dotted,thick] (0.68,-1.3) to (0.92,-1.3);
 \draw[dotted,thick] (0.68,1.3) to (0.92,1.3);
 \draw[dotted,thick] (0.68,0) to (0.92,0);
 \end{tikzpicture}
\label{eqn:2ctri}
 \end{alignat}
 
The other relation allows us to express a diagram involving a $2m_{\scol,\tcol}$-vertex and a dot as a linear combination of diagrams without the $2m_{\scol,\tcol}$-vertex. This procedure depends on the parity of the integer $m$ as well:

\begin{alignat}{2}
m \text{ even:} \quad
\begin{tikzpicture}[baseline=-0.5ex,scale=0.7]
\draw[tline] (-0.2,-1) .. controls (-0.2,-0.5) .. (0,0);
\draw[tline] (-0.2,1) .. controls (-0.2,0.5) .. (0,0);
\draw[tline] (-1,0) to (0.8,0);
\draw[tdot] (0.8,0) circle (2pt);
\draw[sline] (-0.5,-1) .. controls (-0.5,-0.5).. (0,0);
\draw[sline] (-0.5,1) .. controls (-0.5,0.5).. (0,0);
\draw[sline] (0.5,-1) .. controls (0.5,-0.5).. (0,0);
\draw[sline] (0.5,1) .. controls (0.5,0.5).. (0,0);
\draw[dotted,thick] (0.03,-0.8) to (0.27,-0.8);
\draw[dotted,thick] (0.03,0.8) to (0.27,0.8);
\end{tikzpicture}
\enspace & = \enspace
\begin{tikzpicture}[baseline=-0.5ex,scale=0.7]
\draw[tline] (-0.2,-1) to (-0.2,1);
\draw[tline] (0.5,-1) to (0.5,1);
\draw[tline] (-1,0) to (0,0);
\draw[sline] (-0.5,-1) to (-0.5,1);
\draw[sline] (0,0) to (1,0);
\draw[sline] (1,-1) to (1,1);
\draw (0,0) -- (0,0) node[fill=red!20, draw, rounded corners, font=\scriptsize] {$JW_{m-1}$};
\draw[dotted,thick] (0.03,-0.8) to (0.27,-0.8);
\draw[dotted,thick] (0.03,0.8) to (0.27,0.8);
\end{tikzpicture}
 & \qquad\qquad m \text{ odd:} \quad
\begin{tikzpicture}[baseline=-0.5ex,scale=0.7]
\draw[tline] (-0.2,-1) .. controls (-0.2,-0.5) .. (0,0);
\draw[tline] (-0.2,1) .. controls (-0.2,0.5) .. (0,0);
\draw[tline] (0.5,-1) .. controls (0.5,-0.5).. (0,0);
\draw[tline] (0.5,1) .. controls (0.5,0.5).. (0,0);
\draw[tline] (-1,0) to (0,0);
\draw[sline] (-0.5,-1) .. controls (-0.5,-0.5).. (0,0);
\draw[sline] (-0.5,1) .. controls (-0.5,0.5).. (0,0);
\draw[sline] (0,0) to (0.8,0);
\draw[sdot] (0.8,0) circle (2pt);
\draw[dotted,thick] (0.03,-0.8) to (0.27,-0.8);
\draw[dotted,thick] (0.03,0.8) to (0.27,0.8);
\end{tikzpicture}
\enspace & = \enspace
\begin{tikzpicture}[baseline=-0.5ex,scale=0.7]
\draw[tline] (-0.2,-1) to (-0.2,1);
\draw[tline] (-1,0) to (0,0);
\draw[tline] (0,0) to (1,0);
\draw[tline] (1,-1) to (1,1);
\draw[sline] (-0.5,-1) to (-0.5,1);
\draw[sline] (0.5,-1) to (0.5,1);
\draw (0,0) -- (0,0) node[fill=red!20, draw, rounded corners,font=\scriptsize] {$JW_{m-1}$};
\draw[dotted,thick] (0.03,-0.8) to (0.27,-0.8);
\draw[dotted,thick] (0.03,0.8) to (0.27,0.8);
\end{tikzpicture}
\label{eqn:2cdot}
\end{alignat}

Note that the Jones-Wenzl morphism $JW_{m-1}$ is an $\DR$-linear combination of graphs consisting only of dots and trivalent vertices. This morphism will be discussed briefly in the next paragraph.

\subsection{Jones-Wenzl morphisms}

\subsubsection{Gau{\ss}'s $q$-numbers}
Before we can define the Jones-Wenzl projector it is helpful to recall Gau{\ss}'s $q$-numbers. We use them to give formulas for the Jones-Wenzl projector for all dihedral groups simultaneously. Gau{\ss}'s $q$-numbers are defined (see e.g. \cite{Ja:QG}) by
\[[n]_q:= \frac{q^n-q^{-n}}{q-q^{-1}}=q^{-n+1}+^{-n+3}+ \cdots + q^{n-3} + q^{n-1} \in \DZ[q^{\pm 1}].\]
It is convenient to define $[0]:= 0$. In most cases we omit the index and write $[n]$ instead of $[n]_q$.
Observe that $[n]$ is the character of $L(n-1)$, the simple $\mathfrak s\mathfrak l_2(\DC)$-module of dimension $n$. Via the Clebsch-Gordon formula we obtain two out of many useful identities for $q$-numbers: 
\begin{align}
[2][n] & \enspace = \enspace [n+1] + [n-1]. \label{eqn:q1} \\
[n]^2 & \enspace = \enspace [n-1][n+1] + [1].\label{eqn:q5}
\end{align}

We can specialise $q$ to a value $\zeta$ which we denote by $[n]_{\zeta}$. If $\zeta = e^{2\pi i/2m}\in \DC$, i.e. a primitive $2m$-th root of unity, we obtain algebraic integers $[n]_\zeta \in \DR$. By choice $\zeta$ is primitive and therefore we have $\zeta^m=-1$.  Thus, we get the following identities:
\begin{alignat}{3}
[m]_\zeta \enspace & = \enspace 0, &\qquad
[m-i]_\zeta \enspace &= \enspace [i]_\zeta, &\qquad
[m+i]_\zeta \enspace &= \enspace -[i]_\zeta.\label{eqn:q7}
\end{alignat}

Recall the geometric representation of a dihedral group as in \cref{eqn:grep} and note that $[2]_\zeta = \zeta + \zeta^{-1}= 2\cos \left(\frac{\pi}{m} \right)=-a_{s,t}$. Hence, the geometric representation can be encoded in the Cartan matrix
\begin{align}\label{eqn:grepq}
\begin{pmatrix}
2 & -[2]\\
-[2]&2
\end{pmatrix}
\end{align}
simultaneously for all dihedral groups. For a certain dihedral group of type $I_2(m)$ we only have to specialise $q$ to a primitive $2m$-th root of unity. 

\begin{lem}\label{lem:actionroots}
For $i\geq 1$ we have $i_\scol = \ldots \tcol\scol \in W$ of length $i$ and thus
\[ 
i_\scol(\alpha_\tcol) = 
\begin{cases}
[i]_\zeta \alpha_\tcol + [i+1]_\zeta \alpha_\scol \quad\quad \text{ if $i$ is odd,} \\
[i+1]_\zeta \alpha_\tcol + [i]_\zeta \alpha_\scol \quad\quad  \text{ if $i$ is even.}
\end{cases}
\]
\end{lem}

\begin{prf}
We prove this by induction on $i$. For $i=1$ we have:
\[\scol(\alpha_\tcol) = \alpha_\tcol - \langle \alpha_\tcol,\alpha_\scol^\vee \rangle\ \alpha_\scol = \alpha_\tcol + [2]_\zeta \alpha_\scol.\]
Now let $i > 1$. There are two cases to distinguish depending on the parity of $i$. If $i$ is odd then
\begin{align*}
(i+1)_\scol(\alpha_t)
\enspace  = \enspace
\tcol(i_\scol(\alpha_t))
\enspace & = \enspace
\tcol([i]_\zeta \alpha_t + [i+1]_\zeta \alpha_s) \\
\enspace & = \enspace
(-[i]\zeta+ [2]_\zeta [i+1]_\zeta)\alpha_\tcol + [i+1]_\zeta \alpha_\scol \\
\enspace & = \enspace
(-[i]_\zeta + [i]_\zeta + [i+2]_\zeta)\alpha_\tcol + [i+1]_\zeta\alpha_\scol \\
\enspace & = \enspace
[i+2]_\zeta \alpha_\tcol + [i+1]_\zeta \alpha_\scol.
\end{align*}
The other case is similar which finishes the proof.
\end{prf}

\subsubsection{Jones-Wenzl projectors}

We give a short introduction to Jones-Wenzl projectors with a recursive computation formula in the Temperley-Lieb algebra setting where its origins lie (see \cite{J,We}). For a detailed discussion we refer the reader to \cite{CK} or \cite{E:DC}.

The \emph{Temperley-Lieb algebra} $TL_n$ on $n$ strands is a diagram algebra over $\DZ[\delta]$. It has a basis consisting of crossingless matching with $n$ points on bottom and top each. The multiplication is given by vertical concatenation of diagrams such that circles evaluate to scalar $-\delta$. The crossingless matching on $n$ strands is the unit in $TL_n$ and is denoted by $\mathbf 1_n$. This algebra is contained in the \emph{Temperley-Lieb category} which is closely related to the quantum group $U:=U_q(\mathfrak s \mathfrak l_2)$ of $\mathfrak s\mathfrak l_2$ via the base change $\delta \mapsto [2]_q=q+q^{-1}$ (see \cite{E:DC} for technical details).

Let $V_k$ be the irreducible representation of highest weight $q^k$ and let $V=V_1$. The highest non-zero projection from $V^{\otimes n}$ to $V_n$ is known as the \emph{Jones-Wenzl projector} $JW_n\in TL_n$. 

\begin{rem}
This definition of the Jones-Wenzl projector only works while working over the complex numbers. This categorification over the integers is harder and more subtle and was carefully treated in \cite{CK,FSS,R}.
\end{rem}

The following proposition states some important properties of $JW_n$ (see \cite[Propositon 3.2.2.]{KaLi},\cite[Lemma 13.2.]{Lic} and \cite[Section 2.2.]{CK}).

\begin{prop}\label{prop:defJW}
The Jones-Wenzl projector $JW_n$ satisfies the following properties:
\begin{itemize}
\item $JW_n$ is the unique map which is killed when any cap is attached on top or any cup on bottom, and for which the coefficient of $\mathbf{1}_n$ is $1$.
\item $JW_n$ is invariant under horizontal/vertical reflection.
\item The ideal $\langle JW_n \rangle \unlhd TL_n$ has rank $1$. 
\item Any element $x\in TL_n$ acts on $JW_n$ by its coefficient of $\mathbf{1}_n$.
\end{itemize}
\end{prop}

Note that the first property gives an alternative way of defining $JW_n$. For our computations we need the following recursive formula for $JW_n$(see \cite[Theorem 3.5.]{FK}):

\begin{alignat}{4}
\begin{tikzpicture}[baseline=-0.5ex,bline]
\draw (-0.55,-1) to (-0.55,1);
\draw (0.55,-1) to (0.55,1);
\draw (0.3,-1) to (0.3,1);
\draw (0.05,-1) to (0.05,1);
\draw (0,0) -- (0,0) node[fill=red!20, draw] {$JW_{n+1}$};
\draw[dotted,thick] (-0.37,-0.85) to (-0.13,-0.85);
\draw[dotted,thick] (-0.37,0.85) to (-0.13,0.85);
\end{tikzpicture}
\enspace &= \enspace 
\begin{tikzpicture}[baseline=-0.5ex,bline]
\draw (-0.4,-1) to (-0.4,1);
\draw (0.65,-1) to (0.65,1);
\draw (0.4,-1) to (0.4,1);
\draw (0.15,-1) to (0.15,1);
\draw (0,0) -- (0,0) node[fill=red!20, draw] {$JW_{n}$};
\draw[dotted,thick] (-0.245,-0.85) to (-0.005,-0.85);
\draw[dotted,thick] (-0.245,0.85) to (-0.005,0.85);
\end{tikzpicture}
\enspace & + \sum_{i=1}^n \frac{[i]_q}{[n+1]_q}
\enspace & \enspace 
\begin{tikzpicture}[baseline=-0.5ex,bline]
\draw (-0.4,-1) to (-0.4,0.3);
\draw (0.15,-1) to (0.15,0.3);
\draw (0.4,-1) to (0.4,0.3);
\draw (0.9,0.3) arc (0:180:0.25);
\draw (0.9,0.3) to (0.9,-1);
\draw (-0.4,1) arc (180:360:0.25);
\draw (-0.4,1.2) to (-0.4,1.2) node {i};
\draw (-0.4,0.3) .. controls (-0.4,0.5) and (-0.65,0.8) .. (-0.65,1);
\draw (0.15,0.3) .. controls (0.15,0.5) and (0.4,0.8) .. (0.4,1);
\draw (0,0) -- (0,0) node[fill=red!20, draw] {$JW_{n}$};
\draw[dotted,thick] (-0.245,-0.85) to (-0.005,-0.85);
\end{tikzpicture}.\label{thm:jw}
\end{alignat}

In the following we state the Jones-Wenzl projectors for $n=1,2,3$.
\begin{align*}
JW_1 \enspace &\textstyle{=} \enspace 
	\begin{tikzpicture}[baseline=-0.5ex,bline]
	\draw (0,-0.4) to (0,0.4);
	\end{tikzpicture}\hspace{4cm}
JW_2 \enspace \textstyle{=} \enspace
	\begin{tikzpicture}[baseline=-0.5ex,bline]
	\draw (0,-0.4) to (0,0.4);
	\draw (0.3,-0.4) to (0.3,0.4);
	\end{tikzpicture}
	\enspace \textstyle{+ \enspace\frac{1}{[2]}} \enspace
	\begin{tikzpicture}[baseline=-0.5ex,bline]
	\draw (0,0.4) arc (180:360:0.15);
	\draw (0.3,-0.4) arc (0:180:0.15);
	\end{tikzpicture}  \\[0.5cm]
JW_3 \enspace &\textstyle{=} \enspace
	\begin{tikzpicture}[baseline=-0.5ex,bline]
	\draw (0,-0.4) to (0,0.4);
	\draw (0.3,-0.4) to (0.3,0.4);
	\draw (0.6,-0.4) to (0.6,0.4);
	\end{tikzpicture}
	\enspace \textstyle{+\enspace \frac{[2]}{[3]}} \enspace
	\begin{tikzpicture}[baseline=-0.5ex,bline]
	\draw (0,0.4) arc (180:360:0.15);
	\draw (0.3,-0.4) arc (0:180:0.15);
	\draw (-0.3,-0.4) to (-0.3,0.4);
	\end{tikzpicture}
	\enspace \textstyle{+\enspace \frac{[2]}{[3]}} \enspace
	\begin{tikzpicture}[baseline=-0.5ex,bline]
	\draw (-0.3,0.4) arc (180:360:0.15);
	\draw (0,-0.4) arc (0:180:0.15);
	\draw (0.3,-0.4) to (0.3,0.4);
	\end{tikzpicture}
	\enspace \textstyle{+ \enspace\frac{1}{[3]}} \enspace
	\begin{tikzpicture}[baseline=-0.5ex,bline]
	\draw (0.3,0.4) arc (180:360:0.15);
	\draw (0.3,-0.4) arc (0:180:0.15);
	\draw (0.6,-0.4) .. controls (0.6,-0.1) and (0,0.1) .. (0,0.4);
	\end{tikzpicture}
	\enspace \textstyle{+ \enspace\frac{1}{[3]}} \enspace
	\begin{tikzpicture}[baseline=-0.5ex,bline]
	\draw (0,0.4) arc (180:360:0.15);
	\draw (0.6,-0.4) arc (0:180:0.15);
	\draw (0,-0.4) .. controls (0,-0.1) and (0.6,0.1) .. (0.6,0.4);
	\end{tikzpicture}
\end{align*}

Any crossingless matching in $TL_{m-1}$ divides the planar stripe in $m$ regions which we can colour alternatingly with red and blue (for $\scol$ and $\tcol$). This results in the definition of the \emph{two-coloured Temperley-Lieb algebra} which we omit here, details can be found in \cite{E:DC,EW:SC}. Since our chosen realisation is symmetric we can treat a blue circle surrounded by red just as a red circle surrounded by blue and thus evaluate both to the same value. For each diagram there are two possible colourings. Each of those coloured crossingless matchings yields a coloured graph. Deformation retract each region into a tree consisting of trivalent and univalent vertices, colour those resulting trees according to the region. In that way we can associate to each coloured Jones-Wenzl projector a \emph{Jones-Wenzl morphism} $JW_n$. We state the Jones-Wenzl morphisms (with red appearing in the far left region) for $n=1,2,3$.

\begin{align*}
JW_1 \enspace & = \enspace 
\begin{tikzpicture}[scale=0.7]
\draw[sline] (-0.5,0) -- (-0.2,0);
\draw[sdot] (-0.2,0) circle (2pt);
\draw[tline] (0.5,0) -- (0.2,0);
\draw[tdot] (0.2,0) circle (2pt);
\end{tikzpicture} \hspace{4cm}
JW_2 \enspace  = \enspace
\begin{tikzpicture}[baseline=-0.5ex,scale=0.7]
\draw[sline] (-0.6,0) -- (-0.3,0);
\draw[sdot] (-0.3,0) circle (2pt);
\draw[sline] (0.6,0) -- (0.3,0);
\draw[sdot] (0.3,0) circle (2pt);
\draw[tline] (0,-0.6) -- (0,0.6);
\end{tikzpicture}
\enspace + \frac{1}{[2]}\enspace 
\begin{tikzpicture}[baseline=-0.5ex,scale=0.7]
\draw[tline] (0,-0.6) -- (0,-0.3);
\draw[tdot] (0,-0.3) circle (2pt);
\draw[tline] (0,0.6) -- (0,0.3);
\draw[tdot] (0,0.3) circle (2pt);
\draw[sline] (-0.4,0) -- (0.4,0);
\end{tikzpicture} \\[0.5cm]
JW_3 \enspace & = \enspace
\begin{tikzpicture}[baseline=-0.5ex,scale=0.7]
\draw[tline] (-0.15,-0.6) -- (-0.15,0.6);
\draw[sline] (0.15,-0.6) -- (0.15,0.6);
\draw[sdot] (-0.4,0) circle (2pt);
\draw[sline] (-0.7,0) -- (-0.4,0);
\draw[tdot] (0.4,0) circle (2pt);
\draw[tline] (0.4,0) -- (0.7,0);
\end{tikzpicture}
\enspace + \frac{[2]}{[3]}\enspace
\begin{tikzpicture}[baseline=-0.5ex,scale=0.7]
\draw[sline] (0.15,0.6) -- (0.15,0.3);
\draw[sdot] (0.15,0.3) circle (2pt);
\draw[sline] (0.15,-0.6) -- (0.15,-0.3);
\draw[sdot] (0.15,-0.3) circle (2pt);
\draw[tline] (0,0) -- (0.7,0);
\draw[sdot] (-0.4,0) circle (2pt);
\draw[sline] (-0.7,0) -- (-0.4,0);
\draw[tline] (-0.15,0.6) .. controls (-0.15,0.4).. (0,0);
\draw[tline] (-0.15,-0.6) .. controls (-0.15,-0.4).. (0,0);
\end{tikzpicture}
\enspace + \frac{[2]}{[3]}\enspace
\begin{tikzpicture}[baseline=-0.5ex,scale=0.7]
\draw[tline] (-0.15,0.6) -- (-0.15,0.3);
\draw[tdot] (-0.15,0.3) circle (2pt);
\draw[tline] (-0.15,-0.6) -- (-0.15,-0.3);
\draw[tdot] (-0.15,-0.3) circle (2pt);
\draw[sline] (0,0) -- (-0.7,0);
\draw[tdot] (0.4,0) circle (2pt);
\draw[tline] (0.7,0) -- (0.4,0);
\draw[sline] (0.15,0.6) .. controls (0.15,0.4).. (0,0);
\draw[sline] (0.15,-0.6) .. controls (0.15,-0.4).. (0,0);
\end{tikzpicture}
\enspace + \frac{1}{[3]}\enspace
\begin{tikzpicture}[baseline=-0.5ex,scale=0.7]
\draw[sline] (0.15,0.6) -- (0.15,0.3);
\draw[sdot] (0.15,0.3) circle (2pt);
\draw[tline] (-0.15,-0.6) -- (-0.15,-0.3);
\draw[tdot] (-0.15,-0.3) circle (2pt);
\draw[sline] (-0.7,0) .. controls (0,0) and (0.15,-0.3) .. (0.15,-0.6);
\draw[tline] (0.7,0) .. controls (0,0) and (-0.15,0.3) .. (-0.15,0.6);
\end{tikzpicture}
\enspace + \frac{1}{[3]}\enspace
\begin{tikzpicture}[baseline=-0.5ex,scale=0.7]
\draw[tline] (-0.15,0.6) -- (-0.15,0.3);
\draw[tdot] (-0.15,0.3) circle (2pt);
\draw[sline] (0.15,-0.6) -- (0.15,-0.3);
\draw[sdot] (0.15,-0.3) circle (2pt);
\draw[sline] (-0.7,0) .. controls (0,0) and (0.15,0.3) .. (0.15,0.6);
\draw[tline] (0.7,0) .. controls (0,0) and (-0.15,-0.3) .. (-0.15,-0.6);
\end{tikzpicture}
\end{align*}

The Jones-Wenzl morphism above is not yet a morphism in $\mathcal D$ but can be plugged into another diagram to obtain a graph in the planar stripe (see \cref{eqn:2cdot}, for technical details see \cite{E:DC}). We often write $JW$ instead of $JW_{m-1}$ when specialised to $\zeta$. Using the two colour relations (see \cref{eqn:2ctri,eqn:2cdot}) we obtain the following relation between the $m_{\scol,\tcol}$-valent vertex and the Jones-Wenzl projector $JW_{m-1}$. Since the latter is an idempotent in the Temperley-Lieb algebra (see \cite{E:DC}), the right-hand sides in \ref{eqn:idem} are also idempotents. Thus, the $m_{\scol,\tcol}$-valent vertex can be used to construct idempotents.

\begin{alignat}{2}
m \text{ even:}
\begin{tikzpicture}[baseline=-0.5ex,scale=0.5]
 \draw[tline] (0.4,-1.5) .. controls (0.4,-1) .. (0.8,-0.7);
 \draw[tline] (0.8,0.7) .. controls (0.4,1) .. (0.4,1.5);
 \draw[tline] (0.8,-0.7) .. controls (0.4,0) .. (0.8,0.7);
 \draw[tline] (0.8,-0.7) .. controls (1.6,0) .. (0.8,0.7);
 \draw[tline] (1.6,-1.5) .. controls (1.6,-1) .. (0.8,-0.7);
 \draw[tline] (0.8,0.7) .. controls (1.6,1) .. (1.6,1.5);
 \draw[sline] (0.8,-0.7) .. controls (1.2,0) .. (0.8,0.7);
 \draw[sline] (1.2,-1.5) .. controls (1.2,-1) .. (0.8,-0.7);
 \draw[sline] (0.8,0.7) .. controls (1.2,1) .. (1.2,1.5);
 \draw[sline] (0,-1.5) .. controls (0,-1) .. (0.8,-0.7);
 \draw[sline] (0.8,0.7) .. controls (0,1) .. (0,1.5);
 \draw[sline] (0.8,-0.7) .. controls (0,0) .. (0.8,0.7);
 \draw[dotted,thick] (0.68,-1.3) to (0.92,-1.3);
 \draw[dotted,thick] (0.68,1.3) to (0.92,1.3);
 \draw[dotted,thick] (0.68,0) to (0.92,0);
\end{tikzpicture}
\enspace & = \enspace
\begin{tikzpicture}[baseline=-0.5ex,scale=0.7]
\draw[tline] (-0.3,-1) to (-0.3,1);
\draw[sline] (-1,0) to (0,0);
\draw[tline] (0,0) to (1,0);
\draw[tline] (1,-1) to (1,1);
\draw[sline] (-1,-1) to (-1,1);
\draw[sline] (0.3,-1) to (0.3,1);
\draw (0,0) -- (0,0) node[fill=red!20, draw, rounded corners,font=\scriptsize] {$JW_{m-1}$};
\draw[dotted,thick] (-0.12,-0.85) to (0.12,-0.85);
\draw[dotted,thick] (-0.12,0.85) to (0.12,0.85);
\end{tikzpicture}
&\qquad\qquad m \text{ odd:}\quad
\begin{tikzpicture}[baseline=-0.5ex,scale=0.5]
 \draw[tline] (0.8,-1.5) to (0.8,1.5);
 \draw[sline] (0.4,-1.5) .. controls (0.4,-1) .. (0.8,-0.7);
 \draw[sline] (0.8,-0.7) .. controls (1.4,0) .. (0.8,0.7);
 \draw[sline] (0.8,0.7) .. controls (0.4,1) .. (0.4,1.5);
 \draw[sline] (1.4,-1.5) .. controls (1.4,-1) .. (0.8,-0.7);
 \draw[sline] (0.8,-0.7) .. controls (0.4,0) .. (0.8,0.7);
 \draw[sline] (0.8,0.7) .. controls (1.4,1) .. (1.4,1.5);
 \draw[dotted,thick] (0.9,0) to (1.14,0);
 \draw[dotted,thick] (0.98,1.5) to (1.22,1.5);
 \draw[dotted,thick] (0.98,-1.5) to (1.22,-1.5);
\end{tikzpicture}
\enspace & = \enspace
\begin{tikzpicture}[baseline=-0.5ex,scale=0.7]
\draw[tline] (-0.3,-1) to (-0.3,1);
\draw[sline] (-1,0) to (0,0);
\draw[sline] (0,0) to (1,0);
\draw[sline] (1,-1) to (1,1);
\draw[sline] (-1,-1) to (-1,1);
\draw (0,0) -- (0,0) node[fill=red!20, draw, rounded corners,font=\scriptsize] {$JW_{m-1}$};
\draw[dotted,thick] (-0.12,-1) to (0.12,-1);
\draw[dotted,thick] (-0.12,1) to (0.12,1);
\end{tikzpicture}
\label{eqn:idem}
\end{alignat}

The two defining properties of $JW_n$ (see \cref{prop:defJW} or \cite[Claim 4.4]{E:DC}) are crucial and we use them repeatedly in our category $\mathcal D$. The first one is that the coefficient of one single graph is $1$, namely the graph which becomes the identity morphism if trivalent vertices are attached to both sides (precisely the first summands in the example above). The second and more important property is that $JW$ is killed by "cups" and "caps"

\begin{align}
\begin{tikzpicture}[baseline=-0.5ex]
\draw[tline] (0,-0.5) -- (0,0.5);
\draw[tdot] (0,0.5) circle (2pt);
\draw[sline] (-0.25,-0.5) -- (-0.25, 0.4);
\draw[sline] (0.25,-0.5) -- (0.25, 0.4);
\draw[sline] (-0.25,0.4) .. controls (-0.25,0.6) .. (0,0.7);
\draw[sline] (0.25,0.4) .. controls (0.25,0.6) .. (0,0.7);
\draw[sline] (0,0.7) -- (0,0.9);
\draw[tline] (-0.65,0) -- (0.65,0);
\draw[white] (0,-0.5) -- (0,-0.9);
\draw (0,0) -- (0,0) node[fill=red!20, draw, rounded corners] {$JW$};
\end{tikzpicture}
\enspace = \enspace 0.
\end{align}

We refer to this property as \emph{death by pitchfork}.

\subsection{Light leaf and double leaf morphism}\label{sec:light}
The concept of light leaves and double leaves was introduced by Libedinsky in \cite{L:SB} for Soergel bimodules. Elias and Williamson transferred those morphisms into the combinatorial and diagrammatic setting of $\mathcal D$ (see\cite{EW:SC}) where they used the double leaves to prove the equivalence in \cref{thm:equiv}. For all (technical) details and results we refer to \cite{EW:SC}.
The set of light leaf morphisms $LL_{\underline{x},e}$ forms a (left) $R$-basis of $\op{Hom}_{\mathcal D}(\underline{x},\op{\textbf{1}})$. This basis can be constructed inductively and is parametrised by the subexpressions of $\underline x$ which are equal to $e$. Composing the these light leaf morphisms we obtain the set of double leaf morphisms $\DL\DL_{\underline{x},\underline{y}}$ which forms a (left) $R$-basis of $\op{Hom}_{\mathcal D}(\underline{x},\underline{y})$ (see \cite[Theorem 6.11]{EW:SC}). Conesequently, $\DL\DL_{\underline{x},\underline{y}}$ is constructed inductively, too, and it is indexed by certain subexpressions of both $\underline x$ and $\underline y$. Hence all $\op{Hom}$-spaces in $\mathcal D$ are free graded (left) $R$-modules.

\subsection{Equivalence of categories}
Following \cite{E:DC} we define an $\DR$-linear monoidal functor $\mathcal F: \mathcal D_{(W,S)} \to \bsbim$ which maps a sequence $\underline w$ to the Bott-Samelson bimodule $B_{\underline w}$. On the morphisms $\mathcal F$ acts as indicated in \cref{tab:genmorph}. The next result is due to Elias (see \cite{E:DC}):
\begin{thm}\label{thm:equiv}
The $\DR$-linear monoidal functor $\mathcal F$ is well-defined and yields an equivalence of monoidal categories: 
\[\mathcal F:\mathcal D_{(W,S)} \to \bsbim.\]
Moreover, $\mathcal F$ induces an equivalence of monoidal additive categories:
\[\mathcal F: \op{Kar}(\mathcal D_{(W,S)}) \to \sbim.\]
\end{thm}

\begin{rem}
The result in \cref{thm:equiv} were preceded by results by Elias and Khovanov in type $A$ (see \cite{EK}) and by Libedinsky (see \cite{L:rightangled}) for the right-angled case. Furthermore, the above results were generalised by Elias and Williamson (see \cite[Theorem 6.28.]{EW:SC}) to all Coxeter groups with a fixed reflection-faithful realisation.
\end{rem}

\section{Koszul and Quasi-Hereditary Algebras}

\begin{defn}[\cite{BGS}]
A graded ring $A= \bigoplus_{i\geq 0} A_i$ is called \emph{Koszul} if $A_0$ is semisimple and $A_0$ admits a \emph{linear resolution} as a graded left $A$-module, i.e. 
 \[\cdots \to P^2 \to P^1 \to P^0 \twoheadrightarrow A_0 \]
 where the maps are grading preserving and $P^i$ is generated by $(P^i)_i$ as a left $A$-module.
\end{defn}

\begin{rem}
Equivalently one could ask that every simple $A$-module admits a linear resolution.
\end{rem}

\begin{defn}[\cite{BGS}]
A graded ring $A= \bigoplus_{i\geq 0} A_i$ is called \emph{quadratic} if $A_0$ is semisimple and $A$ is generated over $A_0$ by $A_1$ with relations in degree $2$.
\end{defn}

For a given $k$-vector space $V$ the tensor algebra $T_kV$, the symmetric algebra $SV$ and the exterior algebra $\Lambda V$ are quadratic rings (with the usual grading, that is $\op{deg} X =1$). In particular $k[X]$ and $k[X]/(X^2)$ are quadratic. If there are some finiteness conditions we can define a dual version of a quadratic ring:

\begin{defn}[\cite{BGS}]
A graded ring $A= \bigoplus_{i\geq 0} A_i$ is called \emph{left finite} if all $A_i$ are finitely generated as left $A_0$-modules. For a left finite quadratic ring $A=T_{A_0}A_1/(R)$ its \emph{quadratic dual} is defined as $A^!=T_{A_0}A_1^*/(R^\perp)$ where $R^\perp \subset A_1^*\otimes_{A_0}A_1^*\cong (A_1\otimes_{A_0} A_1)^*$ with respect to the standard pairing. Here, $T_{A_0} A_1$ denotes the free tensor algebra of the $A_0$-bimodule $A_1$. The ring $A$ is \emph{self-dual} if $A\cong A^!$.
\end{defn}

For any positively graded algebra $A=\bigoplus_{i\geq 0}A_i$ the degree $0$ part $A_0$ is an $A$-module. Consider the graded ring $E(A):=\op{Ext}_A^\bullet(A_0,A_0)$ of self-extensions of $A_0$. We call $E(A)$ the \emph{Koszul dual} of $A$. Beilinson, Ginzburg and Soergel showed that there exist isomorphisms relating a Koszul algebra $A$ to its Koszul dual $E(A)$ and quadratic dual $A^!$ (see \cite[Corollary 2.3.3., Theorem 2.10.1., Theorem 2.10.2.]{BGS}).
\begin{thm}[\cite{BGS}]\label{thm:BGS}
Let $A=\bigoplus_{i\geq0} A_i$ be a Koszul ring. Then $A$ is quadratic. If additionally $A$ is left finite, then there are canonical isomorphisms $E(A)\cong (A^!)^{op}$ and $E(E(A))\cong A$.
\end{thm}

\begin{rem}
If $A$ is a positively graded algebra there exists a quadratic duality functor on the bounded derived categories
\begin{align*}
K: \mathcal D^b(A\mhy\op{gMod})\to D^b(A^!\mhy\op{gMod})
\end{align*}
which is an equivalence of categories if and only if $A$ is Koszul. This is treated completely and more general in \cite[Thm. 1.2.6.]{BGS} and \cite[Thm. 30]{MOS}. In conclusion Koszul algebras are certain quadratic algebras with additional nice homological properties. For example, under the above equivalence the standard $t$-structure maps to the non-standard $t$-structure on the dual side given by linear complexes of graded projective modules (\cite[Thm. 12]{MOS}). If $A$ is even Koszul self-dual, then this yields a second interesting $t$-structure in $\mathcal D^b(A\mhy\op{gMod})$.
\end{rem}

The algebras we are interested in admit a quasi-hereditary structure which simplifies the later proofs in a crucial way.

\begin{defn}[\cite{Don:S}]Let $A$ be a finite dimensional algebra over a field $k$ with a finite partially ordered set $(\Lambda,\leq)$ indexing the simple left $A$-modules $\{L(\lambda)\}_{\lambda\in\Lambda}$ and let $P(\lambda)$ be the indecomposable projective cover of $L(\lambda)$. A collection of left $A$-modules $\{\Delta(\lambda)\}_{\lambda\in\Lambda}$ defines a quasi-hereditary structure on $(A,(\Lambda,\leq))$ if
\begin{itemize}
\item for $\lambda \in\Lambda $ there exists a surjective A-module homomorphism $\pi:P(\lambda)\twoheadrightarrow\Delta(\lambda)$ such that $\op{ker}\,\pi$ has a $\Delta$-filtration with subquotients isomorphic to some $\Delta(\mu_i)$ where $\mu_i > \lambda$,
\item for $\lambda \in\Lambda $ there exists a surjective A-module homomorphism $\pi^\prime:\Delta(\lambda)\twoheadrightarrow L(\lambda)$ such that $\op{ker}\,\pi^\prime$ has a composition series with composition factors only isomorphic to some $L(\mu_i)$ where $\mu_i < \lambda$.
\end{itemize}
The $\Delta(\lambda)$ are called (left) \emph{standard modules}. (Right) standard modules can be defined similarly.
\end{defn}

The following theorem is due to \'{A}goston, Dlab and Luka\'{a}cs (see \cite[Theorem 1]{ADL}):
\begin{thm}\label{thm:qhk}
Let $(A,(\Lambda,\leq))$ be a graded quasi-hereditary algebra. If both left and right standard modules admit linear resolutions (i.e. A is standard Koszul) then $A$ is Koszul (i.e. all simple modules admit linear resolutions).
\end{thm}

\begin{rem}
We would like to point out the close relationship between the double leaf morphisms in \cref{sec:light} and the filtration obtained by the standard modules.
\end{rem}

\section{Realisation as a Path algebra of a Quiver with Relations}
Recall that $(W,\{\scol,\tcol\})$ is a Coxeter system of type $I_2(m), m \geq 3$ with its geometric representation $\mathfrak{h}=\mathbb R \alpha_\scol^v \oplus \mathbb R \alpha_\tcol^v$. In this section we analyse the structure of the graded endomorphism algebra $\mathcal A:= \op{End}^\bullet_{\smod}\left( \textbf{B}\right)\cong \op{End}^\bullet_{\sbim}\left( \textbf{B}\right)\otimes_R \DR$ where $\textbf B:= \bigoplus_{x\in W} B_x$. Define $\zeta:=e^{2\pi i/2m}\in  \DC$.

\begin{figure}
\[\xymatrix@!=0.8pc{&w_0\ar@{<->}[ld]\ar@{<->}[rd]& \\
{}_\tcol (m-1)\ar@{<->}[d]\ar@{<->}[rrd]&&{}_\scol (m-1)\ar@{<->}[lld]|!{[ll];[d]}\hole\ar@{<->}[d]\\
\quad \scriptscriptstyle{\vdots} \quad \ar@{<->}[d]\ar@{<->}[rrd] & \scriptscriptstyle{\vdots} & \quad \scriptscriptstyle{\vdots} \quad \ar@{<->}[lld]|!{[ll];[d]}\hole\ar@{<->}[d] \\
\quad \tcol\scol\quad\ar@{<->}[d]\ar@{<->}[rrd]&& \quad \scol\tcol \quad\ar@{<->}[d]\ar@{<->}[lld]|!{[ll];[d]}\hole\\
\quad \tcol \quad\ar@{<->}[dr]&& \quad \scol \quad\ar@{<->}[dl]\\
& e& } \]
\caption{The Hasse graph of type $I_2(m)$}.\label{fig:hasse}
\end{figure}
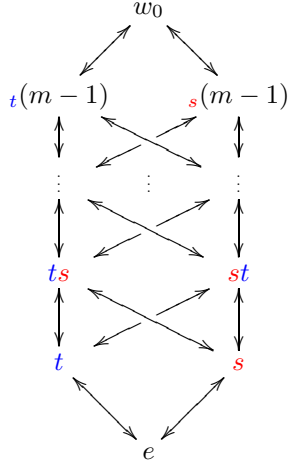

\begin{defn}\label{def:gg}
Let $Q_{m}:=(Q_0,Q_1)$ be the directed quiver with the following vertex and arrow sets 
\begin{align*}
Q_0 & :=W= \{e,{}_\scol k,{}_\tcol k,w_0\mid 1 \leq k \leq m -1\}, \\
Q_1 & :=\{w\to w^\prime \mid |l(w)-l(w^\prime)|=1\}.
\end{align*}  
We write $(w,w^\prime)$ for an arrow $w\to w^\prime$. For the fixed representation $\mathfrak h$ of $I_2(m)$ define the set $R_{m}^\mathfrak h \subset \DR Q_{m}$ consisting of the following relations for all $2\leq i \leq m-1$ and $0 \leq j \leq m-2$ (plus the ones with the roles of $\scol$ and $\tcol$ switched):

\begin{small}
\begin{align}
(e,\scol,e) = {}& 0, \label{eqn:1}\\
(\scol,\scol\tcol,\scol) = {}& 0, \label{eqn:2}\\
({}_\scol i,{}_\scol (i+1),{}_\scol i)= {}& \textstyle{\frac{[i-1]_\zeta}{[i]_\zeta}}({}_\scol i,{}_\scol (i-1),{}_\scol i),\label{eqn:3} \\
(\scol,\tcol\scol,\scol) = {} & \textstyle{-[2]_\zeta}(\scol,e,\scol),\label{eqn:4}\\
\begin{split}
({}_\scol i,{}_\tcol (i+1),{}_\scol i)= {}& {\textstyle([i-1]_\zeta-[i+1]_\zeta)} ({}_\scol i,{}_\scol (i-1),{}_\scol i) -\textstyle{\frac{[i+1]_\zeta}{[i]_\zeta}}({}_\scol i,{}_\tcol (i-1),{}_\scol i),
\end{split}\label{eqn:5}\\
(\scol,\scol\tcol,\tcol) = {}& (\scol,e,\tcol), \label{eqn:6}\\
(\scol,\tcol\scol,\tcol) = {}& (\scol,e,\tcol), \label{eqn:7}\\
({}_\scol i, {}_\scol (i+1), {}_\tcol i) = {}& \textstyle{\frac{1}{[i]_\zeta}} ({}_\scol i,{}_\scol (i-1),{}_\tcol i)
 + ({}_\scol i,{}_\tcol (i-1),{}_\tcol i),\label{eqn:8}\\
({}_\scol i, {}_\tcol (i+1), {}_\tcol i) = {}& ({}_\scol i,{}_\scol (i-1),{}_\tcol i) + \textstyle{\frac{1}{[i]_\zeta}} ({}_\scol i,{}_\tcol (i-1),{}_\tcol i), \label{eqn:9}\\
({}_\scol j,{}_\scol (j+1),{}_\scol (j+2) )= {}& ({}_\scol j ,{}_\tcol (j+1),{}_\scol (j+2)),\label{eqn:10}\\
({}_\scol j, {}_\scol (j+1) ,{}_\tcol (j+2) )= {}& ({}_\scol j ,{}_\tcol (j+1) ,{}_\tcol (j+2)),\label{eqn:11}\\
({}_\scol (j+2) , {}_\scol (j+1) ,{}_\scol j) = {}& ({}_\scol (j+2) , {}_\tcol (j+1) ,{}_\scol j),\label{eqn:12}\\
({}_\scol (j+2) , {}_\scol (j+1),{}_\tcol j) = {}& ({}_\scol (j+2) , {}_\tcol (j+1) ,{}_\tcol j).\label{eqn:13}
\end{align}
\end{small}

We refer to $Q_{m}$ as the Hasse graph of type $I_2(m)$ (c.f. \cref{fig:hasse}) and the set $R_m^\mathfrak h$ is called the set of dihedral relations. Define $\textbf P_m$ to be the $\DR$-algebra $\DR Q_{m}/(R_{m}^\mathfrak h)$.
\end{defn}

The algebra $\Pm$ inherits the naturally grading of $\DR Q_m$ by path length since all relations are homogenous (of degree $2$).

\begin{rem}\label{rem:rel}
For $i=1$ relation \ref{eqn:2} could be included in relation \ref{eqn:3} whilst relation \ref{eqn:4} cannot be included in relation \ref{eqn:5}. Moreover, for $i=m-1$ relation \ref{eqn:3} and \ref{eqn:5} agree and relation \ref{eqn:8} and \ref{eqn:9}.
\end{rem}

\begin{rem}
The above formulas are exactly the relations which hold for the morphisms between the Soergel modules which we show later in detail. The coefficients are integral in type $A_2$, but in general they are real numbers. Stroppel showed in \cite{St:CO} that for types $B_2$ and $G_2$ one may obtain rational (or even integral coefficients) whilst choosing different bases for the morhpism spaces. However, those relations are not symmetric any more.
\end{rem}

By definition the vertices of $Q_{m}$ are indexed by the dihedral group of type $I_2(m)$. Hence we have the Bruhat order on the vertex set. For a vertex $x\in Q_0$ define $V_{\leq x}:=\{y\in Q_0 \mid y\leq x\}$. Analysing the relations carefully we obtain the following lemma and propositions:

\begin{lem}\label{lem:dim}
The path algebra $\textbf{P}_m$ has (finite) dimension $\sum_{x,y} |V_{\leq x}\cap V_{\leq y}|$.
\end{lem}

\begin{rem}\label{rem:paths}
\cref{lem:dim} means in particular that paths from $x$ to $y$ in $\Pm$ cannot be the same if the lowest vertex they pass through is different.
\end{rem}

\begin{prop}\label{prop:simple}
The isomorphism classes of simple graded $\Pm$-modules up to grading shift are in bijection with the vertices of $Q_{m}$ via $w \mapsto~L(w)$ where 
$L(w)_{w^\prime}:=
\begin{cases} \mathbb R &\text{if } w=w^\prime \\
0 & \text{if } w\neq w^\prime
\end{cases}$ \quad and the obvious maps.
\end{prop}

\begin{prop}\label{prop:op}
There exists an isomorphism $\Pm \cong \Pm^{op}$.
\end{prop}

We can now state the main result of this section.

\begin{thm}\label{thm:isom}
There is an isomorphism of graded algebras $\textbf{P}_m \cong \mathcal A$.
\end{thm}

\subsection{Proof of Theorem \ref{thm:isom}}

Consider the assignment $w \mapsto B_{\underline w}$ for $w \in Q_0\setminus \{w_0\}$ and $w_0\mapsto B_{{}_\scol\underline{m}}\oplus B_{{}_\tcol\underline{m}}$. We extend this assignment to the arrows for $0\leq j \leq m-1$.

\begin{alignat*}{3}
\left({}_\scol j \to {}_\scol (j+1)\right) \enspace & \longmapsto \enspace
\overbrace{\begin{tikzpicture}[baseline=-0.5ex]
\draw[sline] (0,-0.4) -- (0,0.4);
\draw[tline] (0.2,-0.4) -- (0.2,0.4);
\draw[sline] (0.4,-0.4) -- (0.4,0.4);
\draw[dotted, thick] (0.58,0) -- (0.82,0);
\draw[bbline] (1,0)--(1,0.4);
\draw[bbdot] (1,0) circle (2pt);
\end{tikzpicture}}^{j+1 \text{ lines}}
\enspace \quad &\text{and}  \quad 
\left({}_\scol (j+1) \to {}_\scol j\right) \enspace & \longmapsto \enspace
\underbrace{\begin{tikzpicture}[baseline=-0.5ex]
\draw[sline] (0,-0.4) -- (0,0.4);
\draw[tline] (0.2,-0.4) -- (0.2,0.4);
\draw[sline] (0.4,-0.4) -- (0.4,0.4);
\draw[dotted, thick] (0.58,0) -- (0.82,0);
\draw[bbline] (1,0)--(1,-0.4);
\draw[bbdot] (1,0) circle (2pt);
\end{tikzpicture}}_{j+1 \text{ lines}}\\[0.5cm]
\left({}_\scol j \to {}_\tcol (j+1)\right) \enspace & \longmapsto \enspace
\overbrace{\begin{tikzpicture}[baseline=-0.5ex]
\draw[tline] (0,0)--(0,0.4);
\draw[tdot] (0,0) circle (2pt);
\draw[sline] (0.2,-0.4) -- (0.2,0.4);
\draw[tline] (0.4,-0.4) -- (0.4,0.4);
\draw[dotted, thick] (0.58,0) -- (0.82,0);
\draw[bbline] (1,-0.4) -- (1,0.4);
\end{tikzpicture}}^{j+1 \text{ lines}}
\enspace \quad& \text{and}  \quad
\left({}_\tcol (j+1) \to {}_\scol j \right) \enspace & \longmapsto \enspace
\underbrace{\begin{tikzpicture}[baseline=-0.5ex)]
\draw[tline] (0,0)--(0,-0.4);
\draw[tdot] (0,0) circle (2pt);
\draw[sline] (0.2,-0.4) -- (0.2,0.4);
\draw[tline] (0.4,-0.4) -- (0.4,0.4);
\draw[dotted, thick] (0.58,0) -- (0.82,0);
\draw[bbline] (1,-0.4) -- (1,0.4);
\end{tikzpicture}}_{j+1 \text{ lines}}
\end{alignat*}

The colour of the black line depends on the parity of $j$ and on the colour of the far left line. Switching colours yields the other half of the assignment. Arrows adjacent to $w_0$ are treated differently since $B_{w_0}$ is embedded diagonally in $B_{{}_\scol\underline{m}}\oplus B_{{}_\tcol\underline{m}}$:

\begin{alignat*}{3}
\left({}_\scol(m-1) \to w_0  \right) \enspace & \longmapsto \enspace
\begin{pmatrix}
\begin{tikzpicture}[baseline=-0.5ex]
\draw[sline] (0,-0.4) -- (0,0.4);
\draw[tline] (0.2,-0.4) -- (0.2,0.4);
\draw[sline] (0.4,-0.4) -- (0.4,0.4);
\draw[dotted, thick] (0.58,0) -- (0.82,0);
\draw[bbline] (1,0)--(1,0.4);
\draw[bbdot] (1,0) circle (2pt);
\end{tikzpicture} \\ \\
\begin{tikzpicture}[baseline=-0.5ex]
\draw[tline] (0,0)--(0,0.4);
\draw[tdot] (0,0) circle (2pt);
\draw[sline] (0.2,-0.4) -- (0.2,0.4);
\draw[tline] (0.4,-0.4) -- (0.4,0.4);
\draw[dotted, thick] (0.58,0) -- (0.82,0);
\draw[bbline] (1,-0.4) -- (1,0.4);
\end{tikzpicture}
\end{pmatrix},
\enspace &&\\[0.5cm] 
\left(w_0 \to {}_\scol(m-1)\right) \enspace & \longmapsto \enspace
\frac{1}{2}\begin{pmatrix}
\begin{tikzpicture}[baseline=-0.5ex]
\draw[sline] (0,-0.4) -- (0,0.4);
\draw[tline] (0.2,-0.4) -- (0.2,0.4);
\draw[sline] (0.4,-0.4) -- (0.4,0.4);
\draw[dotted, thick] (0.58,0) -- (0.82,0);
\draw[bbline] (1,0)--(1,-0.4);
\draw[bbdot] (1,0) circle (2pt);
\end{tikzpicture} & 
\begin{tikzpicture}[baseline=-0.5ex]
\draw[tline] (0,0)--(0,-0.4);
\draw[tdot] (0,0) circle (2pt);
\draw[sline] (0.2,-0.4) -- (0.2,0.4);
\draw[tline] (0.4,-0.4) -- (0.4,0.4);
\draw[dotted, thick] (0.58,0) -- (0.82,0);
\draw[bbline] (1,-0.4) -- (1,0.4);
\end{tikzpicture}
\end{pmatrix},&&
\end{alignat*}

For $\scol$ and $\tcol$ switched the assignment is similar. Note that the images of the arrows are only morphisms between Bott-Samelson bimodules. In order to get a morphism between the corresponding Soergel bimodules we have to pre-/post-compose with the idempotents associated to the Soergel bimodules, that is
For $e,\scol,\tcol$ the idempotents are the identity since the corresponding Bott-Samelson bimodules are indecomposable. The idempotent corresponding to $B_{{}_\scol i}$ for $2 \leq i \leq m-1$ is the Jones-Wenzl morphism $JW_i$ where the quantum numbers must be specialised to the appropriate $2m$-th root of unity, hence the coefficients depend on the dihedral group. In other words,

\begin{alignat*}{2}
2\leq i \leq m-1:\quad e_{{}_\scol i} &\mapsto \begin{tikzpicture}[baseline=-0.5ex]
\draw[sline] (-0.5,-0.4) -- (-0.5,0.4);
\draw[tline] (-0.25,-0.4) -- (-0.25,0.4);
\draw[bbline] (0.25,-0.4) -- (0.25,0.4);
\draw[bbline] (0.5,-0.4) -- (0.5,0.4);
\draw[dotted, thick] (-0.12,0.35) -- (0.12,0.35);
\draw[dotted, thick] (-0.12,-0.35) -- (0.12,-0.35);
\draw[bbline] (0.5,0) -- (0,0);
\draw[sline] (-0.5,0) -- (0,0);
\draw (0,0) -- (0,0) node[fill=red!20, draw] {$\scriptstyle{JW_i}$};
\end{tikzpicture}
&\qquad\qquad 
e_{w_0} & \mapsto
\frac{1}{2}
\begin{pmatrix}
\begin{tikzpicture}[baseline=-0.5ex]
\draw[sline] (-0.5,-0.4) -- (-0.5,0.4);
\draw[tline] (-0.25,-0.4) -- (-0.25,0.4);
\draw[bbline] (0.25,-0.4) -- (0.25,0.4);
\draw[bbline] (0.5,-0.4) -- (0.5,0.4);
\draw[dotted, thick] (-0.12,0.35) -- (0.12,0.35);
\draw[dotted, thick] (-0.12,-0.35) -- (0.12,-0.35);
\draw[bbline] (0.5,0) -- (0,0);
\draw[sline] (-0.5,0) -- (0,0);
\draw (0,0) -- (0,0) node[fill=red!20, draw] {$\scriptstyle{JW}$};
\end{tikzpicture}
 & \begin{tikzpicture}[baseline=-0.5ex]
\draw[tline] (-0.5,-0.4) .. controls (-0.4,0) .. (0,0);
\draw[sline] (-0.25,-0.4) .. controls (-0.2,0).. (0,0);
\draw[bbline] (0.5,-0.4) .. controls (0.4,0) .. (0,0);
\draw[sline] (-0.5,0.4) .. controls (-0.4,0) .. (0,0);
\draw[tline] (-0.25,0.4) .. controls (-0.2,0).. (0,0);
\draw[bbline] (0.5,0.4) .. controls (0.4,0) .. (0,0);
\draw[dotted, thick] (-0.005,0.35) -- (0.235,0.35);
\draw[dotted, thick] (-0.005,-0.35) -- (0.235,-0.35);
\draw (0,0) -- (0,0) node[fill=red!20, draw,ellipse,scale=.8] {$2m$};
\end{tikzpicture}\\ 
 & \\
\begin{tikzpicture}[baseline=-0.5ex]
\draw[sline] (-0.5,-0.4) .. controls (-0.4,0) .. (0,0);
\draw[tline] (-0.25,-0.4) .. controls (-0.2,0).. (0,0);
\draw[bbline] (0.5,-0.4) .. controls (0.4,0) .. (0,0);
\draw[tline] (-0.5,0.4) .. controls (-0.4,0) .. (0,0);
\draw[sline] (-0.25,0.4) .. controls (-0.2,0).. (0,0);
\draw[bbline] (0.5,0.4) .. controls (0.4,0) .. (0,0);
\draw[dotted, thick] (-0.005,0.35) -- (0.235,0.35);
\draw[dotted, thick] (-0.005,-0.35) -- (0.235,-0.35);
\draw (0,0) -- (0,0) node[fill=red!20, draw,ellipse,scale=.8] {$2m$};
\end{tikzpicture}
 &\begin{tikzpicture}[baseline=-0.5ex]
\draw[tline] (-0.5,-0.4) -- (-0.5,0.4);
\draw[sline] (-0.25,-0.4) -- (-0.25,0.4);
\draw[bbline] (0.25,-0.4) -- (0.25,0.4);
\draw[bbline] (0.5,-0.4) -- (0.5,0.4);
\draw[dotted, thick] (-0.12,0.35) -- (0.12,0.35);
\draw[dotted, thick] (-0.12,-0.35) -- (0.12,-0.35);
\draw[bbline] (0.5,0) -- (0,0);
\draw[tline] (-0.5,0) -- (0,0);
\draw (0,0) -- (0,0) node[fill=red!20, draw] {$\scriptstyle{JW}$};
\end{tikzpicture}
\end{pmatrix} =: JW^\Delta.
\end{alignat*}

where $e_x$ denotes the trivial path at vertex $x$. Similarly for $\tcol$ with all colours switched. The idempotent of $B_{w_0}$ is a $2\times 2$-matrix since we consider $B_{w_0}$ a as submodule of $B_{{}_\scol\underline{m}}\oplus B_{{}_\tcol\underline{m}}$. By the classification in \cref{thm:clsbim} the indecomposable bimodule $B_{w_0}$ is a direct summand of $B_{{}_\scol\underline{m}}$ and $B_{{}_\tcol\underline{m}}$ with the Jones-Wenzl projector as idempotent. A direct diagrammatic calculation shows that $JW^\Delta$ is indeed an idempotent and thus  $e_{w_0}\mapsto JW^\Delta$ is well-defined. This assignment yields a homomorphism of graded $\DR$-algebras $\varphi^\prime:\DR Q_m \to \mathcal A$. The proof of \cref{thm:isom} is divided into three steps:
\begin{enumerate}[label=(\Roman*)]
\item The map $\varphi^\prime$ is surjective. 
\item $R_m^\mathfrak h \subseteq \ker \varphi^\prime$, so $\varphi^\prime$ induces a surjection $\varphi: \textbf{P}_m \to \mathcal A$.
\item By dimension arguments we can deduce that $\varphi$ is an isomorphism of graded algebras over $\DR$.
\end{enumerate}

\subsection*{Step \MakeUppercase{\romannumeral 1}}\label{subsec:step1}

In order to show that $\varphi^\prime$ is surjective, it suffices to show that for $x,y \in W$ every morphism in $\op{Hom}_\smod^\bullet (B_x,B_y)$ is generated as an algebra by elements of degree $1$. By \cref{thm:dim} we know that between $B_x$ and $B_y$ there exists a morphism of degree $1$ if and only if $|l(x) - l(y)|=1$. Moreover, the homogeneous part of degree $1$ is at most one-dimensional, hence $\op{im}\,\varphi^\prime$ contains all maps of degree $1$ by construction. The double leaves $\DL\DL_{\underline{x},\underline{y}}$ form an $R$-basis for $\op{Hom}_{\mathcal D}(\underline{x},\underline{y})$, thus, using the equivalence in \cref{thm:equiv}, we can deduce that every element in $\op{Hom}_{\bsbim}(B_{\underline{x}},B_{\underline{y}})$ is a sum of maps factoring through various $B_{\underline w}$ such that $\underline{w}$ is a reduced sub-expression of both $\underline{x}$ and $\underline{y}$. We say \emph{modulo lower terms} if we only consider maps which do not factor through $w$ such that $w < x$ and $w < y$. Consequently the next proposition implies the surjectivity of $\varphi^\prime$.
\begin{prop}\label{prop:gen}
For $x\geq y \in W$ the morphism space $\op{Hom}^\bullet_{\sbim}(B_x,B_y)$ modulo lower terms is free of rank $1$ as $R$-left/right module. If $x$ and $y$ are not comparable the morphism space $\op{Hom}^\bullet_{\sbim}(B_x,B_y)= 0$ modulo lower terms.
\end{prop}

\begin{prf}
Before we prove the general case, consider the example $x= \scol\tcol\scol$ (which covers already all the important cases). Then there exist eight decorated $01$-sequences of $x$ with their corresponding light leaves:

\begin{center}
\begin{minipage}{5cm}
\begin{align}
U0\;U0\;U0 \enspace & \rightsquigarrow \enspace
\begin{tikzpicture}[baseline=-0.5ex]
\draw[white] (0.2,-0.3) -- (0.2,0.3);
\draw[sline] (0,-0.3) -- (0,-0.05);
\draw[sdot] (0,-0.05) circle (2pt);
\draw[tline] (0.2,-0.3) -- (0.2,-0.05);
\draw[tdot] (0.2,-0.05) circle (2pt);
\draw[sline] (0.4,-0.3) -- (0.4,-0.05);
\draw[sdot] (0.4,-0.05) circle (2pt);
\end{tikzpicture}  \nonumber \\
U1\;U0\;D1 \enspace & \rightsquigarrow \enspace \,
\begin{tikzpicture}[baseline=-0.5ex];
\draw[tline] (0.2,-0.3) -- (0.2,-0.05);
\draw[tdot] (0.2,-0.05) circle (2pt);
\draw[sline] (0.4,-0.1) arc (0:180:0.2);
\draw[sline] (0,-0.3) -- (0,-0.1);
\draw[sline] (0.4,-0.3) -- (0.4,-0.1);
\end{tikzpicture} \label{eqn:101}\\
U1\;U0\;D0 \enspace & \rightsquigarrow \enspace \,
\begin{tikzpicture}[baseline=-0.5ex];
\draw[tline] (0.2,-0.3) -- (0.2,-0.05);
\draw[tdot] (0.2,-0.05) circle (2pt);
\draw[sline] (0.4,-0.1) arc (0:180:0.2);
\draw[sline] (0,-0.3) -- (0,-0.1);
\draw[sline] (0.4,-0.3) -- (0.4,-0.1);
\draw[sline] (0.2,0.3) -- (0.2,0.1);
\end{tikzpicture} \label{eqn:100}\\
U0\;U0\;U1 \enspace & \rightsquigarrow \enspace 
\begin{tikzpicture}[baseline=-0.5ex]
\draw[sline] (0,-0.3) -- (0,-0.05);
\draw[sdot] (0,-0.05) circle (2pt);
\draw[tline] (0.2,-0.3) -- (0.2,-0.05);
\draw[tdot] (0.2,-0.05) circle (2pt);
\draw[sline] (0.4,-0.3) -- (0.4,0.3);
\end{tikzpicture}\nonumber
\end{align}
\end{minipage}
\begin{minipage}{5cm}
\begin{align*}
U0\;U1\;U0 \enspace & \rightsquigarrow \enspace 
\begin{tikzpicture}[baseline=-0.5ex]
\draw[sline] (0,-0.3) -- (0,-0.05);
\draw[sdot] (0,-0.05) circle (2pt);
\draw[tline] (0.2,-0.3) -- (0.2,0.3);
\draw[sline] (0.4,-0.3) -- (0.4,-0.05);
\draw[sdot] (0.4,-0.05) circle (2pt);
\end{tikzpicture}\\
U1\;U1\;U0 \enspace & \rightsquigarrow \enspace \,
\begin{tikzpicture}[baseline=-0.5ex]
\draw[sline] (0,-0.3) -- (0,0.3);
\draw[tline] (0.2,-0.3) -- (0.2,0.3);
\draw[sline] (0.4,-0.3) -- (0.4,-0.05);
\draw[sdot] (0.4,-0.05) circle (2pt);
\end{tikzpicture}\\
U0\;U1\;U1 \enspace & \rightsquigarrow \enspace 
\begin{tikzpicture}[baseline=-0.5ex]
\draw[sline] (0,-0.3) -- (0,-0.05);
\draw[sdot] (0,-0.05) circle (2pt);
\draw[tline] (0.2,-0.3) -- (0.2,0.3);
\draw[sline] (0.4,-0.3) -- (0.4,0.3);
\end{tikzpicture}\\
U1\;U1\;U1 \enspace & \rightsquigarrow \enspace \,
\begin{tikzpicture}[baseline=-0.5ex]
\draw[sline] (0,-0.3) -- (0,0.3);
\draw[tline] (0.2,-0.3) -- (0.2,0.3);
\draw[sline] (0.4,-0.3) -- (0.4,0.3);
\end{tikzpicture}
\end{align*}
\end{minipage}
\end{center}

The light leaves in (\ref{eqn:101}) and (\ref{eqn:100}) yield the zero map when pre-composed with the idempotent corresponding to $\scol\tcol\scol$ due to death by pitchfork and all other light leaves are generated by degree $1$ maps. Now assume $l(x) \geq 4$ and set $j:=l(y) \leq l(x)=:i$. This means we have to extend the sequences from above with $0$'s and $1$'s on the right side and decorate the new entries properly with either $D$'s or $U$'s. Using successively the fact that the decorated sequences in (\ref{eqn:101}) and (\ref{eqn:100}) yield morphisms which get cancelled by the idempotents, we obtain that the only two possible sequences are

\begin{align}
\underbrace{U0\cdots U0}_{i-j \text{ times}}\;\underbrace{U1\cdots U1}_{j \text{ times}} \enspace & \rightsquigarrow \enspace 
\begin{tikzpicture}[baseline=-0.5ex]
\draw[sline] (0,-0.3) -- (0,-0.05);
\draw[sdot] (0,-0.05) circle (2pt);
\draw[tline] (0.2,-0.3) -- (0.2,-0.05);
\draw[tdot] (0.2,-0.05) circle (2pt);
\draw[bbline] (0.8,-0.3) -- (0.8,-0.05);
\draw[bbdot] (0.8,-0.05) circle (2pt);
\draw[bbline] (1,-0.3) -- (1,0.3);
\draw[bbline] (1.6,-0.3) -- (1.6,0.3);
\draw[dotted, thick] (0.38,-0.2) -- (0.62,-0.2);
\draw[dotted, thick] (1.18,0) -- (1.42,0);
\end{tikzpicture},\\
\underbrace{U0\cdots U0}_{i-j-1 \text{ times}}\;\underbrace{U1\cdots U1}_{i \text{ times}}\;U0 \enspace & \rightsquigarrow \enspace
\begin{tikzpicture}[baseline=-0.5ex]
\draw[sline] (0,-0.3) -- (0,-0.05);
\draw[sdot] (0,-0.05) circle (2pt);
\draw[tline] (0.2,-0.3) -- (0.2,-0.05);
\draw[tdot] (0.2,-0.05) circle (2pt);
\draw[bbline] (0.8,-0.3) -- (0.8,-0.05);
\draw[bbdot] (0.8,-0.05) circle (2pt);
\draw[bbline] (1,-0.3) -- (1,0.3);
\draw[bbline] (1.6,-0.3) -- (1.6,0.3);
\draw[dotted, thick] (0.38,-0.2) -- (0.62,-0.2);
\draw[dotted, thick] (1.18,0) -- (1.42,0);
\draw[bbline] (1.8,-0.3) -- (1.8,-0.05);
\draw[bbdot] (1.8,-0.05) circle (2pt);
\end{tikzpicture}.
\end{align}

Both of those remaining morphisms are generated by degree $1$ elements which proves the first statement. The second statement is clear since every proper sub-expression is strictly smaller than both (this holds dihedral groups but not in arbitrary Coxeter groups).
\end{prf}

\subsection*{Step \MakeUppercase{\romannumeral 2}}
We check that every relation from $R_m^\mathfrak h$ holds for the morphisms between Soergel bimodules (with trivialised right action). For this this we use the following properties.
\begin{enumerate}[label=(\roman*)]
\item Death by pitchfork, \label{pro:1}
\item the polynomial sliding relations in \ref{eqn:poly}, \label{pro:2}
\item non-constant polynomials act as $0$ from the right, \label{pro:3}
\item isotopy invariance, \label{pro:4}
\item two-colour associativity, \label{pro:5}
\item $JW := JW_{m-1}$ is rotation invariant. \label{pro:6}
\end{enumerate}

It is important to note that all calculations and equations in this section are meant as morphisms between indecomposable Soergel bimodules (with trivialised right action). Therefore every morphism represented by a string diagram should be pre-/post-composed with the corresponding idempotent. We omit this composition in order to make the presentation clearer.

Before we start proving the relations one by one, we state two useful results.

\begin{lem}\label{lem:dbldot}
For a Coxeter system of type $I_2(m)$ we have
\[
\begin{tikzpicture}[baseline=-0.5ex,scale=0.8]
\draw[sline] (-0.2,-0.5) -- (-0.2,0.5);
\draw[sline] (0.2,-0.5) -- (0.2,0.5);
\draw[tline] (0,-0.5) -- (0,-0.2);
\draw[tline] (0,0.5) -- (0,0.2);
\draw[tdot] (0,-0.2) circle (2pt);
\draw[tdot] (0,0.2) circle (2pt);
\draw[dotted, thick] (-0.3,0) -- (-0.54,0);
\draw[dotted, thick] (0.3,0) -- (0.54,0);
\draw (0,0.7) -- (0,0.7) node[inner xsep=15pt, rectangle,fill=red!20,draw]{$\textbf{e}_1$};
bf\draw (0,-0.7) -- (0,-0.7) node[inner xsep=15pt, rectangle,fill=red!20,draw]{$\textbf{e}_2$};
\end{tikzpicture}
\enspace = \enspace 0,
\]

where $\textbf{e}_1$ and $\textbf{e}_2$ are idempotents corresponding to indecomposable Soergel bimodules.
\end{lem}

\begin{prf}
Without loss of generality we can assume that the morphism is given by

\begin{align*}
\begin{tikzpicture}[baseline=-0.5ex]
\draw[sline](-0.2,-0.5) -- (-0.2,0.5);
\draw[sline] (0.2,-0.5) -- (0.2,0.5);
\draw[tline] (0,-0.5) -- (0,-0.3);
\draw[tline] (0,0.5) -- (0,0.3);
\draw[tdot] (0,-0.3) circle (2pt);
\draw[tdot] (0,0.3) circle (2pt);
\end{tikzpicture}
\enspace = \enspace 
\begin{tikzpicture}[baseline=-0.5ex]
\draw[sline](-0.4,-0.5) -- (-0.4,0.5);
\draw[sline] (0.4,-0.5) -- (0.4,0.5);
\draw[sline] (-0.4,0) -- (-0.2,0);
\draw[sline] (0.4,0) -- (0.2,0);
\draw[sdot] (-0.2,0) circle (2pt);
\draw[sdot] (0.2,0) circle (2pt);
\draw[tline] (0,-0.5) -- (0,-0.3);
\draw[tline] (0,0.5) -- (0,0.3);
\draw[tdot] (0,-0.3) circle (2pt);
\draw[tdot] (0,0.3) circle (2pt);
\end{tikzpicture}
\enspace \stackrel{\ref{pro:2}}{=} \enspace 
\tiny{\frac{1}{2}}
\left( 
\begin{tikzpicture}[baseline=-0.5ex]
\draw[sline] (-0.5,-0.5) .. controls (-0.5,-0.4).. (-0.4,0);
\draw[sline] (-0.5,0.5) .. controls (-0.5,0.4).. (-0.4,0);
\draw[sline] (0.5,-0.5) .. controls (0.5,-0.4).. (0.4,0);
\draw[sline] (0.5,0.5) .. controls (0.5,0.4).. (0.4,0);
\draw[sline] (-0.4,0) -- (0.4,0);
\draw[tline] (0,-0.5) -- (0,-0.3);
\draw[tline] (0,0.5) -- (0,0.3);
\draw[tdot] (0,-0.3) circle (2pt);
\draw[tdot] (0,0.3) circle (2pt);
\draw[bbline] (-0.25,0.38) -- (-0.25,0.38) node {$ {\tiny \alpha_\tcol}$};
\end{tikzpicture}
\enspace + \enspace
\begin{tikzpicture}[baseline=-0.5ex]
\draw[sline] (-0.5,-0.5) .. controls (-0.5,-0.4).. (-0.4,0);
\draw[sline] (-0.5,0.5) .. controls (-0.5,0.4).. (-0.4,0);
\draw[sline] (0.5,-0.5) .. controls (0.5,-0.4).. (0.4,0);
\draw[sline] (0.5,0.5) .. controls (0.5,0.4).. (0.4,0);
\draw[sline] (-0.4,0) -- (0.4,0);
\draw[tline] (0,-0.5) -- (0,-0.3);
\draw[tline] (0,0.5) -- (0,0.3);
\draw[tdot] (0,-0.3) circle (2pt);
\draw[tdot] (0,0.3) circle (2pt);
\draw[bbline] (-0.25,-0.38) -- (-0.25,-0.38) node {$ {\tiny \alpha_\tcol}$};
\end{tikzpicture}
\right)
\enspace \stackrel{(\ref{eqn:genass})}{=} \enspace 
\tiny{\frac{1}{2}}
\left( 
\begin{tikzpicture}[baseline=-0.5ex]
\draw[sline] (0,-0.1) -- (0,0.1);
\draw[sline] (0,-0.1) .. controls (-0.5,-0.25).. (-0.5,-0.5);
\draw[sline] (0,0.1) .. controls (-0.5,0.25).. (-0.5,0.5);
\draw[sline] (0,-0.1) .. controls (0.5,-0.25).. (0.5,-0.5);
\draw[sline] (0,0.1) .. controls (0.5,0.25).. (0.5,0.5);
\draw[tline] (0,-0.5) -- (0,-0.3);
\draw[tline] (0,0.5) -- (0,0.3);
\draw[tdot] (0,-0.3) circle (2pt);
\draw[tdot] (0,0.3) circle (2pt);
\draw[bbline] (-0.25,0.38) -- (-0.25,0.38) node {$ {\tiny \alpha_\tcol}$};
\end{tikzpicture}
\enspace + \enspace
\begin{tikzpicture}[baseline=-0.5ex]
\draw[sline] (0,-0.1) -- (0,0.1);
\draw[sline] (0,-0.1) .. controls (-0.5,-0.25).. (-0.5,-0.5);
\draw[sline] (0,0.1) .. controls (-0.5,0.25).. (-0.5,0.5);
\draw[sline] (0,-0.1) .. controls (0.5,-0.25).. (0.5,-0.5);
\draw[sline] (0,0.1) .. controls (0.5,0.25).. (0.5,0.5);
\draw[tline] (0,-0.5) -- (0,-0.3);
\draw[tline] (0,0.5) -- (0,0.3);
\draw[tdot] (0,-0.3) circle (2pt);
\draw[tdot] (0,0.3) circle (2pt);
\draw[bbline] (-0.25,-0.38) -- (-0.25,-0.38) node {$ {\tiny \alpha_\tcol}$};
\end{tikzpicture}
\right)
\enspace \stackrel{\ref{pro:1}}{=} \enspace 0.
\end{align*}
The general case follows easily with the above considerations.
\end{prf}

\begin{prop}\label{prop:polyslide}
As a morphism $B_{{}_\scol i} \to B_{{}_\scol i}$ in the category of Soergel modules the following holds for $i\geq 2$:
\[  \begin{tikzpicture}[baseline=-0.5ex]
\draw[bbline] (-0,-0.4) -- (-0,0.4);
\draw[bbline] (-0.2,-0.4) -- (-0.2,0.4);
\draw[sline] (-0.8,-0.4) -- (-0.8,0.4);
\draw[tline] (-1,-0.15) -- (-1,0.15);
\draw[tdot] (-1,-0.15) circle (2pt);
\draw[tdot] (-1,0.15) circle (2pt);
\draw[dotted, thick] (-0.38,0) -- (-0.62,0);
\end{tikzpicture} 
\enspace = \enspace \textstyle{-[2]_\zeta} \enspace
\begin{tikzpicture}[baseline=-0.5ex]
\draw[sline] (0,-0.4) -- (0,-0.2);
\draw[sline] (0,0.4) -- (0,0.2);
\draw[sdot] (0,-0.2) circle (2pt);
\draw[sdot] (0,0.2) circle (2pt);
\draw[tline] (0.2,-0.4) -- (0.2,0.4);
\draw[bbline] (0.8,-0.4) -- (0.8,0.4);
\draw[dotted, thick] (0.38,0) -- (0.62,0);
\end{tikzpicture}
\enspace + \enspace \textstyle{[i-1]_\zeta - [i+1]_\zeta}\enspace
\begin{tikzpicture}[baseline=-0.5ex]
\draw[sline] (0,-0.4) -- (0,0.4);
\draw[tline] (0.2,-0.4) -- (0.2,0.4);
\draw[bbline] (0.8,-0.4) -- (0.8,-0.2);
\draw[bbdot] (0.8,-0.2) circle (2pt);
\draw[bbline] (0.8,0.4) -- (0.8,0.2);
\draw[bbdot] (0.8,0.2) circle (2pt);
\draw[dotted, thick] (0.38,0) -- (0.62,0);
\end{tikzpicture}.
\]
\end{prop}

\begin{prf}
By property \ref{pro:2} we can slide the polynomial $\alpha_\tcol$ successively through the strings and we obtain for some $\lambda_n$:
\begin{align*}
\begin{tikzpicture}[baseline=-0.5ex]
\draw[bbline] (-0,-0.4) -- (-0,0.4);
\draw[bbline] (-0.2,-0.4) -- (-0.2,0.4);
\draw[sline] (-0.8,-0.4) -- (-0.8,0.4);
\draw[tline] (-1,-0.15) -- (-1,0.15);
\draw[tdot] (-1,-0.15) circle (2pt);
\draw[tdot] (-1,0.15) circle (2pt);
\draw[dotted, thick] (-0.38,0) -- (-0.62,0);
\end{tikzpicture}
\enspace & = \enspace -[2]_\zeta \enspace 
\begin{tikzpicture}[baseline=-0.5ex]
\draw[sline] (0,-0.4) -- (0,-0.2);
\draw[sline] (0,0.4) -- (0,0.2);
\draw[sdot] (0,-0.2) circle (2pt);
\draw[sdot] (0,0.2) circle (2pt);
\draw[tline] (0.2,-0.4) -- (0.2,0.4);
\draw[bbline] (0.8,-0.4) -- (0.8,0.4);
\draw[dotted, thick] (0.38,0) -- (0.62,0);
\end{tikzpicture}
\enspace + \enspace \sum_{n=2}^{i-1} \lambda_n \;
\begin{tikzpicture}[baseline=-0.5ex,bbline]
\draw[sline] (-0.64,-0.4) -- (-0.64,0.4);
\draw (0.64,-0.4) -- (0.64,0.4);
\draw (-0.2,-0.4) -- (-0.2,0.4);
\draw (0.2,-0.4) -- (0.2,0.4);
\draw (0,-0.4) -- (0,-0.2);
\draw (0,0.4) -- (0,0.2);
\draw[bbdot] (0,-0.2) circle (2pt);
\draw[bbdot] (0,0.2) circle (2pt);
\draw[dotted, thick] (-0.3,0) -- (-0.54,0);
\draw[dotted, thick] (0.3,0) -- (0.54,0);
\draw (0,0.53) -- (0,0.53) node {$\footnotesize n$};
\draw (0,-0.53) -- (0,-0.53) node {\phantom{*}};
\end{tikzpicture}
\enspace + \enspace
\begin{tikzpicture}[baseline=-0.5ex,bbline]
\draw[sline] (-0.34,-0.4) -- (-0.34,0.4);
\draw (0.64,-0.4) -- (0.64,0.4);
\draw (0.45,0) -- (0.45,0) node {$\small x$};
\draw (0.26,-0.4) -- (0.26,0.4);
\draw[dotted, thick] (-0.16,0) -- (0.08,0);
\end{tikzpicture}
\intertext{where $x = (i-1)_\scol(\alpha_\tcol)$ is of positive degree. By \cref{lem:dbldot} we have:}
\enspace & = \enspace -[2]_\zeta \enspace
\begin{tikzpicture}[baseline=-0.5ex]
\draw[sline] (0,-0.4) -- (0,-0.2);
\draw[sline] (0,0.4) -- (0,0.2);
\draw[sdot] (0,-0.2) circle (2pt);
\draw[sdot] (0,0.2) circle (2pt);
\draw[tline] (0.2,-0.4) -- (0.2,0.4);
\draw[bbline] (0.8,-0.4) -- (0.8,0.4);
\draw[dotted, thick] (0.38,0) -- (0.62,0);
\end{tikzpicture}
\enspace + \enspace
\begin{tikzpicture}[baseline=-0.5ex,bbline]
\draw[sline] (-0.34,-0.4) -- (-0.34,0.4);
\draw (0.64,-0.4) -- (0.64,0.4);
\draw (0.45,0) -- (0.45,0) node {$\small x$};
\draw (0.26,-0.4) -- (0.26,0.4);
\draw[dotted, thick] (-0.16,0) -- (0.08,0);
\end{tikzpicture}
\end{align*}

From now on assume that $i$ is odd:
\begin{align*}
\begin{tikzpicture}[baseline=-0.5ex]
\draw[sline] (-0,-0.4) -- (-0,0.4);
\draw[tline] (-0.2,-0.4) -- (-0.2,0.4);
\draw[sline] (-0.8,-0.4) -- (-0.8,0.4);
\draw[tline] (-1,-0.15) -- (-1,0.15);
\draw[tdot] (-1,-0.15) circle (2pt);
\draw[tdot] (-1,0.15) circle (2pt);
\draw[dotted, thick] (-0.38,0) -- (-0.62,0);
\end{tikzpicture}
\enspace & = \enspace -[2]_\zeta \enspace
\begin{tikzpicture}[baseline=-0.5ex]
\draw[sline] (0,-0.4) -- (0,-0.2);
\draw[sline] (0,0.4) -- (0,0.2);
\draw[sdot] (0,-0.2) circle (2pt);
\draw[sdot] (0,0.2) circle (2pt);
\draw[tline] (0.2,-0.4) -- (0.2,0.4);
\draw[sline] (0.8,-0.4) -- (0.8,0.4);
\draw[dotted, thick] (0.38,0) -- (0.62,0);
\end{tikzpicture}
\enspace + \enspace
\begin{tikzpicture}[baseline=-0.5ex,bbline]
\draw[sline] (-0.34,-0.4) -- (-0.34,0.4);
\draw[sline] (0.64,-0.4) -- (0.64,0.4);
\draw (0.45,0) -- (0.45,0) node {$\small x$};
\draw[tline] (0.26,-0.4) -- (0.26,0.4);
\draw[dotted, thick] (-0.16,0) -- (0.08,0);
\end{tikzpicture}\\[0.5cm]
\enspace & \stackrel{\ref{pro:2}}{=} \enspace -[2]_\zeta \enspace
\begin{tikzpicture}[baseline=-0.5ex]
\draw[sline] (0,-0.4) -- (0,-0.2);
\draw[sline] (0,0.4) -- (0,0.2);
\draw[sdot] (0,-0.2) circle (2pt);
\draw[sdot] (0,0.2) circle (2pt);
\draw[tline] (0.2,-0.4) -- (0.2,0.4);
\draw[sline] (0.8,-0.4) -- (0.8,0.4);
\draw[dotted, thick] (0.38,0) -- (0.62,0);
\end{tikzpicture}
\enspace + \enspace
\begin{tikzpicture}[baseline=-0.5ex]
\draw[sline] (-0.34,-0.4) -- (-0.34,0.4);
\draw[sline] (0.46,-0.4) -- (0.46,0.4);
\draw (0.84,0) -- (0.84,0) node {$\small \scol(x)$};
\draw[tline] (0.26,-0.4) -- (0.26,0.4);
\draw[dotted, thick] (-0.16,0) -- (0.08,0);
\end{tikzpicture}
\enspace + \enspace \partial_\scol(x) \enspace
\begin{tikzpicture}[baseline=-0.5ex]
\draw[sline] (-0.34,-0.4) -- (-0.34,0.4);
\draw[sline] (0.46,-0.4) -- (0.46,-0.2);
\draw[sline] (0.46,0.4) -- (0.46,0.2);
\draw[sdot] (0.46,-0.2) circle (2pt);
\draw[sdot] (0.46,0.2) circle (2pt);
\draw[tline] (0.26,-0.4) -- (0.26,0.4);
\draw[dotted, thick] (-0.16,0) -- (0.08,0);
\end{tikzpicture}\\[0.5cm]
\enspace & \stackrel{\ref{pro:3}}{=} \enspace -[2]_\zeta \enspace
\begin{tikzpicture}[baseline=-0.5ex]
\draw[sline] (0,-0.4) -- (0,-0.2);
\draw[sline] (0,0.4) -- (0,0.2);
\draw[sdot] (0,-0.2) circle (2pt);
\draw[sdot] (0,0.2) circle (2pt);
\draw[tline] (0.2,-0.4) -- (0.2,0.4);
\draw[sline] (0.8,-0.4) -- (0.8,0.4);
\draw[dotted, thick] (0.38,0) -- (0.62,0);
\end{tikzpicture}
\enspace + \enspace \partial_\scol(x) \enspace
\begin{tikzpicture}[baseline=-0.5ex]
\draw[sline] (-0.34,-0.4) -- (-0.34,0.4);
\draw[sline] (0.46,-0.4) -- (0.46,-0.2);
\draw[sline] (0.46,0.4) -- (0.46,0.2);
\draw[sdot] (0.46,-0.2) circle (2pt);
\draw[sdot] (0.46,0.2) circle (2pt);
\draw[tline] (0.26,-0.4) -- (0.26,0.4);
\draw[dotted, thick] (-0.16,0) -- (0.08,0);
\end{tikzpicture},
\end{align*}
where we used that $s$ acts grading preserving and thus $s(x)$ acts trivially from the right since $x$ was of positive degree. Recall that $i$ is odd and hence by \cref{lem:actionroots} we have:
\begin{align*}
\partial_\scol(x) \enspace = \enspace
\partial_\scol \left( [i]_\zeta \alpha_\tcol + [i-1]_\zeta \alpha_\scol\right)
\enspace = \enspace 
-[2]_\zeta [i]_\zeta + 2[i-1]_\zeta 
\enspace \stackrel{\ref{eqn:q1}}{=} \enspace
[i-1]_\zeta - [i+1]_\zeta.
\end{align*}
The other case ($i$ even) can be treated similarly using $\partial_\tcol$ instead of $\partial_\scol$.
\end{prf}

\subsubsection{Relation \ref{eqn:1} and \ref{eqn:2}}
Clearly we have:
\begin{alignat*}{3}
\begin{tikzpicture}[baseline=-0.5ex]
\draw[sline] (1,-0.15)--(1,0.15);
\draw[sdot] (1,-0.15) circle (2pt);
\draw[sdot] (1,0.15) circle (2pt);
\end{tikzpicture}
\enspace \stackrel{\ref{pro:3}}{=} &\enspace 0 \enspace& \stackrel{\ref{pro:3}}{=} \enspace
\begin{tikzpicture}[baseline=-0.5ex]
\draw[tline] (0.8,-0.4) -- (0.8,0.4);
\draw[sline] (1,-0.15)--(1,0.15);
\draw[sdot] (1,-0.15) circle (2pt);
\draw[sdot] (1,0.15) circle (2pt);
\end{tikzpicture}.
\end{alignat*}

\subsubsection{Relation \ref{eqn:3} and \ref{eqn:5}}
The desired relations for $i\in \{2,\ldots,m-2\}$ are direct consequences of the following lemma:
\begin{lem}\label{lem:jw}
For $i\in \{2,\ldots,m-1\}$ the following holds:
\begin{align*}
\begin{tikzpicture}[baseline=-0.5ex]
\draw[sline] (-0.5,-0.4) -- (-0.5,0.4);
\draw[tline] (-0.25,-0.4) -- (-0.25,0.4);
\draw[bbline] (0.25,-0.4) -- (0.25,0.4);
\draw[bbdot] (0.55,0) circle (2pt);
\draw[dotted, thick] (-0.12,0.35) -- (0.12,0.35);
\draw[dotted, thick] (-0.12,-0.35) -- (0.12,-0.35);
\draw[bbline] (0.55,0) -- (0,0);
\draw[sline] (-0.5,0) -- (0,0);
\draw (0,0) -- (0,0) node[fill=red!20, draw] {$\scriptstyle{JW_i}$};
\end{tikzpicture}
\enspace & =  \enspace
\textstyle{\frac{[i-1]_\zeta}{[i]_\zeta}} \enspace
\begin{tikzpicture}[baseline=-0.5ex]
\draw[sline] (0,-0.4) -- (0,0.4);
\draw[tline] (0.2,-0.4) -- (0.2,0.4);
\draw[bbline] (0.8,-0.4) -- (0.8,-0.2);
\draw[bbdot] (0.8,-0.2) circle (2pt);
\draw[bbline] (0.8,0.4) -- (0.8,0.2);
\draw[bbdot] (0.8,0.2) circle (2pt);
\draw[dotted, thick] (0.38,0) -- (0.62,0);
\end{tikzpicture},\\[0.5cm]
\begin{tikzpicture}[baseline=-0.5ex]
\draw[tdot] (-0.55,0) circle (2pt);
\draw[sline] (-0.25,-0.4) -- (-0.25,0.4);
\draw[bbline] (0.25,-0.4) -- (0.25,0.4);
\draw[bbline] (0.5,-0.4) -- (0.5,0.4);
\draw[dotted, thick] (-0.12,0.35) -- (0.12,0.35);
\draw[dotted, thick] (-0.12,-0.35) -- (0.12,-0.35);
\draw[bbline] (0.5,0) -- (0,0);
\draw[tline] (-0.55,0) -- (0,0);
\draw (0,0) -- (0,0) node[fill=red!20, draw] {$\scriptstyle{JW_i}$};
\end{tikzpicture}
\enspace & =  \enspace \textstyle{-\frac{[i-1]_\zeta}{[i]_\zeta}}\enspace
\begin{tikzpicture}[baseline=-0.5ex]
\draw[sline] (0,-0.4) -- (0,-0.2);
\draw[sline] (0,0.4) -- (0,0.2);
\draw[sdot] (0,-0.2) circle (2pt);
\draw[sdot] (0,0.2) circle (2pt);
\draw[tline] (0.2,-0.4) -- (0.2,0.4);
\draw[bbline] (0.8,-0.4) -- (0.8,0.4);
\draw[dotted, thick] (0.38,0) -- (0.62,0);
\end{tikzpicture}
\enspace + \enspace \textstyle{([i-1]_\zeta -[i+1]_\zeta)}\enspace
\begin{tikzpicture}[baseline=-0.5ex]
\draw[sline] (0,-0.4) -- (0,0.4);
\draw[tline] (0.2,-0.4) -- (0.2,0.4);
\draw[bbline] (0.8,-0.4) -- (0.8,-0.2);
\draw[bbdot] (0.8,-0.2) circle (2pt);
\draw[bbline] (0.8,0.4) -- (0.8,0.2);
\draw[bbdot] (0.8,0.2) circle (2pt);
\draw[dotted, thick] (0.38,0) -- (0.62,0);
\end{tikzpicture}.
\intertext{In particular:}
\begin{tikzpicture}[baseline=-0.5ex]
\draw[sline] (-0.5,-0.4) -- (-0.5,0.4);
\draw[tline] (-0.25,-0.4) -- (-0.25,0.4);
\draw[bbline] (0.25,-0.4) -- (0.25,0.4);
\draw[bbdot] (0.5,0) circle (2pt);
\draw[dotted, thick] (-0.12,0.35) -- (0.12,0.35);
\draw[dotted, thick] (-0.12,-0.35) -- (0.12,-0.35);
\draw[bbline] (0.5,0) -- (0,0);
\draw[sline] (-0.5,0) -- (0,0);
\draw (0,0) -- (0,0) node[fill=red!20, draw] {$\scriptstyle{JW}$};
\end{tikzpicture}
\enspace & = \enspace \textstyle{[2]_\zeta} \enspace
\begin{tikzpicture}[baseline=-0.5ex]
\draw[sline] (0,-0.4) -- (0,0.4);
\draw[tline] (0.2,-0.4) -- (0.2,0.4);
\draw[bbline] (0.8,-0.4) -- (0.8,-0.2);
\draw[bbdot] (0.8,-0.2) circle (2pt);
\draw[bbline] (0.8,0.4) -- (0.8,0.2);
\draw[bbdot] (0.8,0.2) circle (2pt);
\draw[dotted, thick] (0.38,0) -- (0.62,0);
\end{tikzpicture}\enspace = \enspace
\begin{tikzpicture}[baseline=-0.5ex]
\draw[tdot] (-0.5,0) circle (2pt);
\draw[sline] (-0.25,-0.4) -- (-0.25,0.4);
\draw[bbline] (0.25,-0.4) -- (0.25,0.4);
\draw[bbline] (0.5,-0.4) -- (0.5,0.4);
\draw[dotted, thick] (-0.12,0.35) -- (0.12,0.35);
\draw[dotted, thick] (-0.12,-0.35) -- (0.12,-0.35);
\draw[bbline] (0.5,0) -- (0,0);
\draw[tline] (-0.5,0) -- (0,0);
\draw (0,0) -- (0,0) node[fill=red!20, draw] {$\scriptstyle{JW}$};
\end{tikzpicture}.
\end{align*}
\end{lem}

\begin{prf}
Recall the recursive formula for $JW_i$ (see (\ref{thm:jw})) in the Temperley-Lieb algebra setting in which we specialise $q$ to $\zeta= e^{2\pi i/2m}$. Transforming this Jones-Wenzl projector into a morphism in $\mathcal D$ and putting a trivalent vertex on the left side and a dot on the right side yields a morphism between the Bott-Samelson bimodules. If we pre-/post-compose with the idempotents corresponding to the indecomposable Soergel bimodules every summand is cancelled by property (\ref{pro:1}) except the following two summands (the lines from the last inductive step are dashed):

\begin{alignat*}{2}
\textstyle{\frac{[i-1]_\zeta}{[i]_\zeta}}\enspace
\begin{tikzpicture}[baseline=-0.5ex,bline]
\draw (0,-0.4) to (0,0.4);
\draw (0.6,-0.4) to (0.6,0.4);
\draw[dotted,thick] (0.42,0) to (0.18,0);
\draw (1.2,-0.4) arc (0:180:0.2);
\draw (0.8,0.4) arc (180:360:0.2);
\end{tikzpicture}
\enspace & = \enspace \textstyle{\frac{[i-1]_\zeta}{[i]_\zeta}}\enspace
\begin{tikzpicture}[baseline=-0.5ex,bline]
\draw (-0.2,-0.4) to (-0.2,0.4);
\draw (-0.6,-0.4) to (-0.6,0.4);
\draw (0.3,-0.4) to (0.3,0);
\draw[dotted,thick] (-0.52,0) to (-0.28,0);
\draw[densely dashed] (0,0.4) arc (180:360:0.2);
\draw[densely dashed] (0.7,0) arc (0:180:0.2);
\draw[densely dashed] (0.7,0) -- (0.7,-0.4);
\end{tikzpicture};\qquad & \qquad
\begin{tikzpicture}[baseline=-0.5ex,bline]
\draw (-0.6,-0.4) to (-0.6,0.4);
\draw (0,-0.4) to (0,0.4);
\draw[dotted,thick] (-0.42,0) to (-0.18,0);
\end{tikzpicture}
\enspace & =  \enspace
\begin{tikzpicture}[baseline=-0.5ex,bline]
\draw (-0.2,-0.4) to (-0.2,0.4);
\draw (-0.6,-0.4) to (-0.6,0.4);
\draw[densely dashed] (0,-0.4) to (0,0.4);
\draw[dotted,thick] (-0.52,0) to (-0.28,0);
\end{tikzpicture}.
\end{alignat*}

The first summand is a cup attached to the identity of $JW_{i-1}$ at position $i-1$ and thus its scalar is $\textstyle{\frac{[i-1]_\zeta}{[i]_\zeta}}$. The second summand corresponds to the identity and therefore its scalar is $1$. Hence in $\mathcal D$ we have

\[
\begin{tikzpicture}[baseline=-0.5ex]
\draw[sline] (-0.5,-0.4) -- (-0.5,0.4);
\draw[tline] (-0.25,-0.4) -- (-0.25,0.4);
\draw[bbline] (0.25,-0.4) -- (0.25,0.4);
\draw[bbdot] (0.55,0) circle (2pt);
\draw[dotted, thick] (-0.12,0.35) -- (0.12,0.35);
\draw[dotted, thick] (-0.12,-0.35) -- (0.12,-0.35);
\draw[bbline] (0.55,0) -- (0,0);
\draw[sline] (-0.5,0) -- (0,0);
\draw (0,0) -- (0,0) node[fill=red!20, draw] {$\scriptstyle{JW_i}$};
\end{tikzpicture}
\enspace =  \enspace  \textstyle{\frac{[i-1]_\zeta}{[i]_\zeta}} \enspace
\begin{tikzpicture}[baseline=-0.5ex]
\draw[sline] (-0.5,-0.4) -- (-0.5,0.4);
\draw[sline] (-0.5,0) -- (-0.35,0);
\draw[sdot] (-0.35,0) circle (2pt);
\draw[tline] (-0.2,-0.4) -- (-0.2,0.4);
\draw[bbline] (0.4,-0.4) -- (0.4,-0.2);
\draw[bbline] (0.4,0.2) -- (0.4,0.4);
\draw[bbdot] (0.65,0) circle (2pt);
\draw[bbdot] (0.4,0.2) circle (2pt);
\draw[bbdot] (0.4,-0.2) circle (2pt);
\draw[dotted, thick] (-0.095,0) -- (0.145,0);
\draw[bbline] (0.65,0) -- (0.25,0);
\draw[bbline] (0.25,-0.4) -- (0.25,0.4);
\end{tikzpicture}
\enspace +  \enspace
\underbrace{\begin{tikzpicture}[baseline=-0.5ex]
\draw[sline] (-0.5,-0.4) -- (-0.5,0.4);
\draw[sline] (-0.5,0) -- (-0.35,0);
\draw[sdot] (-0.35,0) circle (2pt);
\draw[tline] (-0.2,-0.4) -- (-0.2,0.4);
\draw[bbdot] (0.7,0) circle (2pt);
\draw[bbdot] (0.45,0) circle (2pt);
\draw[dotted, thick] (-0.095,0) -- (0.145,0);
\draw[bbline] (0.7,0) -- (0.45,0);
\draw[bbline] (0.25,-0.4) -- (0.25,0.4);
\end{tikzpicture}}_{=0 \text{ by \ref{pro:3}}}\\[0.5cm]
\enspace =  \enspace \textstyle{\frac{[i-1]_\zeta}{[i]_\zeta}}  \enspace
\begin{tikzpicture}[baseline=-0.5ex]
\draw[sline] (0,-0.4) -- (0,0.4);
\draw[tline] (0.2,-0.4) -- (0.2,0.4);
\draw[bbline] (0.8,-0.4) -- (0.8,0.4);
\draw[bbline] (1,-0.4) -- (1,-0.2);
\draw[bbdot] (1,-0.2) circle (2pt);
\draw[bbline] (1,0.4) -- (1,0.2);
\draw[bbdot] (1,0.2) circle (2pt);
\draw[dotted, thick] (0.38,0) -- (0.62,0);
\end{tikzpicture}.
\]

Similarly, we have

\begin{align*}
\begin{tikzpicture}[baseline=-0.5ex]
\draw[tdot] (-0.55,0) circle (2pt);
\draw[sline] (-0.25,-0.4) -- (-0.25,0.4);
\draw[bbline] (0.25,-0.4) -- (0.25,0.4);
\draw[bbline] (0.5,-0.4) -- (0.5,0.4);
\draw[dotted, thick] (-0.12,0.35) -- (0.12,0.35);
\draw[dotted, thick] (-0.12,-0.35) -- (0.12,-0.35);
\draw[bbline] (0.5,0) -- (0,0);
\draw[tline] (-0.55,0) -- (0,0);
\draw (0,0) -- (0,0) node[fill=red!20, draw] {$\scriptstyle{JW_i}$};
\end{tikzpicture}
\enspace = & \enspace  \textstyle{\frac{[i-1]_\zeta}{[i]_\zeta}}  \enspace 
\begin{tikzpicture}[baseline=-0.5ex]
\draw[bbline] (-0,-0.4) -- (-0,0.4);
\draw[bbline] (-0.2,-0.4) -- (-0.2,0.4);
\draw[tline] (-0.8,-0.4) -- (-0.8,0.4);
\draw[sline] (-1,-0.4) -- (-1,-0.2);
\draw[sdot] (-1,-0.2) circle (2pt);
\draw[sline] (-1,0.4) -- (-1,0.2);
\draw[sdot] (-1,0.2) circle (2pt);
\draw[dotted, thick] (-0.38,0) -- (-0.62,0);
\end{tikzpicture}
\enspace + \enspace 
\begin{tikzpicture}[baseline=-0.5ex]
\draw[bbline] (-0,-0.4) -- (-0,0.4);
\draw[bbline] (-0.2,-0.4) -- (-0.2,0.4);
\draw[sline] (-0.8,-0.4) -- (-0.8,0.4);
\draw[tline] (-1,-0.15) -- (-1,0.15);
\draw[tdot] (-1,-0.15) circle (2pt);
\draw[tdot] (-1,0.15) circle (2pt);
\draw[dotted, thick] (-0.38,0) -- (-0.62,0);
\end{tikzpicture}.
\end{align*}

The scalar in front of the first summand can be computed with the recursive formula for $JW_i$:
\begin{align}
\lambda_i \enspace \begin{tikzpicture}[baseline=-0.5ex,bline]
\draw (0,-0.4) to (0,0.4);
\draw (-0.6,-0.4) to (-0.6,0.4);
\draw[dotted,thick] (-0.42,0) to (-0.18,0);
\draw (-0.8,-0.4) arc (0:180:0.2);
\draw (-1.2,0.4) arc (180:360:0.2);
\end{tikzpicture}
\enspace & = \enspace \lambda_{i-1} \enspace
\begin{tikzpicture}[baseline=-0.5ex,bline]
\draw[densely dashed] (0,-0.4) to (0,0.4);
\draw (-0.2,-0.4) to (-0.2,0.4);
\draw (-0.6,-0.4) to (-0.6,0.4);
\draw[dotted,thick] (-0.52,0) to (-0.28,0);
\draw (-0.8,-0.4) arc (0:180:0.2);
\draw (-1.2,0.4) arc (180:360:0.2);
\end{tikzpicture}
\enspace + \enspace I_i\mu_{i-1} \enspace
\begin{tikzpicture}[baseline=-0.5ex,bline]
\draw (-0.2,-0.4) to (-0.2,0.4);
\draw (-0.6,-0.4) to (-0.6,0.4);
\draw[dotted,thick] (-0.52,0) to (-0.28,0);
\draw (-0.8,-0.4) arc (0:180:0.2);
\draw[densely dashed] (-1.2,0.4) arc (180:360:0.2);
\draw (0,0.4) arc (180:360:0.2);
\draw[densely dashed] (0.8,0.4) arc (0:180:0.2);
\draw[densely dashed] (0.8,0.4) -- (0.8,-0.4);
\end{tikzpicture},\label{eqn:lambdaind}
\end{align}

where $\lambda_i$ is the coefficient of $\begin{tikzpicture}[baseline=-0.5ex,bline,scale=0.5]
\draw (0,-0.4) to (0,0.4);
\draw (-0.6,-0.4) to (-0.6,0.4);
\draw[dotted,thick] (-0.42,0) to (-0.18,0);
\draw (-0.8,-0.4) arc (0:180:0.2);
\draw (-1.2,0.4) arc (180:360:0.2);
\end{tikzpicture}$ in $JW_i$, $\mu_i$ is the coefficient of $\begin{tikzpicture}[baseline=-0.5ex,bline,scale=0.5]
\draw (-0.6,0.4) .. controls (-0.6,0.1) and (0,-0.1) .. (0,-0.4);
\draw (-1.2,0.4) .. controls (-1.2,0.1) and (-0.6,-0.1) .. (-0.6,-0.4);
\draw[dotted,thick] (-0.72,0) to (-0.48,0);
\draw (-0.4,0.4) arc (180:360:0.2);
\draw (-0.8,-0.4) arc (0:180:0.2);
\end{tikzpicture}$ in $JW_i$ and $I_{i}$ is the coefficient in the sum from the induction step.

Inductively it can be shown that $\mu_i = \frac{1}{[i]_\zeta}$ using the formula for $JW_i$. Clearly $I_i=\frac{1}{[i]}_\zeta$ and hence by another induction it follows that $\lambda_i=\frac{[i-1]_\zeta}{[i]_\zeta}$ using
\begin{align*}
\lambda_i  &=  \lambda_{i-1} + I_i\mu_{i-1}= \frac{[(i-1)-1]_\zeta[(i-1)+1]_\zeta +1}{[i]_\zeta [i-1]_\zeta} \stackrel{\ref{eqn:q5}}{=} \frac{[i-1]^2_\zeta}{[i]_\zeta [i-1]_\zeta}= \frac{[i-1]_\zeta}{[i]_\zeta}.
\end{align*}

The second summand in \ref{eqn:lambdaind} arises again from the identity and hence its scalar is $1$. With \cref{prop:polyslide} we obtain
\begin{align*}
\begin{tikzpicture}[baseline=-0.5ex]
\draw[tdot] (-0.55,0) circle (2pt);
\draw[sline] (-0.25,-0.4) -- (-0.25,0.4);
\draw[bbline] (0.25,-0.4) -- (0.25,0.4);
\draw[bbline] (0.5,-0.4) -- (0.5,0.4);
\draw[dotted, thick] (-0.12,0.35) -- (0.12,0.35);
\draw[dotted, thick] (-0.12,-0.35) -- (0.12,-0.35);
\draw[bbline] (0.5,0) -- (0,0);
\draw[tline] (-0.55,0) -- (0,0);
\draw (0,0) -- (0,0) node[fill=red!20, draw] {$\scriptstyle{JW_i}$};
\end{tikzpicture}
\enspace = & \enspace \textstyle{\left(\frac{[i-1]_\zeta}{[i]_\zeta} -[2]_\zeta \right)} \enspace
\begin{tikzpicture}[baseline=-0.5ex]
\draw[sline] (0,-0.4) -- (0,-0.2);
\draw[sline] (0,0.4) -- (0,0.2);
\draw[sdot] (0,-0.2) circle (2pt);
\draw[sdot] (0,0.2) circle (2pt);
\draw[tline] (0.2,-0.4) -- (0.2,0.4);
\draw[bbline] (0.8,-0.4) -- (0.8,0.4);
\draw[dotted, thick] (0.38,0) -- (0.62,0);
\end{tikzpicture}
\enspace + \enspace \textstyle{([i-1]_\zeta - [i+1]_\zeta)}\enspace
\begin{tikzpicture}[baseline=-0.5ex]
\draw[sline] (0,-0.4) -- (0,0.4);
\draw[tline] (0.2,-0.4) -- (0.2,0.4);
\draw[bbline] (0.8,-0.4) -- (0.8,-0.2);
\draw[bbdot] (0.8,-0.2) circle (2pt);
\draw[bbline] (0.8,0.4) -- (0.8,0.2);
\draw[bbdot] (0.8,0.2) circle (2pt);
\draw[dotted, thick] (0.38,0) -- (0.62,0);
\end{tikzpicture}.
\end{align*}

where the scalar in front of the first summand is
\begin{align*} \left(\frac{[i-1]_\zeta}{[i]_\zeta} -[2]_\zeta \right)
= \frac{[i-1]_\zeta - [2]_\zeta [i]_\zeta}{[i]_\zeta}
\stackrel{\ref{eqn:q1}}{=} \frac{[i-1]_\zeta - [i+1]_\zeta - [i-1]_\zeta}{[i]_\zeta}
= - \frac{[i+1]_\zeta}{[i]_\zeta}. \end{align*}
\end{prf}

For $i=m-1$ the Soergel bimodule is embedded diagonally and therefore we have

\begin{align*}
\frac{1}{2}\begin{pmatrix}
\begin{tikzpicture}[baseline=-0.5ex]
\draw[sline] (0,-0.4) -- (0,0.4);
\draw[tline] (0.2,-0.4) -- (0.2,0.4);
\draw[sline] (0.4,-0.4) -- (0.4,0.4);
\draw[dotted, thick] (0.58,0) -- (0.82,0);
\draw[bbline] (1,0)--(1,-0.4);
\draw[bbdot] (1,0) circle (2pt);
\end{tikzpicture} &
\begin{tikzpicture}[baseline=-0.5ex]
\draw[tline] (0,0)--(0,-0.4);
\draw[tdot] (0,0) circle (2pt);
\draw[sline] (0.2,-0.4) -- (0.2,0.4);
\draw[tline] (0.4,-0.4) -- (0.4,0.4);
\draw[dotted, thick] (0.58,0) -- (0.82,0);
\draw[bbline] (1,-0.4) -- (1,0.4);
\end{tikzpicture}
\end{pmatrix}
& \frac{1}{2}
\begin{pmatrix}
\begin{tikzpicture}[baseline=-0.5ex]
\draw[sline] (-0.5,-0.4) -- (-0.5,0.4);
\draw[tline] (-0.25,-0.4) -- (-0.25,0.4);
\draw[bbline] (0.25,-0.4) -- (0.25,0.4);
\draw[bbline] (0.5,-0.4) -- (0.5,0.4);
\draw[dotted, thick] (-0.12,0.35) -- (0.12,0.35);
\draw[dotted, thick] (-0.12,-0.35) -- (0.12,-0.35);
\draw[bbline] (0.5,0) -- (0,0);
\draw[sline] (-0.5,0) -- (0,0);
\draw (0,0) -- (0,0) node[fill=red!20, draw] {$\scriptstyle{JW}$};
\end{tikzpicture}
 & \begin{tikzpicture}[baseline=-0.5ex]
\draw[tline] (-0.5,-0.4) .. controls (-0.4,0) .. (0,0);
\draw[sline] (-0.25,-0.4) .. controls (-0.2,0).. (0,0);
\draw[bbline] (0.5,-0.4) .. controls (0.4,0) .. (0,0);
\draw[sline] (-0.5,0.4) .. controls (-0.4,0) .. (0,0);
\draw[tline] (-0.25,0.4) .. controls (-0.2,0).. (0,0);
\draw[bbline] (0.5,0.4) .. controls (0.4,0) .. (0,0);
\draw[dotted, thick] (-0.005,0.35) -- (0.235,0.35);
\draw[dotted, thick] (-0.005,-0.35) -- (0.235,-0.35);
\draw (0,0) -- (0,0) node[fill=red!20, draw,ellipse,scale=.8] {$2m$};
\end{tikzpicture}\\ 
 & \\
\begin{tikzpicture}[baseline=-0.5ex]
\draw[sline] (-0.5,-0.4) .. controls (-0.4,0) .. (0,0);
\draw[tline] (-0.25,-0.4) .. controls (-0.2,0).. (0,0);
\draw[bbline] (0.5,-0.4) .. controls (0.4,0) .. (0,0);
\draw[tline] (-0.5,0.4) .. controls (-0.4,0) .. (0,0);
\draw[sline] (-0.25,0.4) .. controls (-0.2,0).. (0,0);
\draw[bbline] (0.5,0.4) .. controls (0.4,0) .. (0,0);
\draw[dotted, thick] (-0.005,0.35) -- (0.235,0.35);
\draw[dotted, thick] (-0.005,-0.35) -- (0.235,-0.35);
\draw (0,0) -- (0,0) node[fill=red!20, draw,ellipse,scale=.8] {$2m$};
\end{tikzpicture}
 &\begin{tikzpicture}[baseline=-0.5ex]
\draw[tline] (-0.5,-0.4) -- (-0.5,0.4);
\draw[sline] (-0.25,-0.4) -- (-0.25,0.4);
\draw[bbline] (0.25,-0.4) -- (0.25,0.4);
\draw[bbline] (0.5,-0.4) -- (0.5,0.4);
\draw[dotted, thick] (-0.12,0.35) -- (0.12,0.35);
\draw[dotted, thick] (-0.12,-0.35) -- (0.12,-0.35);
\draw[bbline] (0.5,0) -- (0,0);
\draw[tline] (-0.5,0) -- (0,0);
\draw (0,0) -- (0,0) node[fill=red!20, draw] {$\scriptstyle{JW}$};
\end{tikzpicture}
\end{pmatrix}
\begin{pmatrix}
\begin{tikzpicture}[baseline=-0.5ex]
\draw[sline] (0,-0.4) -- (0,0.4);
\draw[tline] (0.2,-0.4) -- (0.2,0.4);
\draw[sline] (0.4,-0.4) -- (0.4,0.4);
\draw[dotted, thick] (0.58,0) -- (0.82,0);
\draw[bbline] (1,0)--(1,0.4);
\draw[bbdot] (1,0) circle (2pt);
\end{tikzpicture} \\ \\
\begin{tikzpicture}[baseline=-0.5ex]
\draw[tline] (0,0)--(0,0.4);
\draw[tdot] (0,0) circle (2pt);
\draw[sline] (0.2,-0.4) -- (0.2,0.4);
\draw[tline] (0.4,-0.4) -- (0.4,0.4);
\draw[dotted, thick] (0.58,0) -- (0.82,0);
\draw[bbline] (1,-0.4) -- (1,0.4);
\end{tikzpicture}
\end{pmatrix}\\[0.5cm]
\enspace & = \enspace \frac{1}{4} \left(
\begin{tikzpicture}[baseline=-0.5ex]
\draw[sline] (-0.5,-0.4) -- (-0.5,0.4);
\draw[tline] (-0.25,-0.4) -- (-0.25,0.4);
\draw[bbline] (0.25,-0.4) -- (0.25,0.4);
\draw[bbline] (0.5,-0.2) -- (0.5,0.2);
\draw[bbdot] (0.5,0.2) circle (2pt);
\draw[bbdot] (0.5,-0.2) circle (2pt);
\draw[dotted, thick] (-0.12,0.35) -- (0.12,0.35);
\draw[dotted, thick] (-0.12,-0.35) -- (0.12,-0.35);
\draw[bbline] (0.5,0) -- (0,0);
\draw[sline] (-0.5,0) -- (0,0);
\draw (0,0) -- (0,0) node[fill=red!20, draw] {$\scriptstyle{JW}$};
\end{tikzpicture}
\enspace + \enspace
\begin{tikzpicture}[baseline=-0.5ex]
\draw[tline] (-0.5,-0.3) .. controls (-0.4,0) .. (0,0);
\draw[tdot] (-0.5,-0.3) circle (2pt);
\draw[sline] (-0.25,-0.4) .. controls (-0.2,0).. (0,0);
\draw[bbline] (0.5,-0.4) .. controls (0.4,0) .. (0,0);
\draw[sline] (-0.5,0.4) .. controls (-0.4,0) .. (0,0);
\draw[tline] (-0.25,0.4) .. controls (-0.2,0).. (0,0);
\draw[bbline] (0.5,0.3) .. controls (0.4,0) .. (0,0);
\draw[bbdot] (0.5,0.3) circle (2pt);
\draw[dotted, thick] (-0.005,0.35) -- (0.235,0.35);
\draw[dotted, thick] (-0.005,-0.35) -- (0.235,-0.35);
\draw (0,0) -- (0,0) node[fill=red!20, draw,ellipse,scale=.8] {$2m$};
\end{tikzpicture}
\enspace + \enspace
\begin{tikzpicture}[baseline=-0.5ex]
\draw[sline] (-0.5,-0.4) .. controls (-0.4,0) .. (0,0);
\draw[tline] (-0.25,-0.4) .. controls (-0.2,0).. (0,0);
\draw[bbline] (0.5,-0.3) .. controls (0.4,0) .. (0,0);
\draw[bbdot] (0.5,-0.3) circle (2pt);
\draw[tline] (-0.5,0.3) .. controls (-0.4,0) .. (0,0);
\draw[tdot] (-0.5,0.3) circle (2pt);
\draw[sline] (-0.25,0.4) .. controls (-0.2,0).. (0,0);
\draw[bbline] (0.5,0.4) .. controls (0.4,0) .. (0,0);
\draw[dotted, thick] (-0.005,0.35) -- (0.235,0.35);
\draw[dotted, thick] (-0.005,-0.35) -- (0.235,-0.35);
\draw (0,0) -- (0,0) node[fill=red!20, draw,ellipse,scale=.8] {$2m$};
\end{tikzpicture}
\enspace + \enspace
\begin{tikzpicture}[baseline=-0.5ex]
\draw[tline] (-0.5,-0.2) -- (-0.5,0.2);
\draw[tdot] (-0.5,0.2) circle (2pt);
\draw[tdot] (-0.5,-0.2) circle (2pt);
\draw[sline] (-0.25,-0.4) -- (-0.25,0.4);
\draw[bbline] (0.25,-0.4) -- (0.25,0.4);
\draw[bbline] (0.5,-0.4) -- (0.5,0.4);
\draw[dotted, thick] (-0.12,0.35) -- (0.12,0.35);
\draw[dotted, thick] (-0.12,-0.35) -- (0.12,-0.35);
\draw[bbline] (0.5,0) -- (0,0);
\draw[tline] (-0.5,0) -- (0,0);
\draw (0,0) -- (0,0) node[fill=red!20, draw] {$\scriptstyle{JW}$};
\end{tikzpicture}
\right)\\[0.5cm]
\enspace & \stackrel{(\ref{pro:4})}{=} \enspace \frac{1}{4} \left(
\begin{tikzpicture}[baseline=-0.5ex]
\draw[sline] (-0.5,-0.4) -- (-0.5,0.4);
\draw[tline] (-0.25,-0.4) -- (-0.25,0.4);
\draw[bbline] (0.25,-0.4) -- (0.25,0.4);
\draw[bbdot] (0.5,0) circle (2pt);
\draw[dotted, thick] (-0.12,0.35) -- (0.12,0.35);
\draw[dotted, thick] (-0.12,-0.35) -- (0.12,-0.35);
\draw[bbline] (0.5,0) -- (0,0);
\draw[sline] (-0.5,0) -- (0,0);
\draw (0,0) -- (0,0) node[fill=red!20, draw] {$\scriptstyle{JW}$};
\end{tikzpicture}
\enspace + \enspace 2 \enspace
\begin{tikzpicture}[baseline=-0.5ex]
\draw[tline] (-0.6,0) -- (0,0);
\draw[tdot] (-0.6,0) circle (2pt);
\draw[sline] (-0.4,-0.4) .. controls (-0.3,0).. (0,0);
\draw[bbline] (0.4,-0.4) .. controls (0.3,0) .. (0,0);
\draw[bbline] (0.4,0.4) .. controls (0.3,0) .. (0,0);
\draw[sline] (-0.4,0.4) .. controls (-0.3,0) .. (0,0);
\draw[tline] (-0.15,0.4) .. controls (-0.1,0).. (0,0);
\draw[tline] (-0.15,-0.4) .. controls (-0.1,0).. (0,0);
\draw[bbline] (0.6,0) -- (0,0);
\draw[bbdot] (0.6,0) circle (2pt);
\draw[dotted, thick] (-0.005,0.35) -- (0.235,0.35);
\draw[dotted, thick] (-0.005,-0.35) -- (0.235,-0.35);
\draw (0,0) -- (0,0) node[fill=red!20, draw,ellipse,scale=.8] {$2m$};
\end{tikzpicture}
\enspace + \enspace
\begin{tikzpicture}[baseline=-0.5ex]
\draw[tdot] (-0.5,0) circle (2pt);
\draw[sline] (-0.25,-0.4) -- (-0.25,0.4);
\draw[bbline] (0.25,-0.4) -- (0.25,0.4);
\draw[bbline] (0.5,-0.4) -- (0.5,0.4);
\draw[dotted, thick] (-0.12,0.35) -- (0.12,0.35);
\draw[dotted, thick] (-0.12,-0.35) -- (0.12,-0.35);
\draw[bbline] (0.5,0) -- (0,0);
\draw[tline] (-0.5,0) -- (0,0);
\draw (0,0) -- (0,0) node[fill=red!20, draw] {$\scriptstyle{JW}$};
\end{tikzpicture}
\right)\\[0.5cm]
\enspace & \stackrel{(\ref{pro:5})}{=} \enspace \frac{1}{4} \left(
\begin{tikzpicture}[baseline=-0.5ex]
\draw[sline] (-0.5,-0.4) -- (-0.5,0.4);
\draw[tline] (-0.25,-0.4) -- (-0.25,0.4);
\draw[bbline] (0.25,-0.4) -- (0.25,0.4);
\draw[bbdot] (0.5,0) circle (2pt);
\draw[dotted, thick] (-0.12,0.35) -- (0.12,0.35);
\draw[dotted, thick] (-0.12,-0.35) -- (0.12,-0.35);
\draw[bbline] (0.5,0) -- (0,0);
\draw[sline] (-0.5,0) -- (0,0);
\draw (0,0) -- (0,0) node[fill=red!20, draw] {$\scriptstyle{JW}$};
\end{tikzpicture}
\enspace + \enspace 3 \enspace
\begin{tikzpicture}[baseline=-0.5ex]
\draw[tdot] (-0.5,0) circle (2pt);
\draw[sline] (-0.25,-0.4) -- (-0.25,0.4);
\draw[bbline] (0.25,-0.4) -- (0.25,0.4);
\draw[bbline] (0.5,-0.4) -- (0.5,0.4);
\draw[dotted, thick] (-0.12,0.35) -- (0.12,0.35);
\draw[dotted, thick] (-0.12,-0.35) -- (0.12,-0.35);
\draw[bbline] (0.5,0) -- (0,0);
\draw[tline] (-0.5,0) -- (0,0);
\draw (0,0) -- (0,0) node[fill=red!20, draw] {$\scriptstyle{JW}$};
\end{tikzpicture}
\right)\\[0.5cm]
\enspace & \stackrel{\ref{lem:jw}}{=} \enspace [2]_\zeta \enspace
\begin{tikzpicture}[baseline=-0.5ex]
\draw[sline] (0,-0.4) -- (0,0.4);
\draw[tline] (0.2,-0.4) -- (0.2,0.4);
\draw[bbline] (0.8,-0.4) -- (0.8,-0.2);
\draw[bbdot] (0.8,-0.2) circle (2pt);
\draw[bbline] (0.8,0.4) -- (0.8,0.2);
\draw[bbdot] (0.8,0.2) circle (2pt);
\draw[dotted, thick] (0.38,0) -- (0.62,0);
\end{tikzpicture}.
\end{align*}

\subsubsection{Relation \ref{eqn:4}}
Clearly, we have:
\[
\begin{tikzpicture}[baseline=-0.5ex]
\draw[tline] (0,-0.15) -- (0,0.15);
\draw[tdot] (0,-0.15) circle (2pt);
\draw[tdot] (0,0.15) circle (2pt);
\draw[sline] (0.2,-0.4) -- (0.2,0.4);
\end{tikzpicture}
\enspace \stackrel{\ref{pro:2}}{=} \enspace -[2]_\zeta \enspace  
\begin{tikzpicture}[baseline=-0.5ex]
\draw[sline] (0,-0.4) -- (0,-0.2);
\draw[sdot] (0,-0.2) circle (2pt);
\draw[sdot] (0,0.2) circle (2pt);
\draw[sline] (0,0.4) -- (0,0.2);
\end{tikzpicture}.
\]

\subsubsection{Relation \ref{eqn:6} and \ref{eqn:7}}
The bimodules $B_{\scol}B_{\tcol}$ and $B_{\tcol}B_{\scol}$ are indecomposable and therefore the corresponding idempotent is trivial, hence by isotopy invariance

\[ 
\begin{tikzpicture}[baseline=-0.5ex]
\draw[sline] (0,-0.4) -- (0,0.2);
\draw[sdot] (0,0.2) circle (2pt);
\draw[tline] (0.2,-0.2) -- (0.2,0.4);
\draw[tdot] (0.2,-0.2) circle (2pt);
\end{tikzpicture}
\enspace \stackrel{\ref{pro:4}}{=} \enspace 
\begin{tikzpicture}[baseline=-0.5ex]
\draw[sline] (0,-0.4) -- (0,-0.2);
\draw[sdot] (0,-0.2) circle (2pt);
\draw[tline] (0,0.2) -- (0,0.4);
\draw[tdot] (0,0.2) circle (2pt);
\end{tikzpicture}
\enspace \stackrel{\ref{pro:4}}{=} \enspace
\begin{tikzpicture}[baseline=-0.5ex]
\draw[tline] (0,0.4) -- (0,-0.2);
\draw[tdot] (0,-0.2) circle (2pt);
\draw[sline] (0.2,-0.4) -- (0.2,0.2);
\draw[sdot] (0.2,0.2) circle (2pt);
\end{tikzpicture}.\]

\subsubsection{Relation \ref{eqn:8} and \ref{eqn:9}}
The relations \ref{eqn:8} and \ref{eqn:9} for $i\in \{2,\ldots,m-2\}$ are direct consequences of the following lemma:

\begin{lem}\label{lem:rel8}
For $i\in\{ 2,\ldots, m-2\}$ the following holds.
\begin{align*}
\begin{tikzpicture}[baseline=-0.5ex]
\draw[sline,sharp corners] (-0.5,-0.4) -- (-0.5,0) -- (0,0);
\draw[tline] (-0.25,-0.4) -- (-0.25,0.4);
\draw[bbline] (0.25,-0.4) -- (0.25,0.4);
\draw[bbline, sharp corners] (0,0) -- (0.5,0) -- (0.5,0.4);
\draw[dotted, thick] (-0.12,0.35) -- (0.12,0.35);
\draw[dotted, thick] (-0.12,-0.35) -- (0.12,-0.35);
\draw (0,0) -- (0,0) node[fill=red!20, draw] {$\scriptstyle{JW_i}$};
\end{tikzpicture}
\enspace & = \enspace \textstyle{\frac{1}{[i]_\zeta}} \enspace
\begin{tikzpicture}[baseline=-0.5ex]
\draw[tline] (0,0.4) -- (0,0.2);
\draw[tdot] (0,0.2) circle (2pt);
\draw[sline] (0.2,-0.4) -- (0.2,0.4);
\draw[tline] (0.4,-0.4) -- (0.4,0.4);
\draw[dotted, thick] (0.58,0) -- (0.82,0);
\draw[bbline] (1,-0.4) -- (1,-0.2);
\draw[bbdot] (1,-0.2) circle (2pt);
\end{tikzpicture}
\enspace + \enspace 
\begin{tikzpicture}[baseline=-0.5ex]
\draw[sline] (0,-0.4) -- (0,-0.2);
\draw[sdot] (0,-0.2) circle (2pt);
\draw[tline] (0.2,-0.4) -- (0.2,0.4);
\draw[sline] (0.4,-0.4) -- (0.4,0.4);
\draw[dotted, thick] (0.58,0) -- (0.82,0);
\draw[bbline] (1,0.4) -- (1,0.2);
\draw[bbdot] (1,0.2) circle (2pt);
\end{tikzpicture},\\[0.5cm]
\begin{tikzpicture}[baseline=-0.5ex]
\draw[tline,sharp corners] (-0.5,0.4) -- (-0.5,0) -- (0,0);
\draw[sline] (-0.25,-0.4) -- (-0.25,0.4);
\draw[bbline] (0.25,-0.4) -- (0.25,0.4);
\draw[bbline, sharp corners] (0,0) -- (0.5,0) -- (0.5,-0.4);
\draw[dotted, thick] (-0.12,0.35) -- (0.12,0.35);
\draw[dotted, thick] (-0.12,-0.35) -- (0.12,-0.35);
\draw (0,0) -- (0,0) node[fill=red!20, draw] {$\scriptstyle{JW_i}$};
\end{tikzpicture}
\enspace & = \enspace 
\begin{tikzpicture}[baseline=-0.5ex]
\draw[tline] (0,0.4) -- (0,0.2);
\draw[tdot] (0,0.2) circle (2pt);
\draw[sline] (0.2,-0.4) -- (0.2,0.4);
\draw[tline] (0.4,-0.4) -- (0.4,0.4);
\draw[dotted, thick] (0.58,0) -- (0.82,0);
\draw[bbline] (1,-0.4) -- (1,-0.2);
\draw[bbdot] (1,-0.2) circle (2pt);
\end{tikzpicture}
\enspace + \enspace \textstyle{\frac{1}{[i]_\zeta}} \enspace
\begin{tikzpicture}[baseline=-0.5ex]
\draw[sline] (0,-0.4) -- (0,-0.2);
\draw[sdot] (0,-0.2) circle (2pt);
\draw[tline] (0.2,-0.4) -- (0.2,0.4);
\draw[sline] (0.4,-0.4) -- (0.4,0.4);
\draw[dotted, thick] (0.58,0) -- (0.82,0);
\draw[bbline] (1,0.4) -- (1,0.2);
\draw[bbdot] (1,0.2) circle (2pt);
\end{tikzpicture},
\intertext{In particular,}
\begin{tikzpicture}[baseline=-0.5ex]
\draw[sline,sharp corners] (-0.5,-0.4) -- (-0.5,0) -- (0,0);
\draw[tline] (-0.25,-0.4) -- (-0.25,0.4);
\draw[bbline] (0.25,-0.4) -- (0.25,0.4);
\draw[bbline,sharp corners] (0,0) -- (0.5,0) -- (0.5,0.4);
\draw[dotted, thick] (-0.12,0.35) -- (0.12,0.35);
\draw[dotted, thick] (-0.12,-0.35) -- (0.12,-0.35);
\draw (0,0) -- (0,0) node[fill=red!20, draw] {$\scriptstyle{JW}$};
\end{tikzpicture}
\enspace & = \enspace\enspace
\begin{tikzpicture}[baseline=-0.5ex]
\draw[tline] (0,0.4) -- (0,0.2);
\draw[tdot] (0,0.2) circle (2pt);
\draw[sline] (0.2,-0.4) -- (0.2,0.4);
\draw[tline] (0.4,-0.4) -- (0.4,0.4);
\draw[dotted, thick] (0.58,0) -- (0.82,0);
\draw[bbline] (1,-0.4) -- (1,-0.2);
\draw[bbdot] (1,-0.2) circle (2pt);
\end{tikzpicture}
\enspace + \enspace 
\begin{tikzpicture}[baseline=-0.5ex]
\draw[sline] (0,-0.4) -- (0,-0.2);
\draw[sdot] (0,-0.2) circle (2pt);
\draw[tline] (0.2,-0.4) -- (0.2,0.4);
\draw[sline] (0.4,-0.4) -- (0.4,0.4);
\draw[dotted, thick] (0.58,0) -- (0.82,0);
\draw[bbline] (1,0.4) -- (1,0.2);
\draw[bbdot] (1,0.2) circle (2pt);
\end{tikzpicture}
\enspace = \enspace
\begin{tikzpicture}[baseline=-0.5ex]
\draw[tline,sharp corners] (-0.5,0.4) -- (-0.5,0) -- (0,0);
\draw[sline] (-0.25,-0.4) -- (-0.25,0.4);
\draw[bbline] (0.25,-0.4) -- (0.25,0.4);
\draw[bbline, sharp corners] (0,0) -- (0.5,0) -- (0.5,-0.4);
\draw[dotted, thick] (-0.12,0.35) -- (0.12,0.35);
\draw[dotted, thick] (-0.12,-0.35) -- (0.12,-0.35);
\draw (0,0) -- (0,0) node[fill=red!20, draw] {$\scriptstyle{JW}$};
\end{tikzpicture}.
\end{align*}
\end{lem}

\begin{prf} Similarly to the proof of \cref{lem:jw}.
\end{prf}

For $i=m-1$ relation \ref{eqn:8} and \ref{eqn:9} agree (see \cref{rem:rel}) and the relation is proven simirlarly to relations \ref{eqn:3} and \ref{eqn:5} since $B_{w_0}$ is embedded diagonally.

\subsubsection{Relation \ref{eqn:10} and \ref{eqn:12}}
Inspecting both sides of the relation (\ref{eqn:10}) gives us
\begin{alignat*}{2}
\begin{tikzpicture}[baseline=-0.5ex]
\draw[sline] (-0.55,-0.4) -- (-0.55,0.4);
\draw[tline] (-0.25,-0.4) -- (-0.25,0.4);
\draw[bbline] (0.25,-0.4) -- (0.25,0.4);
\draw[bbdot] (0.55,-0.25) circle (2pt);
\draw[dotted, thick] (-0.12,0.35) -- (0.12,0.35);
\draw[dotted, thick] (-0.12,-0.35) -- (0.12,-0.35);
\draw[bbline] (0.75,0.4) -- (0.75,0.2);
\draw[bbdot] (0.75,0.2) circle (2pt);
\draw[bbline] (0.55,0) -- (0,0);
\draw[bbline] (0.55,-0.25) -- (0.55,0.4);
\draw[sline] (-0.55,0) -- (0,0);
\draw (0,0) -- (0,0) node[fill=red!20, draw] {$\scriptstyle{JW_j}$};
\end{tikzpicture}
\enspace & = \enspace
\begin{tikzpicture}[baseline=-0.5ex]
\draw[sline] (0,-0.4) -- (0,0.4);
\draw[bbline] (0.4,-0.4) -- (0.4,0.4);
\draw[bbline] (0.6,0.4) -- (0.6,0.2);
\draw[bbline] (0.8,0.4) -- (0.8,0.2);
\draw[dotted, thick] (0.08,0) -- (0.32,0);
\draw[bbdot] (0.6,0.2) circle (2pt);
\draw[bbdot] (0.8,0.2) circle (2pt);
\end{tikzpicture},\quad\quad\quad&
\begin{tikzpicture}[baseline=-0.5ex]
\draw[bbline] (0.55,-0.4) -- (0.55,0.4);
\draw[bbline] (0.25,-0.4) -- (0.25,0.4);
\draw[sline] (-0.25,-0.4) -- (-0.25,0.4);
\draw[tdot] (-0.55,-0.25) circle (2pt);
\draw[dotted, thick] (0.12,0.35) -- (-0.12,0.35);
\draw[dotted, thick] (0.12,-0.35) -- (-0.12,-0.35);
\draw[sline] (-0.75,0.4) -- (-0.75,0.2);
\draw[sdot] (-0.75,0.2) circle (2pt);
\draw[tline] (-0.55,0) -- (0,0);
\draw[tline] (-0.55,-0.25) -- (-0.55,0.4);
\draw[bbline] (--0.55,0) -- (0,0);
\draw (0,0) -- (0,0) node[fill=red!20, draw] {$\scriptstyle{JW_j}$};
\end{tikzpicture}
\enspace & = \enspace
\begin{tikzpicture}[baseline=-0.5ex]
\draw[bbline] (0,-0.4) -- (0,0.4);
\draw[sline] (-0.4,-0.4) -- (-0.4,0.4);
\draw[tline] (-0.6,0.4) -- (-0.6,0.2);
\draw[sline] (-0.8,0.4) -- (-0.8,0.2);
\draw[dotted, thick] (-0.08,0) -- (-0.32,0);
\draw[tdot] (-0.6,0.2) circle (2pt);
\draw[sdot] (-0.8,0.2) circle (2pt);
\end{tikzpicture}
\end{alignat*}

since by property \ref{pro:1} every summand except the identity is cancelled (when pre-/post-composed with the idempotents). The right hand sides are equal by the following lemma (which proves the desired relation).

\begin{lem}
In $\sbim$ we have the following equality as morphisms $B_{{}_\scol i} \to B_{{}_\scol (i+2)}$ for $i\geq 0$
\[
\begin{tikzpicture}[baseline=-0.5ex]
\draw[sline] (0,-0.4) -- (0,0.4);
\draw[bbline] (0.4,-0.4) -- (0.4,0.4);
\draw[bbline] (0.6,0.4) -- (0.6,0.2);
\draw[bbline] (0.8,0.4) -- (0.8,0.2);
\draw[dotted, thick] (0.08,0) -- (0.32,0);
\draw[bbdot] (0.6,0.2) circle (2pt);
\draw[bbdot] (0.8,0.2) circle (2pt);
\end{tikzpicture}
\enspace = \enspace
\begin{tikzpicture}[baseline=-0.5ex]
\draw[bbline] (0,-0.4) -- (0,0.4);
\draw[sline] (-0.4,-0.4) -- (-0.4,0.4);
\draw[tline] (-0.6,0.4) -- (-0.6,0.2);
\draw[sline] (-0.8,0.4) -- (-0.8,0.2);
\draw[dotted, thick] (-0.08,0) -- (-0.32,0);
\draw[tdot] (-0.6,0.2) circle (2pt);
\draw[sdot] (-0.8,0.2) circle (2pt);
\end{tikzpicture}.
\]
\end{lem}

\begin{prf}
The case $i=0$ is trivial, so let us assume $i=1$.
\begin{align*}
\begin{tikzpicture}[baseline=-0.5ex]
\draw[sline] (0,-0.4) -- (0,0.4);
\draw[tline] (0.2,0.4) -- (0.2,0.2);
\draw[sline] (0.4,0.4) -- (0.4,0.2);
\draw[tdot] (0.2,0.2) circle (2pt);
\draw[sdot] (0.4,0.2) circle (2pt);
\end{tikzpicture}
\enspace  \stackrel{\ref{pro:2}}{=} \enspace \tiny{\frac{1}{2}} \enspace \left( 
\begin{tikzpicture}[baseline=-0.5ex]
\draw[sline] (-0.3,-0.4) -- (-0.3,0.4);
\draw[sline] (-0.3,0) .. controls  (0.3,0) and (0.4,0.1) .. (0.4,0.4);
\draw[tline] (0.2,0.4) -- (0.2,0.2);
\draw[tdot] (0.2,0.2) circle (2pt);
\draw (-0.06,0.2) -- (-0.06,0.2) node {$\alpha_\scol$};
\end{tikzpicture}
\enspace + \enspace
\begin{tikzpicture}[baseline=-0.5ex]
\draw[sline] (-0.3,-0.4) -- (-0.3,0.4);
\draw[sline] (-0.3,0) .. controls  (0.3,0) and (0.4,0.1) .. (0.4,0.4);
\draw[tline] (0.2,0.4) -- (0.2,0.2);
\draw[tdot] (0.2,0.2) circle (2pt);
\draw (-0.06,-0.2) -- (-0.06,-0.2) node {$\alpha_\scol$};
\end{tikzpicture} \right)
\enspace  \stackrel{\ref{pro:1}}{=} \enspace \tiny{\frac{1}{2}} \enspace
\begin{tikzpicture}[baseline=-0.5ex]
\draw[sline] (-0.3,-0.4) -- (-0.3,0.4);
\draw[sline] (-0.3,0) .. controls  (0.3,0) and (0.4,0.1) .. (0.4,0.4);
\draw[tline] (0.2,0.4) -- (0.2,0.2);
\draw[tdot] (0.2,0.2) circle (2pt);
\draw (-0.06,0.2) -- (-0.06,0.2) node {$\alpha_\scol$};
\end{tikzpicture}
\enspace  \stackrel{\ref{pro:2}}{=} \enspace \tiny{\frac{\scol\alpha_\scol}{2}}
\begin{tikzpicture}[baseline=-0.5ex]
\draw[sline] (-0.3,-0.4) -- (-0.3,0.4);
\draw[sline] (-0.3,0) .. controls  (0.3,0) and (0.4,0.1) .. (0.4,0.4);
\draw[tline] (0.2,0.4) -- (0.2,0.2);
\draw[tdot] (0.2,0.2) circle (2pt);
\end{tikzpicture} 
\enspace + \enspace \underbrace{\tiny{\frac{\partial_\scol(\alpha_\scol)}{2}}}_{=1}
\begin{tikzpicture}[baseline=-0.5ex]
\draw[sline] (0,-0.4) .. controls  (0,0) and (0.4,0) .. (0.4,0.4);
\draw[sline] (0,0.4) -- (0,0.2);
\draw[sdot] (0,0.2) circle (2pt);
\draw[tline] (0.2,0.4) -- (0.2,0.2);
\draw[tdot] (0.2,0.2) circle (2pt);
\end{tikzpicture}
\enspace & \stackrel{\ref{pro:1}}{=} \enspace 
\begin{tikzpicture}[baseline=-0.5ex]
\draw[sline] (0.4,-0.4) -- (0.4,0.4);
\draw[sline] (0,0.4) -- (0,0.2);
\draw[sdot] (0,0.2) circle (2pt);
\draw[tline] (0.2,0.4) -- (0.2,0.2);
\draw[tdot] (0.2,0.2) circle (2pt);
\end{tikzpicture}.
\end{align*}
The general case follows immediately by using the above equation repeatedly. 
\end{prf}
Relation \ref{eqn:12} can be dealt with in the same way by flipping the diagrams at the horizontal axis.

\subsubsection{Relation \ref{eqn:11} and \ref{eqn:13}}
Translating relation \ref{eqn:11} into the Soergel bimodule setting yields
\[
\begin{tikzpicture}[baseline=-0.5ex]
\draw[sline] (-0.55,-0.4) -- (-0.55,0.4);
\draw[tline] (-0.25,-0.4) -- (-0.25,0.4);
\draw[bbline] (0.25,-0.4) -- (0.25,0.4);
\draw[bbdot] (0.55,-0.25) circle (2pt);
\draw[dotted, thick] (-0.12,0.35) -- (0.12,0.35);
\draw[dotted, thick] (-0.12,-0.35) -- (0.12,-0.35);
\draw[tline] (-0.75,0.4) -- (-0.75,0.2);
\draw[tdot] (-0.75,0.2) circle (2pt);
\draw[bbline] (0.55,0) -- (0,0);
\draw[bbline] (0.55,-0.25) -- (0.55,0.4);
\draw[sline] (-0.55,0) -- (0,0);
\draw (0,0) -- (0,0) node[fill=red!20, draw] {$\scriptstyle{JW_j}$};
\end{tikzpicture}
\enspace = \enspace
\begin{tikzpicture}[baseline=-0.5ex]
\draw[bbline] (0.55,-0.4) -- (0.55,0.4);
\draw[bbline] (0.25,-0.4) -- (0.25,0.4);
\draw[sline] (-0.25,-0.4) -- (-0.25,0.4);
\draw[tdot] (-0.55,-0.25) circle (2pt);
\draw[dotted, thick] (0.12,0.35) -- (-0.12,0.35);
\draw[dotted, thick] (0.12,-0.35) -- (-0.12,-0.35);
\draw[bbline] (0.75,0.4) -- (0.75,0.2);
\draw[bbdot] (0.75,0.2) circle (2pt);
\draw[tline] (-0.55,0) -- (0,0);
\draw[tline] (-0.55,-0.25) -- (-0.55,0.4);
\draw[bbline] (--0.55,0) -- (0,0);
\draw (0,0) -- (0,0) node[fill=red!20, draw] {$\scriptstyle{JW_j}$};
\end{tikzpicture}.
\]

In $JW_j$ on both sides above, every summand is cancelled by property (\ref{pro:1}) except the one corresponding to the identity with scalar $1$ (after pre-/post-composing with the idempotens). Hence, the equation above becomes
\[
\begin{tikzpicture}[baseline=-0.5ex]
\draw[tline] (-0.5,0.2) -- (-0.5,0.4);
\draw[tdot] (-0.5,0.2) circle (2pt);
\draw[bbline] (0.5,0.2) -- (0.5,0.4);
\draw[bbdot] (0.5,0.2) circle (2pt);
\draw[sline] (-0.25,-0.4) -- (-0.25,0.4);
\draw[bbline] (0.25,-0.4) -- (0.25,0.4);
\draw[dotted, thick] (0.12,0) -- (-0.12,0);
\end{tikzpicture}
\enspace = \enspace
\begin{tikzpicture}[baseline=-0.5ex]
\draw[tline] (-0.5,0.2) -- (-0.5,0.4);
\draw[tdot] (-0.5,0.2) circle (2pt);
\draw[bbline] (0.5,0.2) -- (0.5,0.4);
\draw[bbdot] (0.5,0.2) circle (2pt);
\draw[sline] (-0.25,-0.4) -- (-0.25,0.4);
\draw[bbline] (0.25,-0.4) -- (0.25,0.4);
\draw[dotted, thick] (0.12,0) -- (-0.12,0);
\end{tikzpicture}
\]
which clearly holds. Relation \ref{eqn:13} holds for the same reasons after flipping the above diagrams at the horizontal axis.

\subsection*{Step \MakeUppercase{\romannumeral 3}}
The path algebra $\Pm$ has by \cref{lem:dim} dimension $\sum_{x,y} |V_{\leq x}\cap V_{\leq y}|$ over $\DR$, whereas by \cref{thm:dim}
\[\op{dim}_\DR(\mathcal A) \enspace = \enspace \sum_{x,y\in W} \op{dim}_\DR\; \op{Hom}^\bullet_{\smod}(B_x,B_y) \enspace = \enspace \sum_{x,y\in W} |W_{\leq x}\cap W_{\leq y}|. \]
Altogether we conclude that $\varphi$ is an isomorphism.

\hfill\qed

\section{Koszul Self-Duality of $\mathcal A$}
\subsection{Linear resolutions of the Standard Modules}
The indexing set $\Lambda$ of isomorphism classes of simple graded $\textbf P_m$-modules is in bijection with $W$ by \cref{prop:simple}. Hence, the reverse Bruhat order $\leq_r$ turns $\Lambda$ into a finite partially ordered set $(\Lambda,\leq_r)$. For $x\in W$ define the left modules
\begin{align*}
\Delta(x):= P(x)/M(x)
\end{align*}
where $M(x):= \left\langle \text{all paths starting in }x \text{ passing through } y, y \nleq_r x \text{ of length }\geq 1 \right\rangle $. Therefore we have the short exact sequence:

\begin{align}
\xymatrix{0\ar[r] & M(x) \ar[r] & P(x)\ar[r]^\pi  & \Delta(x) \ar[r] & 0}. \label{eqn:ses}
\end{align}

Note that $e_y\Delta(x)=0$ for all $y\nleq_r x$ which means that no path in $\Delta(x)$ ends in such a $y$. A precise formula for the dimensions of the involved modules in this short exact sequence gives the following lemma.

\begin{lem} For $0\leq i \leq m$ we have: \label{lem:dimmod}
\begin{align}
\dim_\DR P({}_\scol i) \enspace &  = \enspace 
\begin{cases}
2m & \text{if } i=0, \\
2i(2m-i) & \text{if } 1\leq i \leq m-1, \\
1+ 2m^2 & \text{if } i=m ;
\end{cases} \label{eqn:dimp} \\
\dim_\DR \Delta({}_\scol i)\enspace &  = \enspace 
\begin{cases}
2(m-i) & \text{if } 0\leq i \leq m-1, \\
1 & \text{if } i=m;
\end{cases} \label{eqn:dimstd} \\
\dim_\DR M({}_\scol i) \enspace &  = \enspace 
\begin{cases}
0 & \text{if } i=0, \\
2i(2m-i) - 2(m-i) & \text{if } 1\leq i \leq m-1, \\
2m^2 & \text{if  } i=m.
\end{cases} \label{eqn:dimsub}
\end{align}
\end{lem}

\begin{prf}
For fixed $x\in W$ the indecomposable projective cover $P(x)$ has a basis consisting of all paths starting in $x$. Hence by \cref{lem:dim} and its proof we have
\begin{align}
\dim_\DR P(x) \enspace & =  \enspace \sum_{y\in W} \dim_\DR e_y\Pm e_x  \enspace =  \enspace \sum_{y\in W} |V_{\geq_r y}\cap V_{\geq_r x}|. \label{eqn:pdim1}
\intertext{Note that we reversed the order on the vertices. For $x={}_\scol i$ and $1\leq i \leq m-1$ this becomes}
\enspace &= \enspace \underbrace{1}_{y=e} + \underbrace{2\sum_{j=1}^{i-1}2j}_{1\leq \ell(y)\leq i-1}  + \underbrace{2i + (2i-1)}_{\ell(y)=i} + \underbrace{2\sum_{j=i+1}^{m-1}2i }_{i+1\leq \ell(y)\leq m-1} + \underbrace{2i}_{y=w_0} \label{eqn:pdim2}\\
\enspace &= \enspace 2i(2m-i). \notag
\end{align}
Consider \ref{eqn:pdim2} which reduces for $i=m$ to $1 + \sum_{j=1}^{m-1}4j + 2m = 1 + 2m^2$. The case $i=0$ follows easily from \ref{eqn:pdim1}.
By construction of $\Delta(x)$ we have
\begin{align*}
\dim_\DR \Delta(x) \enspace =  \enspace \sum_{y\in W} |V_{\leq_r y}\cap \{x \}| 
\enspace  = \enspace \sum_{y \leq_r x} 1
\enspace  = \enspace
\begin{cases}
2(m-\ell(x)) & \text{if  } x\neq w_0, \\
1 & \text{if  } x=w_0
\end{cases}
\end{align*}
which yields \ref{eqn:dimstd}. The dimension of $M({}_\scol i)$ follows directly from the above combined with the short exact sequence in \ref{eqn:ses}.
\end{prf}

\begin{thm}\label{thm:qh}
The set $\{\Delta(x)\}_{x\in W}$ defines a quasi-hereditary structure on $(\textbf P_m,(W,\leq_r))$, i.e. for $x\in W$ the (left) module $\Delta(x)$ is the \emph{(left) standard module}.
\end{thm}
\begin{prf}
We prove that $P(x)$ has a $\Delta$-filtration with subquotients isomorphic to $\Delta(y)$ for $y\geq_r x$ (each with multiplicity 1) and $x={}_\scol i$ via induction over $i$ (this proves the first condition for being quasi-hereditary).

\begin{figure}
\[\xymatrix{ & {}_\scol i \ar[dl]_\omega \ar[d]^{\omega^\prime}\\
{}_\tcol (i-1) \ar[d]_\beta \ar[dr]^\alpha & {}_\scol (i-1) \\
{}_\tcol (i-2) & {}_\scol (i-2)}\]
\caption{Exemplary setting of the proof of \cref{thm:qh}}
\label{fig:prf}
\end{figure}
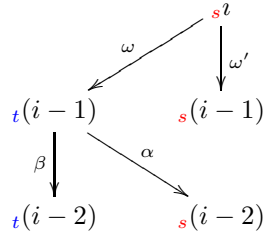

For $x=e$, i.e. $i=0$, there is nothing to show since $P(x)=\Delta(x)$ (by comparing dimensions via \cref{lem:dimmod}). For $i=1$, that is $x=\scol$, we have the short exact sequence
\[\xymatrix{0 \ar[r] & P(e)\langle -1\rangle \ar[r]^f & P(\scol)\ar[r] & \Delta(\scol) \ar[r]& 0}\]
where $f$ is pre-composing with $(\scol,e)$. It is clear that $\op{im} f \subseteq M(\scol)$ and the equality holds by comparing dimensions (by \cref{lem:dimmod}).

Now let $i\geq 2$. Pre-composing with $\omega^\prime:=({}_\scol i,{}_\scol (i-1))$ (as indicated in \cref{fig:prf}) gives us
\begin{align}
\xymatrix{f:\; P({}_\scol (i-1))\langle -1 \rangle \ar[r] & P({}_\scol i)}
\end{align}
which turns out to be injective since basis elements are mapped to pairwise non-equivalent paths . Clearly $\pi f= 0$, hence using $\ref{eqn:ses}$ we obtain the following commutative diagram

\[
\xymatrix{
M({}_\scol i) \ar[r] & P({}_\scol i) \ar[r]^\pi &   \Delta({}_\scol i)\\
&P({}_\scol (i-1))\langle -1 \rangle \ar[u]^{f} \ar@{-->}[ul]^{\exists ! \; \iota} \ar[ur]_0 & }
\]

Consider the composition
\[g:\xymatrix{P({}_\tcol (i-1))\langle -1 \rangle \ar[r]^/.8em/h & M({}_\scol i ) \ar[r]& M({}_\scol i )/\op{im}\,\iota }\]
where $h$ is pre-composing with $\omega := ({}_\scol i, {}_\tcol (i-1))$. The module $M({}_\scol i )$ has generators $\omega$ and $\omega^\prime$ as a left module. We can conclude that $g$ is surjective since $\omega \in \op{im}\; h$ and $\omega^\prime\in \op{im}\;\iota$. Similarly, the module $M({}_\tcol (i-1))$ has generators $\alpha := ({}_\tcol (i-1),{}_\scol (i-2))$ and $\beta := ({}_\tcol (i-1),{}_\tcol (i-2))$. Using the relations \ref{eqn:12} and \ref{eqn:13} we obtain the identities:
\begin{align*}
\alpha \omega \enspace & =  \enspace ({}_\scol i,{}_\tcol (i-1),{}_\scol (i-2)) = ({}_\scol i,{}_\scol (i-1),{}_\scol (i-2)), \\
\beta\omega \enspace & =  \enspace ({}_\scol i,{}_\tcol (i-1),{}_\tcol (i-2)) = ({}_\scol i,{}_\scol (i-1),{}_\tcol (i-2)).
\end{align*}
Hence we can deduce that $M({}_\tcol (i-1))\langle -1\rangle \subseteq \op{ker}\; g$ and therefore

\[\xymatrix{P({}_\tcol (i-1))\langle -1 \rangle \ar@{->>}[r]^g \ar@{->>}[d] &  M({}_\scol i )/\op{im}\;\iota  \\
\Delta({}_\tcol (i-1))\langle -1 \rangle. \ar@{-->}[ru]_{\exists ! \; \varphi}& }\]

Since $g$ is surjective $\varphi$ is surjective, too. Using \cref{lem:dimmod} we can compare dimensions of both sides and deduce that $\varphi$ is an isomorphism. Thus we have the short exact sequence:
\[ \xymatrix{0\ar[r] & P({}_\scol (i-1))\langle -1 \rangle \ar[r] & M({}_\scol i)\ar[r] & \Delta({}_\tcol (i-1))\langle -1 \rangle \ar[r] & 0}.\]
By induction hypothesis $P({}_\scol (i-1))$ has $\Delta$-filtration with subquotients isomorphic to $\Delta(y)$ with $y \geq_r {}_\scol (i-1)$ (each with multiplicity 1). Thus $M({}_\scol i)$ has $\Delta$-filtration with subquotients isomorphic to $\Delta(y)$ with $y >_r {}_\scol i$ and therefore $P({}_\scol i)$ has the desired $\Delta$-filtration.

For $x \in W$ there clearly exists a surjective map $\pi^\prime: \Delta(x)\to L(x)$ and by construction $\op{ker}\;\pi^\prime$ has only composition factors isomorphic to $L(\mu_i)$ where $\mu_i \leq_r x$. Since every non-trivial path from $x$ to $x$ passes through a vertex $y >_r x$ we have $\left[\Delta(x):L(x) \right]=1$ and hence the second condition is satisfied.
\end{prf}

\begin{thm}\label{thm:stdkoszul}
$\textbf P_m$ is Standard Koszul, i.e. every left (resp. right) standard module admits a linear resolution.
\end{thm}

\begin{prf}
By \cref{prop:op} we have an isomorphism $\textbf P_m\cong\textbf P_m^{op}$ and thus left and right modules can be identified. Therefore it is enough to show that the (left) standard modules admit such resolutions. Since $e\in W$ is maximal we have $\Delta(e)=P(e)$ and there is nothing to show.
For $x \in W\setminus \{e\}$ we construct a linear resolution $P_\bullet(x)\stackrel{\epsilon}{\to} \Delta(x) \to 0$ with $P_\bullet(x):=(P_i)_{i\geq 0}$ defined by:
{\small \[P_i := \left( \bigoplus_{\substack{w\geq_r x \\ \ell(x)-\ell(w) = i}}P(w)\right)\langle -i \rangle.\]}
Note that $P_0=P(x)$. The augmentation map $\epsilon:P_0 \to \Delta(x)$ is just the canonical projection $\pi: P(x) \to \Delta(x)$ (see \ref{eqn:ses}). For $x={}_\scol \ell$ the boundary maps $p_i:P_i \to P_{i-1}$ for $i\geq 1$ are defined as follows:

\begin{small}
\[p_i \enspace = \enspace
\begin{cases}
\cdot \begin{pmatrix}
({}_\scol \ell ,{}_\scol (\ell -1)) & ({}_\scol \ell ,{}_\tcol (\ell -1))
\end{pmatrix}& \text{if } i=1< \ell,\\
\cdot \begin{pmatrix} (\scol,e) \end{pmatrix} & \text{if } i=1=\ell,\\
\cdot \begin{pmatrix}
\left( {}_\scol(\ell-i+1),{}_\scol(\ell-i)\right) & (-1)^{i+1}\left( {}_\scol(\ell-i+1),{}_\tcol(\ell-i)\right)\\
(-1)^{i+1} \left( {}_\tcol(\ell-i+1),{}_\scol(\ell-i)\right) & \left( {}_\tcol(\ell-i+1),{}_\tcol(\ell-i)\right)
\end{pmatrix}& \text{if } 2 \leq i < \ell, \\
\cdot \begin{pmatrix}
({}\scol,e)\\
(-1)^{i+1} (\tcol,e)
\end{pmatrix}& \text{if } i= \ell,\\
0 & \text{if } i > \ell.
\end{cases}\]
\end{small}

Recall that we compose paths from right to left (similar as morphisms) and hence the defined $p_i$'s are pre-composing with certain arrows. All indecomposable projective modules $P(w)$ are generated by the corresponding idempotent $e_w$ which lies in degree $0$. Thus by construction each $P_i$ is projective and generated by its degree $i$ component. Clearly, the augmentation map $\epsilon$ is surjecitve.

The fact that for fixed $x$ the above sequence $P_\bullet (x)$ is a complex is a direct consequence of relations \ref{eqn:12} and \ref{eqn:13}. We only present the generic case $2\leq i < \ell(x)$. Each $P(w)$ is generated by $e_w$, so it is enough to check that $p^2(e_w)=0$. Let $w={}_\scol (i-2))$ then we have
\begin{align*}
p(e_w) \enspace = \enspace & \begin{pmatrix} ({}_\scol(i-1),{}_\scol (i-2))\\ (-1)^{*}({}_\tcol(i-1),{}_\scol (i-2)) \end{pmatrix} \\
p^2(e_w) \enspace = \enspace & \begin{pmatrix} ({}_\scol i,{}_\scol(i-1),{}_\scol (i-2)) + (-1)^{*+\bullet}({}_\scol i,{}_\tcol(i-1),{}_\scol (i-2)) \\ (-1)^{\bullet}({}_\tcol i,{}_\scol(i-1),{}_\scol (i-2)) + (-1)^{*}({}_\tcol i,{}_\tcol(i-1),{}_\scol (i-2))
\end{pmatrix} 
\end{align*}

where $*$ and $\bullet$ depend on the parity of $\ell - i$ (resp. $\ell - i -1$). Nevertheless, the parity of $*$ and $\bullet$ is always different and therefore $(-1)^{*+\bullet}=-1$ and $(-1)^* = -(-1)^\bullet$. Thus by relations \ref{eqn:12} and \ref{eqn:13} we can deduce that $p^2(e_w)=0$. The other remaining cases are similar. Hence the sequence is a complex and it remains to show that this complex is exact. 

We only show the generic case $2\leq i \leq \ell(x)-2$ which covers all important arguments. Assume we have $(y,x)\in\op{ker} \;p_i$ where $p_i: P_i \to P_{i-1}$. Without loss of generality 
\begin{alignat*}{2}
P_i & = \left[P({}_\scol j)\bigoplus P({}_\tcol j)\right]\langle -i \rangle, \quad\quad &\quad\quad
P_{i-1} & = \left[P({}_\scol (j+1))\bigoplus P({}_\tcol (j+1))\right]\langle -i+1\rangle 
\end{alignat*} 
for some $2 \leq j \leq \ell(x)-2$. From the definition of the $P_i$ we know that $y$ (resp. $x$) is a path starting in ${}_\scol j$ (resp. ${}_\tcol j$). By assumption we have $p_i((y,x))=0$ and in particular $\pi_1 (p_i((y,x)))=0 \in P({}_\scol (j-1))\langle -i+1\rangle$. Recall that $\pi_1p_i$ is pre-composing with $\omega^\prime=({}_\scol(j+1)),{}_\scol j)$ in the first component and pre-composing with $\omega=({}_\scol(j+1),{}_\tcol j)$ in the second component. Since $(y,x)\in\op{ker}\;p_i$ we can conclude that they must end in the same vertex, say $z$ (see \cref{fig:prfstd}). Note that it is not possible for either $x$ or $y$ to be the trivial path because then the other path would have to have length $\geq 2$ and hence they would be linearly independent. Assuming both $x$ and $y$ pass only through vertices $v$ such that $v\leq_r {}_\scol j$ then $z\notin \{{}_\scol j,{}_\tcol j\}$ (otherwise one path would pass through a vertex $v$ such that $v>_r{}_\scol j$). Moreover, $x\omega$ (resp. $y\omega^\prime$) does not pass through ${}_\scol j$ (resp. ${}_\tcol j$). In particular, the highest vertex $x\omega$ (resp. $y\omega^\prime$) passes through is ${}_\tcol j$ (resp. ${}_\scol j$). Thus $x\omega$ and $y\omega^\prime$ are linearly independent by \cref{rem:paths} (which is formulated in terms of the Bruhat order). This is a contradiction to $(y,x)\in\op{ker}\;p_i$. If only one path passes through a vertex $v$ with $v >_r {}_\scol j$ the same argumentation using \cref{rem:paths} again yields a contradiction to $(y,x)\in\op{ker}\;p_i$. Therefore $x$ and $y$ pass only through vertices $v$ such that $v>_r {}_\scol j$. Hence we can deduce that $y$ starts either with $\alpha$ or $\beta$ and $x$ starts either with $\gamma$ or $\delta$. But that means that $(y,x)\in\op{im}\;p_{i+1}$. Thus the complex is exact and indeed a linear resolution.

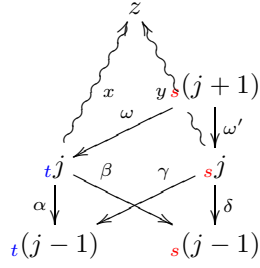
\begin{figure}
\[\xymatrix@!=0.5pc{&z &\\
&& {}_\scol (j+1) \ar[dll]_\omega \ar[d]^{\omega^\prime}\\
{}_\tcol j \ar@{~>}[uur]_x \ar[d]_\alpha\ar[drr]^/-1.5em/\beta && {}_\scol j\ar@{~>}[luu]^y|!{[ll];[u]}\hole \ar[d]^\delta\ar[lld]_/-1.5em/\gamma\\
{}_\tcol (j-1) && {}_\scol (j-1)} \]
\caption{Exemplary setting of the proof of \cref{thm:stdkoszul}}
\label{fig:prfstd}
\end{figure}
\end{prf}

By \cref{thm:qhk} we obtain the following corollary:
\begin{cor}\label{cor:koszul}
$\textbf P_m$ is Koszul.
\end{cor}

\subsection{The quadratic dual of $\Pm$}
By construction $\Pm = \DR Q_m/(R_m^\mathfrak h)$. Since $A:=\DR Q_m = \bigoplus_{i\geq 0} A^i$ is graded by the path length and $R_m^\mathfrak h\subseteq A^2$ is homogeneous $\Pm$ inherits the grading by the path length. We note that $A = T_{A^0}A^1 = \bigoplus_{i \geq 0}A^i$ and thus $\Pm=T_{A^0} V /(R_m^\mathfrak h)$ is the desired quadratic structure for the $A^0$-bimodule $V = A^1$. $V$ is a vector space spanned by the edges in $Q_m$ on which we can define the standard scalar product
\[\langle \alpha,\beta \rangle := 
\begin{cases}
1 & \text{if  } t(\alpha)= s(\beta) \text{ and } t(\beta)=s(\alpha)\\
0 & \text{otherwise}
\end{cases}\]
for edges $\alpha = (s(\alpha),t(\alpha))$ and $\beta = (s(\beta),t(\beta))$. Using this scalar product we can identify $V\cong V ^*$. In the construction of the quadratic dual we have to consider 
\[R^\perp := \{x \mid \forall v\in R_m^\mathfrak h:\;\langle v,x \rangle =0  \}\subseteq V^*\otimes_{A^0} V^* \cong (V \otimes_{A^0} V)^*.\]
By \cref{lem:dim} the algebra $\Pm$ is finite-dimensional  and with the above identification $R^\perp$ is just the usual orthogonal complement of $R_m^\mathfrak h$ inside $A^2$ with respect to the standard scalar product. The vector space $A^2$ has a basis consisting of all paths of length $2$ and therefore it is of the form $A^2= \bigoplus_{x,y\in W} e_yA^2e_x$. Taking duals commutes with finite direct sums, so for fixed $x,y\in W$ the complement of all relations starting in $x$ and ending in $y$ is inside $e_yA^2e_x$ which has at most $\DR$-dimension $4$. With the above considerations it follows that $\Pm^! = T_{A^0}V/(R^\perp)$.

Recall that all defining relations in $\Pm$ are homogeneous of degree $2$ and can be interpreted as linear combinations of paths of length $2$. Using this identification there are three types of relations in $\Pm$: paths with the same starting and terminal point, paths with different starting and terminal points between vertices of the same length and paths with different starting and terminal points between vertices of different lengths. We only present the cases where the starting point $x$ is of the form ${}_\scol i$ for $0\leq i \leq m$. The cases where $x={}_\tcol i$ can be treated similarly since all relations are symmetric in $\scol$ and $\tcol$. For the sake of simplicity define $\overline i := m-i$ for $0\leq i \leq m$. 

\subsubsection{Relations \ref{eqn:1} - \ref{eqn:5}}
These relations have in common that the starting point and terminal point coincide, say $x$. For fixed $x$ all calculations in this section take place in $A(x):=e_xA^2e_x$. We have to distinguish five different cases:
\begin{itemize}
\item $x=e$: We have $\op{dim}_\DR A(e)=2$ and both basis elements $(e,\scol,e)$ and $(e,\tcol,e)$ are in $R_m^\mathfrak h$, hence the orthogonal complement is $0$ and there are no orthogonal relations.
\item $x=\scol$: The vector space $A(\scol)$ has a basis consisting of:
\begin{alignat*}{3}
u_1^\scol & := (\scol,\scol\tcol,\scol), \quad \quad & \quad \quad u_2^\scol &:=(\scol,\tcol\scol,\scol), \quad \quad & \quad \quad u_3^\scol & := (\scol,e,\scol).
\end{alignat*}
The relations are $u_1^\scol=0$ and $u_2^\scol + [2]_\zeta u_3^\scol=0$, thus the orthogonal complement is spanned by $[2]_\zeta u_2^\scol - u_3^\scol$.
\item $x={}_\scol i$ for some $2\leq i \leq m-2$: The vector space $A({}_\scol i)$ has a basis consisting of:
\begin{alignat*}{2}
v_1^{\scol,i} &:= ({}_\scol i,{}_\scol (i+1),{}_\scol i),\quad  & \quad  v_2^{\scol,i} &:= ({}_\scol i,{}_\tcol (i+1),{}_\scol i),\\
v_3^{\scol,i} &:= ({}_\scol i,{}_\scol (i-1),{}_\scol i), \quad  & \quad  v_4^{\scol,i} &:= ({}_\scol i,{}_\tcol (i-1),{}_\scol i).
\end{alignat*}
The relations are $v_1^{\scol,i} - \lambda_i v_3^{\scol,i} = 0$ and $v_2^{\scol,i}-\mu_i v_3^{\scol,i} + \nu_i v_4^{\scol,i}=0$ where 
\[\lambda_i:= \frac{[i-1]_\zeta}{[i]_\zeta},\quad\quad\quad \mu_i:= [i-1]_\zeta-[i+1]_\zeta,\quad\quad\quad \nu_i:= \frac{[i+1]_\zeta}{[i]_\zeta}.\] 
The orthogonal complement is spanned by $\lambda_i v_1^{\scol,i} + \mu_i v_2^{\scol,i} + v_3^{\scol,i}$ and $\nu_i v_2^{\scol,i} -v_4^{\scol,i}$.
\item $x={}_\scol (m-1)={}_\scol \overline{1}$: The vector space $A({}_\scol \overline{1})$ has a basis consisting of:
\begin{alignat*}{3}
w_1^\scol & := ({}_\scol \overline{1},w_0,{}_\scol \overline{1}), \quad  & \quad w_2^\scol & :=({}_\scol \overline{1},{}_\scol \overline{2},{}_\scol \overline{1}),\quad &\quad
w_3^\scol & := ({}_\scol \overline{1},{}_\tcol\overline{2},{}_\scol \overline{1}).
\end{alignat*}
The relation is $0=w_1^\scol-\frac{[m-2]_\zeta}{[m-1]_\zeta}w_2^\scol =w_1^\scol -[2]_\zeta w_2^\scol$ (using \ref{eqn:q7}). Therefore the orthogonal complement is spanned by $[2]_\zeta w_1^\scol +w_2^\scol$ and $w_3^\scol$.
\item $x=w_0$: We have $\op{dim}_\DR A(w_0)=2$ and there are no relations, hence the orthogonal complement is the complete subspace with basis $(w_0,{}_\scol \overline{1},w_0)$ and $(w_0,{}_\tcol \overline{1},w_0)$.
\end{itemize}

\subsubsection{Relations \ref{eqn:6} - \ref{eqn:9}}
All these relations are paths from $x$ to $y$ such that $x\neq y$ and $\ell(x)=\ell(y)$. All computations in this section take place in $A(x,y) := e_yA^2e_x$, i.e. the vector subspace consisting of all paths of length $2$ from $x$ to $y$. For these relations there are three cases to consider:

\begin{itemize}
\item $x=\scol$ and $y=\tcol$: The vector space $A(\scol,\tcol)$ has a basis consisting of:
\begin{alignat*}{3}
\overline{u}_1^\scol & := (\scol,\scol\tcol,\tcol), \quad \quad & \quad \quad \overline{u}_2^\scol &:=(\scol,\tcol\scol,\tcol), \quad \quad & \quad \quad \overline{u}_3^\scol & := (\scol,e,\tcol).
\end{alignat*}
The relations are $\overline{u}_1^\scol-\overline{u}_3^\scol=0$ and $\overline{u}_2^\scol-\overline{u}_3^\scol=0$. Thus the orthogonal complement is spanned by $\overline{u}_1^\scol +\overline{u}_2^\scol+ \overline{u}_3^\scol$.
\item $x={}_\scol i$ and $y={}_\tcol i$ for some $2\leq i \leq m-2$: The vector space $A({}_\scol i,{}_\tcol i)$ has a basis consisting of:
\begin{alignat*}{2}
\overline{v}_1^{\scol,i} &:= ({}_\scol i,{}_\scol (i+1),{}_\tcol i),\quad \quad  & \quad \quad \overline{v}_2^{\scol,i} &:= ({}_\scol i,{}_\tcol (i+1),{}_\tcol i),\\
\overline{v}_3^{\scol,i} &:= ({}_\scol i,{}_\scol (i-1),{}_\tcol i),\quad \quad  & \quad \quad \overline{v}_4^{\scol,i}&:= ({}_\scol i,{}_\tcol (i-1),{}_\tcol i).
\end{alignat*}
The relations are $\overline{v}_1^{\scol,i} - \alpha_i \overline{v}_3^{\scol,i} -\overline{v}_4^{\scol,i}=0$ and $\overline{v}_2^{\scol,i}-\overline{v}_3^{\scol,i}-\alpha_i\overline{v}_4^{\scol,i}=0$ where $\alpha_i:= \frac{1}{[i]_\zeta}$. Hence the orthogonal complement is spanned by $\alpha_i \overline{v}_1^{\scol,i} + \overline{v}_2^{\scol,i}+\overline{v}_3^{\scol,i}$ and $\overline{v}_1^{\scol,i}+\alpha_i \overline{v}_2^{\scol,i}+ \overline{v}_4^{\scol,i}$.
\item $x={}_\scol \overline{1}$ and $y={}_\tcol \overline{1}$: The vector space $A({}_\scol \overline{1},{}_\tcol \overline{1})$ has a basis consisting of:
\begin{alignat*}{3}
\overline{w}_1^\scol & := ({}_\scol \overline{1},w_0,{}_\tcol \overline{1}), \quad & \quad \overline{w}_2^\scol & :=({}_\scol \overline{1},{}_\scol \overline{2},{}_\tcol \overline{1}),\quad &\quad
\overline{w}_3^\scol& := ({}_\scol \overline{1},{}_\tcol \overline{2},{}_\tcol \overline{1}).
\end{alignat*}
The relation is $\overline{w}_1^\scol -\overline{w}_2^\scol- \overline{w}_3^\scol$  (using \ref{eqn:q7}) and therefore the orthogonal complement is spanned by $\overline{w}_1^\scol+\overline{w}_2^\scol$ and $\overline{w}_1^\scol + \overline{w}_3^\scol$.
\end{itemize}

\subsubsection{Relations \ref{eqn:10} - \ref{eqn:13}}\label{sub:rel10to13}
For the relations for which the starting point $x$ and the terminal point $y$ differ in their length, it is immediately clear that $|\ell(x)-\ell(y)|=2$. Hence, $\op{dim}_\DR A(x,y)= 2$, say with basis $\{p_1,p_2\}$. The relations in \ref{eqn:10}- \ref{eqn:13} are then $p_1-p_2 = 0$ and therefore their complements are $p_1 + p_2 = 0$. 

\begin{small}
\begin{center}
\begin{table}
$\begin{array}{c|c}
\text{Relations } & \text{Orthogonal Relations} \\
\hline
\begin{matrix}
(e,\scol,e)\\
(e,\tcol,e)
\end{matrix} & n/a\\
\hline
\begin{matrix}
u_1^\scol\\
u_2^\scol + [2]_\zeta u_3^\scol
\end{matrix} & [2]_\zeta u_2^\scol - u_3^\scol\\
\hline
\begin{matrix}
v_1^{\scol,i} - \lambda_i v_3^{\scol,i}\\
v_2^{\scol,i}-\mu_i v_3^{\scol,i} + \nu_i v_4^{\scol,i}
\end{matrix} & 
\begin{matrix}
\lambda_i v_1^{\scol,i} + \mu_i v_2^{\scol,i} + v_3^{\scol,i} \\
\nu_i v_2^{\scol,i} - v_4^{\scol,i}
\end{matrix}\\
\hline
w_1^\scol -[2]_\zeta w_2^\scol &
\begin{matrix} 
[2]_\zeta w_1^\scol +w_2^\scol \\ w_3^\scol
\end{matrix} \\
\hline
n/a & \begin{matrix}
(w_0,{}_\scol (m-1),w_0)\\
(w_0,{}_\tcol (m-1),w_0)
\end{matrix} \\
\hline
\begin{matrix}
\overline{u}_1^\scol-\overline{u}_3^\scol\\
\overline{u}_2^\scol-\overline{u}_3^\scol
\end{matrix} & \overline{u}_1^\scol+\overline{u}_2^\scol+\overline{u}_3^\scol\\
\hline
\begin{matrix}
\overline{v}_1^{\scol,i} - \alpha_i \overline{v}_3^{\scol,i}  -\overline{v}_4^{\scol,i} \\
\overline{v}_2^{\scol,i} -\overline{v}_3^{\scol,i} -\alpha_i\overline{v}_4^{\scol,i} 
\end{matrix} & 
\begin{matrix}
\alpha_i \overline{v}_1^{\scol,i}  + \overline{v}_2^{\scol,i} +\overline{v}_3^{\scol,i} \\
\overline{v}_1^{\scol,i} +\alpha_i \overline{v}_2^{\scol,i} + \overline{v}_4^{\scol,i} 
\end{matrix}\\
\hline
\overline{w}_1^{\scol} -\overline{w}_2^{\scol} -\overline{w}_3^{\scol}  &
\begin{matrix} 
\overline{w}_1^{\scol} +\overline{w}_2^{\scol}  \\\overline{w}_1^{\scol} +\overline{w}_3^{\scol} 
\end{matrix}\\
\hline
({}_\scol j,{}_\scol (j+1),{}_\scol(j+2))- ({}_\scol j,{}_\tcol (j+1),{}_\scol(j+2)) & ({}_\scol j,{}_\scol (j+1),{}_\scol(j+2))+ ({}_\scol j,{}_\tcol (j+1),{}_\scol(j+2))\\
({}_\scol j, {}_\scol (j+1) ,{}_\tcol(j+2))-  ({}_\scol j,{}_\tcol (j+1) ,{}_\tcol(j+2)) &({}_\scol j, {}_\scol (j+1) ,{}_\tcol(j+2)) +({}_\scol j,{}_\tcol (j+1) ,{}_\tcol(j+2))\\
({}_\scol(j+2), {}_\scol (j+1) ,{}_\scol j) -  ({}_\scol(j+2), {}_\tcol (j+1) ,{}_\scol j) &({}_\scol(j+2), {}_\scol (j+1), {}_\scol j) + ({}_\scol(j+2), {}_\tcol (j+1) ,{}_\scol j) \\
({}_\scol(j+2), {}_\scol (j+1),{}_\tcol j) -  ({}_\scol(j+2), {}_\tcol (j+1) ,{}_\tcol j) &({}_\scol(j+2), {}_\scol (j+1),{}_\tcol j) + ({}_\scol(j+2), {}_\tcol (j+1) ,{}_\tcol j)
\end{array}$ 
\caption{Relations and its orthogonal relations}\label{tab:typeall}
\end{table}
\end{center}
\end{small}

Table \ref{tab:typeall} summarises all relations and orthogonal relations. Define $R^\perp \subseteq A^2$ be the set of all orthogonal relations from Tables \ref{tab:typeall}. Then as aforementioned we have $\Pm^! = A/(R^\perp)$. Mimicking the proof of  \cref{lem:dim} we obtain the following lemma:
\begin{lem}\label{lem:dimquadraticdual}
The quadratic dual $\Pm^!$ has the same dimension as $\Pm$ as $\DR$-vector space.
\end{lem}

\subsection{Self-duality}
In this section we prove the following theorem 

\begin{thm}\label{thm:selfduality}
The algebra $\textbf P_m$ is Koszul self-dual, i.e. $\Pm \cong E(\Pm) \cong \textbf P_m^!$. 
\end{thm}

\begin{prf}
The algebra $\textbf P_m$ is finite-dimenesional and Koszul by \cref{cor:koszul}. Thus there is a canonical isomorphism $E(\Pm)\cong (\Pm^!)^{op}$. Since $\Pm \cong \Pm^{op}$ by \cref{prop:op} it is enough to show that $\Pm \cong \Pm^! = E(\Pm)$.

Define the map $\Theta$ on the vertices of $Q_m$ by $x\mapsto x^{-1}w_0$ which extends to an $A^0$-bimodule homomorphism $\Theta: A \to A$ for $A=\DR Q_m$ (see \cref{fig:hasseQm}). Note that the images of the arrows are only determined up to scalars which we choose as $\pm 1$ as indicated in \cref{fig:scalarsprf}. Note that the scalar of an arrow only depends on the adjacent vertices and is independent of the direction. The pattern is highly regular except for the scalar at the edge $({}_\tcol (m-1),w_0)$. Since all scalars are invertible, the map $\Theta$ is an isomorphism and thus we have a surjection $\Phi: A \stackrel{\Theta}{\longrightarrow} A \twoheadrightarrow \Pm^!$. It suffices to show that $R_m^\mathfrak h \subseteq \op{ker}\,\Phi$ (or equivalently: $\Theta(R_m^\mathfrak h)\subseteq (R^\perp)$) which implies the existence of a surjection $\Psi: \Pm \twoheadrightarrow \Pm^!$, which is an isomorphism by dimension reasons (see \cref{lem:dimquadraticdual}).
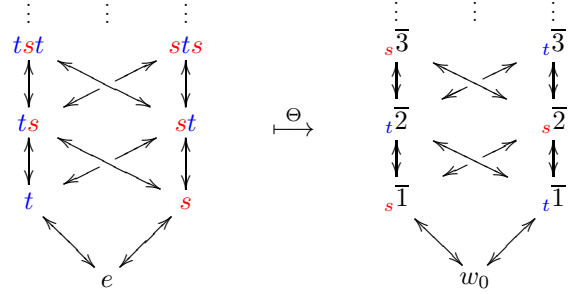
\begin{figure}
\begin{align*}
\begin{minipage}[c]{1cm}
\xymatrix@!=0.45pc{\ar@{{}{ }{}}[d]|\vdots &\ar@{{}{ }{}}[d]|\vdots& \ar@{{}{ }{}}[d]|\vdots\\
\quad \tcol\scol\tcol \quad \ar@{<->}[d]\ar@{<->}[rrd]& &\quad \scol\tcol\scol \quad\ar@{<->}[lld]|!{[ll];[d]}\hole\ar@{<->}[d]\\
\quad\tcol\scol\quad  \ar@{<->}[d]\ar@{<->}[rrd]&&\quad\scol\tcol\quad \ar@{<->}[lld]|!{[ll];[d]}\hole\ar@{<->}[d]\\
\quad \tcol \quad\ar@{<->}[dr]&& \quad \scol \quad\ar@{<->}[dl]\\
& e& }
\end{minipage} \quad
\stackrel{\Theta}{\longmapsto} \quad
\begin{minipage}[c]{1cm}
\xymatrix@!=0.45pc{\ar@{{}{ }{}}[d]|\vdots &\ar@{{}{ }{}}[d]|\vdots& \ar@{{}{ }{}}[d]|\vdots\\
\quad {}_\scol\overline{3} \quad \ar@{<->}[d]\ar@{<->}[rrd]&& \quad {}_\tcol\overline{3}\quad\ar@{<->}[lld]|!{[ll];[d]}\hole\ar@{<->}[d]\\
\quad {}_\tcol\overline{2} \quad  \ar@{<->}[d]\ar@{<->}[rrd]&&\quad {}_\scol\overline{2} \quad \ar@{<->}[lld]|!{[ll];[d]}\hole\ar@{<->}[d]\\
\quad {}_\scol\overline{1} \quad\ar@{<->}[dr]&& \quad {}_\tcol\overline{1} \quad\ar@{<->}[dl]\\
& w_0&}
\end{minipage}
\end{align*}
\caption{The Hasse graph $Q_m$ and its image under $\Theta$}\label{fig:hasseQm}
\end{figure}

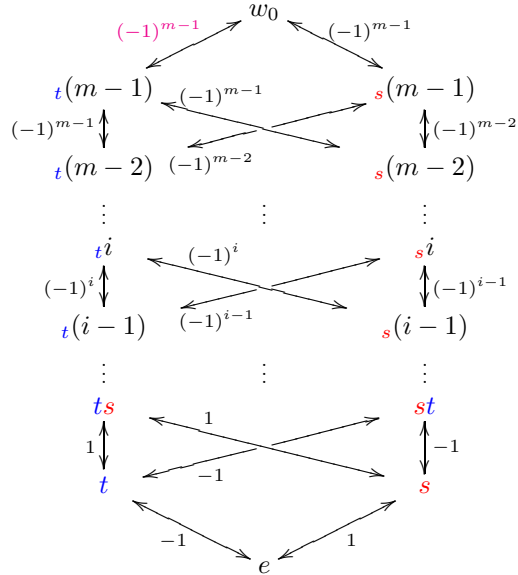
\begin{figure}
\[\xymatrix@!=0.5pc{&&w_0\ar@{<->}[lld]_{\textcolor{magenta}{(-1)^{m-1}}}\ar@{<->}[rrd]^{(-1)^{m-1}}&& \\
{}_\tcol (m-1)\ar@{<->}[d]_{(-1)^{m-1}}\ar@{<->}[rrrrd]^/-0.75cm/{(-1)^{m-1}}&&&&{}_\scol (m-1)\ar@{<->}[lllld]^/0.9cm/{(-1)^{m-2}}|!{[llll];[d]}\hole\ar@{<->}[d]^{(-1)^{m-2}}\\
\quad {}_\tcol (m-2) \quad \ar@{{}{ }{}}[d]|\vdots\ar@{{}{ }{}}[rrrrd]|\vdots&&&& \quad {}_\scol (m-2)\quad \ar@{{}{ }{}}[d]|\vdots \\
\quad {}_\tcol i \quad\ar@{<->}[d]_{(-1)^i}\ar@{<->}[rrrrd]^/-0.8cm/{(-1)^i}&&&& \quad {}_\scol i \quad\ar@{<->}[d]^{(-1)^{i-1}}\ar@{<->}[lllld]^/0.8cm/{(-1)^{i-1}}|!{[llll];[d]}\hole\\
\quad {}_\tcol (i-1) \quad\ar@{{}{ }{}}[d]|\vdots \ar@{{}{ }{}}[rrrrd]|\vdots &&&& \quad {}_\scol (i-1) \quad\ar@{{}{ }{}}[d]|\vdots\\
\quad \tcol\scol\quad\ar@{<->}[d]_1\ar@{<->}[rrrrd]^/-0.8cm/1&&&& \quad \scol\tcol \quad\ar@{<->}[d]^{-1}\ar@{<->}[lllld]^/0.8cm/{-1}|!{[llll];[d]}\hole\\
\quad \tcol \quad\ar@{<->}[drr]_{-1}&& &&\quad \scol \quad\ar@{<->}[dll]^1\\
&& e&& } \]
\caption{A choice of scalars yielding an isomorphism $\Pm \cong \Pm^!$}\label{fig:scalarsprf}
\end{figure}

In order to show that $R_m^\mathfrak h \subseteq \op{ker}\;\Phi$ we check that each relation from $R_m^\mathfrak h$ is mapped to a relation from $R^\perp$ (up to a sign). The map $\Theta: W \to W$ is an isomorphism and for all $w\in W$ we have $\ell(\Theta (w))=m-\ell(w)$. Therefore we can treat the three types from above separately, again. Note that 
\begin{align}\label{eqn:map}
\Theta({}_\scol i) \enspace & = \enspace
\begin{cases}
{}_\scol (m-i) = {}_\scol \overline i & \text{if } i \text{ is even,}\\
{}_\tcol (m-i) = {}_\tcol \overline i & \text{if } i \text{ is odd.}\\
\end{cases} 
\end{align}
Recall the identity in \ref{eqn:q7} for the quantum numbers 
$[i]_\zeta \enspace  = \enspace [m-i]_\zeta \enspace = \enspace \left[\,\overline{i}\,\right]_\zeta$
when $q$ is specialised to $\zeta=e^{2\pi i/2m}$. Therefore there are the following identities for $\alpha_i,\lambda_i,\mu_i$ and $\nu_i$:
\begin{alignat}{3}\label{eqn:scalars}
\alpha_{\overline i} \enspace & = \enspace \alpha_i, \quad\quad & \quad \quad
\lambda_{\overline i} \enspace & = \enspace \nu_i, \quad\quad & \quad \quad
\mu_{\overline{ i}} \enspace & = \enspace - \mu_i.
\end{alignat}

Since the calculations are not difficult but tedious we only compute the images of the relations \ref{eqn:1} - \ref{eqn:5} under $\Theta$ using \ref{eqn:map}. Note that all basis elements are loops and therefore all occurring scalars are $1=(-1)^2$ since the scalar is independent of the direction of the arrow.

\begin{small}
\begin{align*}
(e,\scol,e) \enspace & \longmapsto \enspace (w_0,{}_\tcol (m-1),w_0), \\
(e,\tcol,e) \enspace & \longmapsto \enspace (w_0,{}_\scol (m-1),w_0),\\
u_1^\scol \enspace & \longmapsto \enspace w_3^\tcol, \\
u_2^\scol + [2]_\zeta u_3^\scol \enspace & \longmapsto \enspace w_2^\tcol + [2]_\zeta w_1^\tcol,\\
v_1^{\scol,i} - \lambda_i v_3^{\scol,i} \enspace & \longmapsto \enspace  
\begin{cases}
v_4^{\scol,\overline i} - \lambda_i v_2^{\scol,\overline i} &\stackrel{\ref{eqn:scalars}}{=} v_4^{\scol,\overline i} - \nu_{\overline i} v_2^{\scol,\overline i} \quad  \text{if } i \text{ is even,}\\
v_4^{\tcol,\overline i} - \lambda_i v_2^{\tcol,\overline i} &\stackrel{\ref{eqn:scalars}}{=} v_4^{\tcol,\overline i} - \nu_{\overline i} v_2^{\tcol,\overline i} \quad\, \text{if } i \text{ is odd,}
\end{cases}\\
v_2^{\scol,i} -\mu_i v_3^{\scol,i} + \nu_i v_4^{\scol,i} \enspace & \longmapsto \enspace 
\begin{cases}
v_3^{\scol,\overline i} - \mu_i v_2^{\scol,\overline i} + \nu_i  v_1^{\scol,\overline i} &\stackrel{\ref{eqn:scalars}}{=}
v_3^{\scol,\overline i}  + \mu_{\overline i} v_2^{\scol,\overline i} + \lambda_{\overline i}  v_1^{\scol,\overline i}\quad \text{if } i \text{ is even,}\\
v_3^{\tcol,\overline i} -\mu_i v_2^{\tcol,\overline i} + \nu_i v_1^{\tcol,\overline i}  &\stackrel{\ref{eqn:scalars}}{=}
v_3^{\tcol,\overline i} +\mu_{\overline i} v_2^{\tcol,\overline i} + \lambda_{\overline i} v_1^{\tcol,\overline i}  \quad\, \text{if } i \text{ is odd,}
\end{cases}\\
w_1^\scol -[2]_\zeta w_2^\scol \enspace & \longmapsto \enspace 
\begin{cases}
u_3^\tcol -[2]_\zeta u_2^\tcol & \text{if } m \text{ is even,}\\
u_3^\scol -[2]_\zeta u_2^\scol & \text{if } m \text{ is odd,}
\end{cases} 
\end{align*}
\end{small}

The images of the relations agree with the right-hand side in \cref{tab:typeall} up to a sign, swapping the roles of $\scol$ and $\tcol$ and replacing $i$ with $\overline i$. Therefore their image is contained in $(R^\perp)$. The other relations can be treated analogously as above.
\end{prf}

\begin{rem}
There is a certain degree of freedom in the choice of the scalars for this isomorphism. In types $A_2$ and $B_2$ there are $5$ scalars which can be chosen freely. However, solving such a system of quadratic equations gets more and more complicated and is not very enlightening.
\end{rem}

\begin{rem}
The basis for our approach in this paper was the equivalence of the diagrammatic category $\mathcal D$ and the category of Soergel bimodules $\sbim$ (see \cref{thm:equiv}). However, this equivalence holds in particular for all finite Coxeter groups with the geometric realisation (see \cite{EW:SC}). Therefore at least in theory the techniques used in this paper could be applied an arbitrary finite Coxeter group. Note that this will be much harder since the difficult three-coloured Zamolodzhikov relations will appear and there is no complete diagrammatic classification of the idempotents.
\end{rem}

\bibliography{biblio}
\bibliographystyle{amsalpha}
\end{document}